\documentclass[12pt,reqno]{amsart}
\usepackage[backend=biber,style=alphabetic,giveninits=true,maxnames=4]{biblatex}
\bibliography{OurReferences}
\usepackage{amssymb}
\usepackage{subfiles}
\usepackage{mathtools}
\usepackage[utf8]{inputenc}
\usepackage{float}
\usepackage{color}
\usepackage{tikz}
\usepackage[textsize=footnotesize]{todonotes}
\definecolor{darkgreen}{rgb}{0,0.45,0}
\usepackage{graphicx}
\usepackage{adjustbox}
\usepackage[all,2cell]{xy}\UseAllTwocells\SilentMatrices
\usepackage{bm}
\usepackage[colorlinks,citecolor=darkgreen,linkcolor=darkgreen,urlcolor=gray,final,hyperindex,linktoc=page,pdfpagelabels]{hyperref}
\usepackage[capitalize]{cleveref}
\crefname{equation}{}{}
\usepackage{enumerate}

\newlength{\leftstackrelawd}
\newlength{\leftstackrelbwd}
\def\leftstackrel#1#2{\settowidth{\leftstackrelawd}%
{${{}^{#1}}$}\settowidth{\leftstackrelbwd}{$#2$}%
\addtolength{\leftstackrelawd}{-\leftstackrelbwd}%
\leavevmode\ifthenelse{\lengthtest{\leftstackrelawd>0pt}}%
{\kern-.5\leftstackrelawd}{}\mathrel{\mathop{#2}\limits^{#1}}}

\AtEveryBibitem{\clearfield{doi}}
\AtEveryBibitem{\clearfield{url}}
\AtEveryBibitem{\clearfield{issn}}
\AtEveryBibitem{\clearfield{isbn}}
\DeclareSourcemap{
  \maps[datatype=bibtex]{
    \map[overwrite=true]{
      \step[fieldsource=fjournal]
      \step[fieldset=journal, origfieldval]
    }
  }
}

\usetikzlibrary{decorations.markings, arrows.meta, calc, fit, quotes, cd, math, external}
\usetikzlibrary{strings}

\tikzset{tick/.style={postaction={decorate,decoration={markings,
mark=at position 0.5 with {\draw[-] (0,.4ex) -- (0,-.4ex);}}}}}
\tikzset{tickx/.style={postaction={decorate,decoration={markings,mark=at position 0.5 with
{\fill circle [radius=.28ex];}}}}}

\newcommand{\spn}[5]{\SelectTips{eu}{10}\xymatrix@C=.25in{#1 & #2\ar[l]|-{#4}\ar[r]|-{#5} & #3}}
\newcommand{\cspn}[5]{\SelectTips{eu}{10}\xymatrix@C=.25in{#1\ar[r]|-{#4} & #2 & #3\ar[l]|-{#5}}}

\newtheorem{proposition}{Proposition}[section]
\newtheorem{lemma}[proposition]{Lemma}
\newtheorem{corollary}[proposition]{Corollary}
\newtheorem{theorem}[proposition]{Theorem}
\newtheorem*{theorem*}{Theorem}

\theoremstyle{definition}
\newtheorem{definition}[proposition]{Definition}
\newtheorem{example}[proposition]{Example}

\newtheorem{remark}[proposition]{Remark}

\newcommand{\mlt}{m} 
\newcommand{\uni}{j} 
\newcommand{\comlt}{d} 
\newcommand{\couni}{\epsilon} 
\newcommand{\lcomlt}{\delta} 
\newcommand{\lcouni}{\varepsilon} 
\newcommand{\lmlt}{\mu}
\newcommand{\luni}{\eta}

\newcommand{\gact}{\tau} 
\newcommand{\lcoact}{\rho} 
\newcommand{\lact}{\nu} 
\newcommand{\gcoact}{\chi} 

\newcommand{\cas}[2]{\stackrel{\scriptscriptstyle#2}{#1}} 

\newcommand{\Comon}{\sf Comon}

\newcommand{\id}{{\sf id}}

\newcommand{\Hom}{{\sf Hom}}
\newcommand{\Ob}{\sf Ob}
\newcommand{\Hopf}{{\sf Hopf}}
\newcommand{\sHopf}{{\sf sHopf}}

\newcommand{\Mod}{{\sf Mod}}
\newcommand{\opMod}{{\sf opMod}}

\newcommand{\Grph}{\sf Grph}
\newcommand{\dGrph}{\sf dGrph}
\newcommand{\Span}{{\sf Span}}

\newcommand{\Mon}{{\sf Mon}}

\newcommand{\VopCat}{\Vv\textrm{-}\sf{opCat}}

\newcommand{\op}{\textrm{op}}
\newcommand{\sop}{{*,\textrm{op}}}
\newcommand{\co}{\textrm{co}}
\newcommand{\cop}{\textrm{cop}}
\newcommand{\coop}{\textrm{coop}}
\newcommand{\Tr}{\textrm{Tr}}

\def\ot{\otimes}

\newcommand{\ca}{\mathcal}

\newcommand{\Cc}{\mathcal{C}}
\newcommand{\Dd}{\mathcal{D}}

\newcommand{\Ii}{\mathcal{I}}

\newcommand{\Vv}{\mathcal{V}}
\newcommand{\Ww}{\mathcal{W}}

\def\text#1{{\rm {\rm #1}}}

\newcommand{\Cat}{{\sf Cat}}

\def\ol{\overline}
\def\ul{\underline}

\newcommand{\Nat}{\mathsf{Nat}}
\newcommand{\Set}{\mathbf{Set}}

\def\lim{{\rm lim\,}}

\definecolor{purple(x11)}{rgb}{0.5, 0.0, 0.5}

\def\Set{{\sf Set}}

\def\ot{\otimes}

\newcommand{\HopfCat}[1]{\Hopf\mhyphen #1 \mhyphen\Cat}
\newcommand{\sHopfCat}[1]{\sHopf\mhyphen #1 \mhyphen\Cat}

\def\*C{{}^*\hspace*{-1pt}{\Cc}}

\def\text#1{{\rm {\rm #1}}}

\def\ol{\overline}
\def\ul{\underline}

\newcommand{\dagarr}[3]{\ar@{<..}@<-2pt>[#1]_-{#3} \ar@<2pt>[#1]^-{#2}}

\mathchardef\mhyphen="2D 

\newcommand{\cd}[2][]{\vcenter{\hbox{\xymatrix#1@!C{#2}}}}

\def\atpd{s}
\def\braid{\sigma}
\def\invbraid{{\braid^{-1}}}

\newcommand{\coev}{{\rm coev}}

\newcommand{\ev}{{\rm ev}}

\usepackage{comment,silence}
\title
{A {L}arson-{S}weedler theorem for {H}opf $\Vv$-categories}

\author[Buckley]{Mitchell Buckley}
\address{Data61, CSIRO, Sydney NSW 1466, Australia}
\email{mitchell.buckley@data61.csiro.au}

\author[Fieremans]{Timmy Fieremans}
\address{Vakgroep Wiskunde, Faculteit Ingenieurswetenschappen, Vrije Universiteit Brussel, Pleinlaan 2, B-1050 Brussel, Belgium}
\email{tfierema@vub.ac.be}

\author[Vasilakopoulou]{Christina Vasilakopoulou} 
\address{Departement of Mathematics, University of Patras, 26504 Greece}
\email{cvasilak@math.upatras.gr}

\author[Vercruysse]{Joost Vercruysse}
\address{D\'epartement de Math\'ematiques, Facult\'e des sciences, Universit\'e Libre de Bruxelles, Boulevard du
Triomphe, B-1050 Bruxelles, Belgium}
\email{jvercruy@ulb.ac.be}

\begin{document}

\begin{abstract}
The aim of this paper is to extend the classical Larson-Sweedler theorem, namely that a $k$-bialgebra has a non-singular integral (and in particular is Frobenius) if and only if it is a finite dimensional Hopf algebra, to the `many-object' setting of Hopf categories. To this end, we provide new characterizations of Frobenius $\Vv$-categories and we develop the integral theory for Hopf $\Vv$-categories. Our results apply to Hopf algebras in any braided monoidal category as a special case, and also relate to Turaev's Hopf group algebras and particular cases of weak and multiplier Hopf algebras.
\end{abstract}

\maketitle

\setcounter{tocdepth}{1}
\tableofcontents
\section{Introduction}\label{sec:introduction}

The classical Larson-Sweedler theorem \cite{LS} characterizes
finite dimensional Hopf algebras among finite dimensional bialgebras as those that
possess a non-singular left integral. The existence of such an integral
implies in particular that the Hopf algebra is {\em Frobenius}, and this result has been refined in  \cite{P1971}.
Just like Hopf algebras, Frobenius algebras have both
an algebra and a coalgebra structure, although with different compatibility conditions.
In fact, if some finite dimensional algebra $A$ has a Hopf structure, the induced Frobenius
structure on $A$ has (in general) a different comultiplication and counit than
the ones from the Hopf algebra structure.
One of the reasons why the Larson-Sweedler theorem is so important is that it led to definitions of (locally compact) quantum groups by means of well-behaving integrals rather than antipodes. 
The result of Larson and Sweedler furthermore inspired many other results about the connection between Hopf and Frobenius structures on a given (bi)algebra. For example, in Hopf-Galois theory, a Frobenius structure on the Hopf algebra allows to describe the (Morita) equivalence between modules over the ring of coinvariants and Hopf modules in a more direct and symmetric way, see \cite{CFM1990} and
\cite{CVW2004}. More recently, some new connections between Frobenius and Hopf properties have been observed in \cite{S2019}.

Just as the notion of Hopf algebra has been generalized in several ways, so has the Larson-Sweedler theorem.
For example, Larson-Sweedler type theorems for weak Hopf algebras \cite{V:LSweak}, multiplier Hopf algebras
\cite{VDW:LS} and weak multiplier Hopf algebras \cite{KVD2018} have been formulated; however, as explained in \cite{iovanov}, the weak
Hopf algebras case is quite subtle, since the `difficult' direction
of the Larson-Sweedler theorem only holds under
additional assumptions on the {\em target algebra} of the weak Hopf algebra.

The aim of the current paper is to unify and generalize these results in proving a Larson-Sweedler theorem for {\em Hopf $\Vv$-categories}. Hopf 
categories were introduced in \cite{BCV} and can be understood as a many-object generalization of usual Hopf algebras, in the same way as one can 
understand a groupoid as a many-object version of a group. More precisely, a {\em semi-Hopf $\Vv$-category}, where $\Vv$ is a braided monoidal 
category, is a category $A$ that is enriched over the monoidal category of comonoids in $\Vv$: namely, if we denote the hom-object between two objects 
$x,y$ of $A$ by $A_{x,y}$, it comes equipped with usual composition $\mlt_{xyz}\colon A_{x,y}\ot A_{y,z}\to A_{x,z}$ and identities $\uni_{x}\colon 
I\to A_{x,x}$ but also comultiplications $\lcomlt_{xy}\colon A_{x,y}\to A_{x,y}\ot A_{x,y}$ and counits $\lcouni_{xy}\colon A_{x,y}\to I$. A semi-Hopf 
category is called {\em Hopf} if it admits an antipode given by $\atpd_{xy}:A_{x,y}\to A_{y,x}$ for any two objects $x,y$ in $A$. With appropriate 
axioms, a (semi-) Hopf category with one object is exactly a Hopf monoid (bimonoid) in $\Vv$.
In \cite{Paper1a}, we showed that such Hopf categories have a natural interpretation as {\em oplax Hopf algebras}; in \cite{Bohm2017} an alternative interpretation of Hopf categories was given as Hopf monads in a suitable monoidal bicategory.
Hopf categories have the interesting feature that they are general enough to cover many interesting examples of generalized Hopf-structures, and concrete enough to manipulate them without the need to involve heavy higher categorical machinery. In particular, by `packing' a Hopf category, one obtains interesting examples of weak (multiplier) Hopf algebras, whose target algebra is a direct product of copies of the base ring.

Also in \cite{Paper1a}, we introduced the notion of a Frobenius $\Vv$-category as the natural Frobenius analogue of Hopf categories, and which again 
serves as a many-object generalization of Frobenius monoids. In contrast to the one-object case, where both Hopf algebras and Frobenius algebras 
consist of algebras that also have a coalgebra structure albeit with different compatibility conditions, the coalgebraic structure of a Hopf category 
and a Frobenius category are of a very different nature. Indeed the coalgebraic structure of a Hopf category is `local' in the sense that every 
hom-object $A_{x,y}$ is a comonoid (in the monoidal category $\Vv$). On the other hand, a Frobenius $\Vv$-category is at the same time a 
$\Vv$-enriched category and a $\Vv^{\op}$-enriched category, which means that it comes with cocomposition $A_{x,y}\to A_{x,z}\ot A_{z,y}$ and 
coidentity arrows $A_{x,x}\to I$ for all objects $x,y,z$ of $A$, i.e. the coalgebraic structure of a Frobenius $\Vv$-category is `global'. 
In fact, this difference is exactly one of the advantages of working with this many-object generalization. For example, as said before, when the same 
algebra has both a Frobenius and a Hopf structure, then the two coalgebraic structures are not the same; in the many-object setting, having the 
same comultiplication would not even make sense definition-wise.

In the main result of our paper \cref{oppositeLarsonSweedler} we characterize locally rigid (i.e. such that all Hom-objects $A_{x,y}$ are rigid 
objects in $\Vv$) Frobenius Hopf categories as those semi-Hopf categories that possess non-singular integrals or equivalently as those Hopf categories 
for which the integral spaces are isomorphic to the monoidal unit. Such a Frobenius Hopf category then naturally possesses four different structures 
(category, opcategory, local monoid and local comonoid) that can be combined in different ways to naturally form Hopf and Frobenius structures, see 
\cref{hopffrobconnection}.

As we argue in this paper, this result properly generalizes the classical Larson-Sweedler theorem to the many-object setting by taking 
$\Vv=\mathsf{Vect}_k$. On the other hand, by taking the one-object version of our general theorem, we also obtain a version of the Larson-Sweedler 
theorem for Hopf monoids in braided monoidal categories.
This is in fact a folklore result,
however no explicit written reference can be found in the literature. Moreover, the one object-version of the above mentioned 
\cref{hopffrobconnection} clarifies many results related to `interacting Hopf Frobenius' structures, see e.g.\ \cite{Zanasi2017,CD2019}.

Finally, let us point out that the term ``Hopf category'' has been used in literature before, albeit with different meaning. For example, Crane and Frenkel \cite{Crane} considered Hopf (monoidal) categories to characterize 4-dimensional extended quantum field theories. In a different direction, Turaev \cite{T2010} introduced notions of (crossed) Hopf $G$-categories and (crossed) Hopf $G$-categories (where $G$ is a group) to characterize homotopy quantum field theories. The exact relation between these various structures is not clear at the moment and is the subject of forthcoming work, although a first connection between our work and Turaev's Hopf $G$-categories is discussed in \cref{sec:applications}.

\subsection*{Outline}
In \cref{sec:HopfVCats}, we review some basic properties of Hopf $\Vv$-categories. For example, we study in detail how invertibility of the antipode morphisms is related to the existence of an {\em op-antipode}. We also show how the notion of Hopf categories is closely related to the notion of bi-Galois objects and explain how this leads to the construction of non-trivial examples of Hopf categories (see \cref{Galoisobjectsremark}). After recalling the fundamental theorem for Hopf modules, we also prove the fundamental theorem of Hopf opmodules over Hopf categories (see \cref{fundhopfmod}). 

In \cref{sec:FrobeniusVCats}, we provide equivalent characterizations of the Frobenius $\Vv$-categories from \cite{Paper1a} in terms of self-duality, Casimir elements, trace maps, module isomorphisms and pairs of Frobenius (i.e. two-sided adjoint) functors. These all very naturally generalize the classical ones in the many-object setting, however in non-trivial ways.

The main results of this paper can be found in \cref{sec:mainresults}. After briefly recalling the classical setting, we present a detailed study of 
the integral theory for Hopf categories. As one can expect, this theory becomes much more involved from the 1-object case, since the integral space is 
no longer described as an equalizer but as a more general limit. We investigate the relation between the existence of integrals and Frobenius 
structures on a Hopf category. In particular, we show that a Frobenius Hopf $\Vv$-category also has a local Frobenius structure, i.e. all hom-objects 
$H_{x,y}$ are Frobenius algebras in $\Vv$. Furthermore, this additional local algebra structure is isomorphic to the local algebra structure of the 
dual opcategory $H^\sop$; these four structures, the local and global algebra and coalgebra structures on a single Frobenius and Hopf category fit 
together as explained in \cref{hopffrobconnection}.
We then prove our main result: a generalization of the Larson-Sweedler theorem for Hopf $\Vv$-categories, \cref{oppositeLarsonSweedler2}. We also show that in the particular case of $k$-linear Hopf categories, where $k$ is a commutative base ring for which all projective modules are free, our theorem reduces to a result that subsumes the classical Larson-Sweedler theorem (\cref{classicalLScat}).

In the final \cref{sec:applications}, we present some applications of our result. In particular, in the one-object case we recover the classical Larson Sweedler theorem for Hopf algebras, but also for several of their generalizations, such as monoidal Hom-Hopf algebras and graded Hopf algebras. Other applications to Turaev's Hopf group coalgebras \cite{T2010}, weak (multiplier) Hopf algebras and groupoid algebras are presented as well.

\subsection*{Acknowledgements}

JV wants to thank Paolo Saracco for interesting and motivating discussions on the interaction between Hopf and Frobenius properties. He also thanks the FNRS for support through the MIS grant "Antipode". This work was initiated when both MB and CV were working as postdoctoral researchers at the {\em Universit\'e Libre de Bruxelles} within the framework of the ARC grant ``Hopf algebras and the Symmetries of Non-commutative Spaces'' funded by the ``F\'ed\'eration Wallonie-Bruxelles''. CV was supported by the General Secretariat for Research and Technology (GSRT) and the Hellenic Foundation for Research and Innovation (HFRI).
All authors thank the referees for their careful reading and useful comments.

\section{Hopf \texorpdfstring{$\ca{V}$}{V}-categories}\label{sec:HopfVCats}

In this section we recall some basic notions and constructions relatively to the concept of a Hopf $\Vv$-category, where $\Vv$ is a braided monoidal category; relevant references to that end are \cite{BCV} and \cite{Paper1a}. We assume familiarity with the basics of theory of monoidal
categories, see \cite{BraidedTensorCats}, as well as the theory of (co)monoids, Hopf monoids and Frobenius monoids.

\subsection{Preliminary results}\label{Hopfcats}

 In what follows, $(\ca{V},\otimes,I)$ denotes a monoidal category which, by Mac-Lane's coherence theorem, we will regard as a strict monoidal category without loss of generality.
 
A standard reference for the theory of enriched categories is \cite{Kelly}.
Briefly recall that 
a $\Vv$-\emph{enriched graph} is a family of objects $\{A_{x,y}\}_{x,y\in X}$ in $\Vv$, indexed by its set of objects $X$;
we shall use that notation for hom-objects, rather than the more common $A(x,y)$.
Along with $\Vv$-\emph{graph morphisms}, i.e.\ functions between the sets of objects with
arrows $F_{xy}\colon A_{x,y}\to B_{fx,fy}$ in $\Vv$, enriched graphs form a
category $\Vv$-$\Grph$. It has a subcategory $\Vv$-$\Cat$ of $\Vv$-enriched graphs equipped with composition laws
$\mlt_{xyz}\colon A_{x,y}\otimes A_{y,z}\to A_{x,z}$ (again notice the difference with standard terminology) and identities $\uni_x\colon I\to A_{x,x}$
satisfying the usual associativity and unity conditions. A $\Vv$-functor is then a $\Vv$-graph morphism that respects this structure. If $F\colon\Vv\to\Ww$ is a monoidal functor, it induces a \emph{change of base} functor $\Vv$-$\Cat\to\Ww$-$\Cat$.

We call $k$-\emph{linear} categories those enriched in the category $\Mod_k$ of $k$-modules for a commutative ring $k$.  In what follows, for a $k$-linear category $A$ we usually write composition as simple concatenation, namely $pq\coloneqq\mlt_{xyz}(p\otimes q)$. We also write $1_{x,x}$ for $\uni_x(1)$, the image of $1\in k$ under the identity map $\uni_x\colon k\to A_{x,x}$. Finally for $r\in k$ and $m\in M$, we denote by $r\cdot m$ the scalar multiplication, which gives the natural isomorphism $k\ot M\cong M$.

If $\Vv$ is equipped with a braiding $\braid$, every $\Vv$-graph $A$ has an \emph{opposite} $\Vv$-graph $A^\op$ with the same objects and
hom-objects $A^\op_{x,y}\coloneqq A_{y,x}$. In case $A$ is moreover a $\Vv$-category, $A^\op$ is a also a $\Vv$-category whose composition is 
$$\mlt^{A^\op}_{xyz}\colon A_{y,x}\ot A_{z,y}\xrightarrow{\braid^{\textrm{-}1}}A_{z,y}\ot A_{y,x}\xrightarrow{\mlt_{zyx}}A_{z,x}$$

If $(A,\mlt,\uni)$ is a $\Vv$-category, a (right) $A$-\emph{module} \cite{Lawvereclosedcats} $(N,\gact)$ is a $\Vv$-graph
$\{N_{x,y}\}$ over the same set of objects, equipped with actions
$\gact_{xyz}\colon N_{x,y}\otimes A_{y,z}\to N_{x,z}$
satisfying
\begin{equation*}\label{VMod}
\begin{tikzcd}[row sep=.3in,column sep=.5in]
N_{x,y}\ot A_{y,z}\ot A_{z,w}\ar[r,"{\gact_{xyz}\ot1}"]\ar[d,"1\ot\mlt_{yzw}"'] &
N_{x,z}\ot A_{z,w}\ar[d,"\gact_{xzw}"] \\
N_{x,y}\ot A_{y,w}\ar[r,"\gact_{xyw}"] & N_{x,w}
\end{tikzcd}\quad \text{and} \quad
\begin{tikzcd}[row sep=.3in]
N_{x,y}\ar[r,"1\ot\uni_y"]\ar[dr,"1"'] & N_{x,y}\ot A_{y,y}\ar[d,"\gact_{xyy}"] \\
& N_{x,y}
\end{tikzcd}
\end{equation*}
Morphisms are identity-on-objects $\Vv$-graph maps $\{\varphi_{xy}\colon N_{x,y}\to P_{x,y}\}$
such that
\begin{equation}\label{VModmaps}
\begin{tikzcd}
N_{x,y}\ot A_{y,z}\ar[r,"\gact_{xyz}"]\ar[d,"\varphi_{xy}\ot1"'] & N_{x,z}\ar[d,"\varphi_{xz}"] \\
P_{x,y}\ot A_{y,z}\ar[r,"\gact_{xyz}"'] & P_{x,z}
\end{tikzcd}
\end{equation}
These form a category of (right) $A$-modules which we denote $\Vv$-$\Mod_A$.
Clearly, any $\Vv$-category is both a left and right $A$-module $(A,\mlt)$ called the \emph{regular} $A$-module. 

Notice that in the definition of a right $A$-module for example, the left indexing object of the graph is not playing any role. As a result, for any 
right $A$-module $(N,\gact)$ and any map $h:X\to X$ we can define the {\em $h$-shuffle} of $N$ as the (right) $A$-module $(N^h,\gact^h)$ 
where $N^h_{x,y}:=N_{h(x),y}$ and $\gact^h_{xyz}=\gact_{h(x)yz}$.

Finally, recall~\cite[\S 9]{Monoidalbicatshopfalgebroids} that a $\ca{V}$-\emph{opcategory}
$C$ is a category enriched in the opposite monoidal category $\ca{V}^\op$.
Explicitly, and for future reference, there exist cocomposition and coidentity arrows in $\Vv$
\begin{equation}\label{Vopcat}
\comlt_{xyz}\colon C_{x,z}\to C_{x,y}\otimes C_{y,z},\quad \couni_{x}\colon C_{x,x}\to I
\end{equation}
satisfying coassociativity and counity axioms:
\begin{equation}\label{Vopcat1}
\begin{tikzcd}[column sep=.6in, row sep=.5in]
C_{x,w}\ar[r,"\comlt_{xyw}"]\ar[d,"\comlt_{xzw}"'] & C_{x,y}\otimes C_{y,w}\ar[d,"1\otimes\comlt_{yzw}"] \\
C_{x,z}\otimes C_{z,w}\ar[r,"\comlt_{xyz}\otimes1"'] & C_{x,y}\otimes C_{y,z}\otimes C_{z,w}
\end{tikzcd}\;
\begin{tikzcd}[column sep=.4in, row sep=.5in]
C_{x,y}\otimes C_{y,y}\ar[d,"1\otimes\couni_{y}"'] & C_{x,y}
\ar[l,"\comlt_{xyy}"']\ar[r,"\comlt_{xxy}"]\ar[dl,"\sim"]\ar[dr,"\sim"'] &
C_{x,x}\otimes C_{x,y}\ar[d,"\couni_{x}\otimes1"] \\
C_{x,y}\otimes I && I\otimes C_{x,y}
\end{tikzcd}
\end{equation}
where the coherence isomorphisms in $\Vv$ are suppressed. Similarly, a $\Vv$-\emph{opfunctor} is a $\Vv^\op$-functor.
Together these form a category $\VopCat$. In a dual way to modules, there is a category of $C$-\emph{opmodules} $\Vv$-$\mathsf{opMod}_C$ equipped 
with a $C$-coaction which is compatible with cocomposition and coidentities.

An object $A$ in a monoidal category $\Vv$ has a \emph{left dual} $A^*$ when there exists evaluation and coevaluation morphisms $\ev\colon A^*\ot A\to I$ and $\coev\colon I\to A\ot A^*$ making the following diagrams commute, where the associator and unitors are suppressed:
\begin{equation}\label{eq:dualcond}
 \begin{tikzcd}
A\ar[r,"\coev\ot1"]\ar[dr,"\id"'] & A\ot A^*\ot A\ar[d,"1\ot\ev"] && A^*\ar[r,"1\ot\coev"]\ar[dr,"\id"'] & A^*\ot A\ot A^*\ar[d,"\ev\ot1"] \\
& A &&& A^*
\end{tikzcd}
\end{equation}
Since all duals of an object $A$ in $\ca{V}$ are naturally isomorphic to one another, we will from now on speak about `the' dual of $A$.
Dually,  $A$ is called the \emph{right dual} of $A^*$ and when the monoidal category is braided, there is a bijection between left and right duals 
induced by composing the evaluation and coevaluation by the braiding and its inverse, accordingly. A monoidal category with duals is called 
\emph{rigid} or \emph{autonomous}. Each morphism $f\colon A\to B$ gives rise to $f^*\colon B^*\to A^*$ via
$$
B^*\xrightarrow{1\ot\coev}B^*\ot A\ot A^*\xrightarrow{1\ot f\ot1}B^*\ot B\ot A^*\xrightarrow{\ev\ot 1}A^*.
$$

\begin{remark}\label{rem:classicalduals}
It is a standard fact that if $\Vv$ has left/right duals, then it is 
left/right monoidal closed via $-\otimes B\dashv [B,-]_\ell\cong -\otimes B^*$ and $A\otimes -\dashv[A,-]_r\cong A^{\diamond}\otimes-$, where 
$A^{\diamond}$ denotes a right dual of $A$. Notice that the opposite is not true: in a (left) monoidal closed category $\ca{V}$, we can indeed denote 
$A^*=[A,I]$ which e.g. for $\Vv=\mathbf{Vect}_k$ gives the classical `dual' of a (possible infinite dimensional) vector space, that also comes with 
an evaluation map $[A,I]\otimes A\to I$ from the counit of the tensor-hom adjunction. However, this is not 
necessarily the categorical dual as defined above, unless $A$ is finite dimensional.
\end{remark}

We will henceforth call a $\Vv$-enriched graph or category $A$ \emph{locally rigid} when all hom-objects $A_{xy}$ have duals in $\Vv$ -- but $\Vv$ 
itself is not necessarily rigid.
We will denote by $\cas{\coev}{xy}$ and $\cas{\ev}{xy}$ the corresponding (co)evaluation maps for the hom-object $A_{x,y}\in\Vv$ of a locally rigid $\Vv$-(op)category or simply enriched graph.

\begin{example}\label{daggermodules}
Suppose $G$ is a locally rigid $\Vv$-graph; for example, it is $k$-linear where all hom-objects are finite dimensional for a field $k$, or more generally it is $\Mod_R$-enriched where all hom-objects are finitely generated projective $R$-modules for a commutative ring $R$.
There is the dual $\Vv$-graph $G^*$ with the same objects and hom-objects $G^*_{x,y}$, for example $\Hom_k(G_{x,y},k)$ in the linear case. Notice that `dual' here does not refer to `opposite' as is the usual terminology. 

Morever, there is the opposite dual graph $(G^*)^\op$ henceforth denoted $G^\sop$ given by $(G^\sop)_{x,y}=G^*_{y,x}$. From now on, we will omit the parenthesis and write $G^\sop_{x,y}=G^*_{y,x}$. 

Consider now any {\em left} $A$-module $(M,\gact)$ which is locally rigid as a $\Vv$-graph. Then $M^\sop$ becomes 
in a natural way a {\em right} $A$-module by means of the action
$$M^*_{y,x}\ot A_{y,z} \xrightarrow{1\ot 1\ot\cas{\coev}{zx}}M^*_{y,x}\ot A_{y,z} \ot M_{z,x}\ot M^*_{z,x} \xrightarrow{1\ot \gact_{yzx}\ot 
1}M^*_{y,x}\ot M_{y,x}\ot M^*_{z,x} \xrightarrow{\cas{\ev}{yx}\ot 1}M^*_{z,x}$$
We denote this module by $M^\dagger$.
Similarly, the 
dual of any right $A$-module $N$ is becomes naturally a left $A$-module which we denote by ${}^\dagger N$.

Dually, for a $\Vv$-opcategory $C$, any {\em right} locally rigid $C$-opmodule $(N,\gcoact)$ gives rise to a {\em left} $C$-opmodule 
structure on $N^\sop$ via
\begin{displaymath}
 N^*_{z,x}\xrightarrow{1\ot\cas{\coev}{zy}}N^*_{z,x}\ot N_{z,y}\ot N^*_{z,y}\xrightarrow{1\ot\gcoact_{zxy}\ot1}N^*_{z,x}\ot N_{z,x}\ot 
A_{x,y}\ot N^*_{z,y}\xrightarrow{\cas{\ev}{zx}\ot1\ot1}A_{x,y}\ot N^*_{z,y}
\end{displaymath}
and similarly a left one makes its opposite dual into a right $C$-opmodule.
\end{example}

In \cite[Theorem 5.5]{BCV}, the linear case of the following result is exhibited.

\begin{proposition}\label{A*opcategory}
For a locally rigid $\Vv$-category $A$, its opposite dual $A^\sop$ has the structure of a $\Vv$-opcategory. Dually, $C^\sop$ is a $\Vv$-category for any locally rigid $\Vv$-opcategory $C$.
\end{proposition}

\begin{proof}
Cocomposition and counits are given by applying the functor
$(\textrm{-})^*\colon \Vv^\op\to\Vv$ (restricted to the dualizable objects) to the composition and identities of $A$ as in
\begin{align}
\comlt_{xyz}&\colon A^\sop_{x,z}=A^*_{z,x}\xrightarrow{\;\mlt^*_{zyx}\;}(A_{z,y}\otimes A_{y,x})^*\cong A^*_{y,x} \otimes A^*_{z,y}=A^\sop_{x,y}\ot A^\sop_{y,z} \label{eq:A*opopcategory}\\
\couni_{x}&\colon A^*_{x,x}\xrightarrow{\;\uni_x^*\;}I \nonumber
\end{align}
 Essentially, the strong anti-monoidal functor $(-)^*$ via $\phi_{XY}\colon(X\ot Y)^*\cong Y^*\ot X^*$ turns a $\Vv$-category $A$ into a $\Vv$-opcategory $A^\sop$, and dually a $\Vv$-opcategory $C$ into a $\Vv$-category $C^\sop$ via a process similar to the change of enrichment base. 

 \end{proof}
 
Also in \cite[Proposition 5.4]{BCV}, the linear case of the following can be found.
\begin{proposition}\label{prop:staropmodules}
 For any locally rigid $\Vv$-category $A$, $\Vv$-$\Mod_A\cong\Vv$-$\opMod_{A^\sop}$.
\end{proposition}
\begin{proof}
Suppose $(M,\gact)$ is a right $A$-module. Then it can be given the structure of a right $A^\sop$-opmodule via
\begin{displaymath}
 M_{x,z}\xrightarrow{1\ot\cas{\coev}{zy}}M_{x,z}\ot A_{z,y}\ot A^*_{z,y}\xrightarrow{\gact_{xzy}\ot1}M_{x,y}\ot A^*_{z,y} 
\end{displaymath}
and vice versa, if $(N,\gcoact)$ is a right $A^\sop$-opmodule, then there is an $A$-action
\begin{displaymath}
 N_{x,y}\ot A_{y,z}\xrightarrow{\gcoact_{xzy}\ot1}N_{x,z}\ot A^*_{y,z}\ot A_{y,z}\xrightarrow{1\ot\cas{\ev}{yz}}
 N_{x,z}
\end{displaymath}
These two establish a bijection between $A$-modules and $A^\sop$-opmodules; clearly this works for left (op)modules too. Similarly 
$\Vv$-$\mathsf{opMod}_C\cong\Vv$-$\Mod_{C^\sop}$ for a $\Vv$-opcategory $C$.
\end{proof}

Hopf enriched categories, introduced in~\cite{BCV}, constitute a natural many-object generalization of a Hopf monoid in a braided monoidal category. In what follows, suppose that $(\Vv,\ot,I,\braid)$ is a braided monoidal category, and recall that its category of comonoids $\Comon(\Vv)$ inherits the monoidal structure, via $$C\ot D\xrightarrow{\delta\ot\delta}C\ot C\ot D\ot D\xrightarrow{1\ot\braid\ot1}C\ot D\ot C\ot D.$$
Notice that we use Latin letters to denote `global' operations (those that may relate different hom-objects, i.e. of different indices), and Greek letters to denote `local' operations (those that concern each hom-object object individually).

\begin{definition}\label{def:sHopfcat}
A \emph{semi-Hopf $\Vv$-category} $A$
is a $\Comon(\Vv)$-enriched category. Explicitly, it consists of objects together with a collection of $A_{x,y}\in\Vv$
for any two objects $x,y$, and families of morphisms in $\Vv$
\begin{gather*}
\mlt_{xyz}\colon A_{x,y}\otimes A_{y,z}\to A_{x,z}\qquad \uni_x\colon I\to A_{x,x} \\
\lcomlt_{xy}\colon A_{x,y}\to A_{x,y}\otimes A_{x,y}\qquad \lcouni_{xy}\colon A_{x,y}\to I 
\end{gather*}
which make $A$ a $\Vv$-category, each $A_{x,y}$ a comonoid in $\Vv$, and
satisfy
\begin{equation}
\begin{gathered}
\xymatrix@C=.6in@R=.2in{
 A_{x,y}\otimes  A_{y,z} \ar[r]^-{\lcomlt_{xy}\otimes\lcomlt_{yz}} \ar[dd]_{\mlt_{xyz}} &  A_{x,y}\otimes A_{x,y}\otimes
 A_{y,z}\otimes A_{y,z} \ar[d]^{ 1\otimes\braid\otimes  1} \\ 
&  A_{x,y}\otimes A_{y,z}\otimes A_{x,y}\otimes A_{y,z} \ar[d]^{\mlt_{xyz}\otimes\mlt_{xyz}} \\ 
 A_{x,z} \ar[r]_{\lcomlt_{xz}} &  A_{x,z}\otimes A_{x,z} }\label{Hax1} \\
\xymatrix{
I \ar[r]^{\sim} \ar[d]_-{\uni_{x}} & I\otimes I \ar[d]^-{\uni_{x}\otimes\uni_{x}} \\ 
 A_{x,x} \ar[r]_-{\lcomlt_{xx}} &  A_{x,x}\otimes A_{x,x} } \quad
\xymatrix@C=.5in{
 A_{x,y}\otimes A_{y,z} \ar[r]^-{\lcouni_{xy}\otimes\lcouni_{yz}} \ar[d]_-{\mlt_{xyz}} & I\otimes I
\ar[d]^-{\sim} \\ 
 A_{x,z} \ar[r]_-{\lcouni_{xz}} & I }\quad
\xymatrix{
I \ar[r]^{\id} \ar[d]_-{\uni_{x}} & I \ar[d]^{\id} \\ 
 A_{x,x} \ar[r]_-{\lcouni_{xx}} & I }
\end{gathered}
\end{equation}
\end{definition}

Semi-Hopf $\Vv$-categories with $\Comon(\Vv)$-functors
form the category $\Comon(\Vv)\mhyphen\Cat$ which we also denote $\sHopfCat{\Vv}$.

\begin{example}\label{manyobjectbialgebra}
Every bimonoid in a braided monoidal category $\Vv$ is a one-object semi-Hopf $\Vv$-category.
\end{example}

\begin{example}
If $A$ is a semi-Hopf $\Vv$-category for $(\Vv,\otimes,I,\braid)$, it gives rise to new semi-Hopf $\Vv$-categories $A^\op$, $A^\cop$, $A^{\op,\cop}$ and $A^{\cop,\op}$ as follows, see also \cite[\S 3]{BCV}.
\begin{enumerate}
\item $A^\op_{xy}=A_{yx}$ with composition $A_{yx}\ot A_{zy}\xrightarrow{\braid^{\textrm{-}1}}A_{zy}\ot A_{yx}\xrightarrow{\mlt_{zyx}}A_{zx}$; the monoidal base of the enrichment is $(\ca{V},\braid^{-1})$.
\item $A^{\cop}_{xy}=A_{xy}$ with local comultiplications $\lcomlt_{xy}$ post-composed with the inverse braiding;
again the monoidal base is $(\ca{V},\braid^{-1})$.
\item $A^{\op,\cop} = (A^{\op})^{\cop}$ has hom-objects $A^{\op,\cop}_{xy}=A_{yx}$, composition is pre-composed with the inverse braiding, comultiplication is post-composed with the usual braiding, and the monoidal base is $(\ca{V},\braid)$.
\item $A^{\cop,\op} = (A^{\cop})^{\op}$ has hom-objects $A^{\cop,\op}_{xy}=A_{yx}$, composition is pre-composed with the usual braiding, comultiplication is post-composed with the inverse braiding, and the monoidal base is $(\ca{V},\braid)$.
\end{enumerate}
Clearly, if $\Vv$ is symmetric then one no longer needs to distinguish between the braiding and its inverse. 
\end{example}

We now turn to Hopf categories and their basic properties. 

\begin{definition}\cite[Def.~3.3]{BCV}\label{def:Hopfcat}
A \emph{Hopf $\Vv$-category} $H$ is a semi-Hopf $\Vv$-category equip\-ped with a family of maps
$\atpd_{xy} \colon  H_{x,y} \to  H_{y,x}$ satisfying
\begin{equation}\label{HopfCatAntipodeEquations}
\begin{gathered}
\cd[@C-2em]{
 &  H_{x,y}\otimes H_{x,y} \ar[rr]^{ 1 \otimes \atpd_{xy}} & &  H_{x,y}\otimes H_{y,x} \ar[dr]^{\mlt_{xyx}} & \\
  H_{x,y} \ar[rr]^{\lcouni_{xy}} \ar[ur]^{\lcomlt_{xy}} & & I \ar[rr]^{\uni_{x}} & &  H_{x,x} 
} \\ 
\cd[@C-2em]{
 &  H_{x,y}\otimes H_{x,y} \ar[rr]^{\atpd_{xy} \otimes 1} & &  H_{y,x}\otimes H_{x,y} \ar[dr]^{\mlt_{yxy}} & \\
  H_{x,y} \ar[rr]^{\lcouni_{xy}} \ar[ur]^{\lcomlt_{xy}} & & I \ar[rr]^{\uni_{y}} & &  H_{y,y} \rlap{\ .}
} 
\end{gathered}
\end{equation}
This $\Vv$-graph map $\atpd\colon H\to H^\op$ is called the \emph{antipode} of $H$.
\end{definition}
If only the upper (respectively lower) diagram commutes, $\atpd$ is called a \emph{right} (respectively \emph{left}) \emph{antipode} of $H$.

\begin{definition}\label{def:opantipode}
An {\em op-antipode} for a semi-Hopf $\Vv$-category $H$ is an antipode for $H^\op$, i.e.\ a 
family of maps $\ol\atpd_{xy} \colon  H_{y,x} \to  H_{x,y}$ satisfying the following two conditions:
\begin{displaymath}
\begin{gathered}
\cd[@C=3em]{
 &  
 H_{x,y}\otimes H_{x,y} \ar[r]^{ 1 \otimes \ol\atpd_{yx}} & 
 H_{x,y}\otimes H_{y,x} \ar[r]^{ \braid^{-1}} & 
 H_{y,x}\otimes H_{x,y} \ar[dr]^{\mlt_{yxy}} & \\
 H_{x,y} \ar[rr]^{\lcouni_{xy}} \ar[ur]^{\lcomlt_{xy}} & & I \ar[rr]^{\uni_{y}} & &  H_{y,y} 
} \\ 
\cd[@C=3em]{
 &  
 H_{x,y}\otimes H_{x,y} \ar[r]^{ \ol\atpd_{yx} \otimes 1 } & 
 H_{y,x}\otimes H_{x,y} \ar[r]^{ \braid^{-1}} & 
 H_{x,y}\otimes H_{y,x} \ar[dr]^{\mlt_{xyx}} & \\
 H_{x,y} \ar[rr]^{\lcouni_{xy}} \ar[ur]^{\lcomlt_{xy}} & & I \ar[rr]^{\uni_{x}} & &  H_{x,x} 
} 
\end{gathered}
\end{displaymath}
A left (right) op-antipode for $H$ is a left (right) antipode for $H^\op$.
\end{definition}

\begin{remark}\cite{BCV}\label{antipodeproperties}
The following properties of antipodes can be deduced from the definitions
\begin{align*}
\atpd_{xz} \circ \mlt_{xyz} &= \mlt_{zyx} \circ \braid \circ (\atpd_{xy} \otimes \atpd_{yz})& \lcomlt_{yx} \circ \atpd_{xy} &= \braid \circ (\atpd_{xy} \otimes \atpd_{xy}) \circ \lcomlt_{xy} \\
\atpd_{xx} \circ \uni_x & = \uni_x, &\lcouni_{yx} \circ \atpd_{xy} & = \lcouni_{xy} 
\end{align*}
Since op-antipodes are antipodes for $H^\op$, they also obey respective formulas. 
\end{remark}

\begin{lemma}\label{atpdbijopatpd}
If $H$ is a Hopf $\Vv$-category with antipode $\atpd$, then $\ol\atpd$ is an op-antipode
if and only if each $\ol\atpd_{yx}$ is inverse to $\atpd_{xy}$.
\end{lemma}

\begin{proof}
If $\ol\atpd_{xy}$ is as in \cref{def:opantipode}, then on one side we find
\begin{align*}
\ol\atpd_{yx} \circ \atpd_{xy} &= \ol\atpd_{yx} \circ \atpd_{xy} \circ (\lcouni_{xy} \ot A_{x,y}) \circ \lcomlt_{xy}\\
&= (\lcouni_{xy} \ot A_{x,y}) \circ (A_{x,y} \ot\ol\atpd_{yx}) \circ (A_{x,y} \ot\atpd_{xy}) \circ \lcomlt_{xy} \\
&= m_{xxy} \circ (\uni_{xx} \ot A_{x,y}) \circ (\lcouni_{xy} \ot A_{x,y}) \circ (A_{x,y} \ot\ol\atpd_{yx}) \circ (A_{x,y} \ot\atpd_{xy}) \circ \lcomlt_{xy}\\
&= m_{xxy} \circ (m_{xyx} \ot A_{x,y}) \circ (A_{x,y} \ot \atpd_{xy} \ot A_{x,y}) \circ  (A_{x,y} \ot A_{x,y} \ot\ol\atpd_{yx})\\
&~~ \circ (A_{x,y} \ot A_{x,y} \ot\atpd_{xy})\circ (\lcomlt_{xy} \ot A_{x,y}) \circ\lcomlt_{xy}\\
&= m_{xxy} \circ (m_{xyx} \ot A_{x,y}) \circ (A_{x,y} \ot A_{y,x} \ot\ol\atpd_{yx}) \circ (A_{x,y} \ot \atpd_{xy} \ot A_{x,y})\\
&~~ \circ (A_{x,y} \ot A_{x,y} \ot \atpd_{xy})\circ ( A_{x,y}\ot \lcomlt_{xy} ) \circ\lcomlt_{x,y}\\
&= m_{xyy} \circ (A_{x,y} \ot m_{yxy}) \circ (A_{x,y} \ot A_{y,x} \ot\ol\atpd_{yx}) \circ (A_{x,y} \ot \braid^{-1})  \\
&~~ \circ (A_{x,y} \ot \lcomlt_{yx}) \circ (A_{x,y} \ot \atpd_{xy}) \circ\lcomlt_{xy}\\
&= m_{xyy} \circ (A_{x,y} \ot \uni_{y})\circ (A_{x,y} \ot \lcouni_{yx}) \circ (A_{x,y} \ot \atpd_{xy}) \circ\lcomlt_{xy}\\
&= m_{xyy} \circ (A_{x,y} \ot \uni_{y})\circ (A_{x,y} \ot \lcouni_{xy}) \circ\lcomlt_{xy}\\
&= A_{x,y}
\end{align*}
So $\ol\atpd$ is left inverse to $s$ -- recall these are identity-on-objects graph morphisms.
A similar argument shows that $\ol\atpd$ is also right inverse to $s$ and the one direction is established.

Now suppose that an antipode $\atpd$ has inverses $\ol s_{yx}$ for each $s_{xy}$. Then these indeed form an op-antipode; for example, the left axiom is verified by
\begin{align*}
m_{yxy} \circ (\ol\atpd_{xy}\ot A_{x,y}) \circ \sigma^{-1} \circ \comlt_{xy}
& = m_{yxy} \circ (\ol\atpd_{xy}\ot A_{x,y}) \circ \sigma^{-1} \circ \comlt_{xy} \circ \atpd_{yx} \circ \ol{\atpd}_{xy}\\
& = m_{yxy} \circ (\ol\atpd_{xy}\ot A_{x,y}) \circ \sigma^{-1} \circ \sigma \circ (\atpd_{yx} \ot \atpd_{yx})\\
&~~ \circ \comlt_{yx} \circ \ol{\atpd}_{xy}\\
& = m_{yxy} \circ (A_{x,y}\ot \atpd_{yx}) \circ \comlt_{yx} \circ \ol{\atpd}_{xy}\\
& = \uni_{yy} \circ \couni_{xy} \circ \ol{\atpd}_{xy}\\
& = \uni_{yy} \circ \couni_{yx}\\
\end{align*}
where $\ol\atpd$ satisfies conditions dual to those in \cref{antipodeproperties}
merely by being inverse to $\atpd$. That $\ol\atpd$ a right antipode is proved dually.
\end{proof}

If $H$ and $K$ are Hopf $\Vv$-categories, a $\Comon(\Vv)$-functor $F\colon H\to K$
is called a \emph{Hopf $\Vv$-functor} if $\atpd_{fxfy}\circ F_{xy}
=F_{yx}\circ\atpd_{xy}$ for all $x,y\in X$. It is shown in~\cite[Prop.~3.10]{BCV} that any
$\Comon(\Vv)$-functor between Hopf $\Vv$-categories
automatically satisfies that condition; hence we have a full subcategory
$\HopfCat{\Vv}$ of $\Comon(\Vv)\mhyphen\Cat$.

\begin{example}\label{manyobjectHopfalgebra}
Every Hopf algebra $H$ in a braided monoidal $\Vv$ is a one-object Hopf $\Vv$-category;
this fulfils its purpose as a many-object generalization. In particular, each endo-hom object $H_{x,x}$ of an arbitrary Hopf $\Vv$-category $H$ is a Hopf monoid in $\Vv$.
\end{example}

\begin{remark}
 It was shown in \cite{Paper1a} that $\sHopfCat{\Vv}$ and $\HopfCat{\Vv}$ are in fact categories of \emph{oplax bimonoids} and \emph{Hopf oplax bimonoids} in a symmetric monoidal bicategory $\Span|\Vv$. This exhibits a more elaborate sense in which Hopf structure can be generalized in higher categorical settings, and Hopf categories are example of such. 
\end{remark}

\begin{example}
\cite[3.12]{BCV} The `linearization' functor $L \colon \Set \to \Mod_k$
which sends each set to the free $k$-module on that set, is a strong monoidal functor.
Hence it induces a change-of-base functor between $\HopfCat{\Set}$ and $\HopfCat{\Mod_k}$, namely ordinary Hopf categories which are the same as groupoids, and $k$-linear Hopf categories. As a result, every groupoid $G$ determines a $k$-linear
Hopf category $H$ with $H_{x,y}\coloneqq LG_{x,y}$, the free $k$-module on the set of morphisms $x\to y$ in $G$.
\end{example}

\begin{proposition}\label{packed}
Suppose that $H$ is a Hopf $\Vv$-category with finitely many objects, where $\Vv$ has a zero object and biproducts, that are also preserved by the tensor. The \emph{packed form} of $H$,
$$\hat H= \coprod\limits_{x,y} H_{x,y}$$
is a weak Hopf algebra.
\end{proposition}

\begin{proof}
Suppose that $\Vv$ is a monoidal category
with coproducts that commute with the tensor product --
as is the case for any monoidal closed category -- and a zero object.
For an arbitrary $\Vv$-graph $G$, we get a new graph $\hat{G}\ot\hat{G}$ given by
\begin{displaymath}
\hat{G}\ot \hat{G}=\coprod_{x,y}G_{x,y}\ot\coprod_{z,u}G_{z,u}\cong\coprod_{x,y,z,u} G_{x,y}\ot G_{z,u}\rlap{\ .}
\end{displaymath}
Now if $(A,\mlt,\uni)$ is a $\Vv$-category, first of all we can define families of maps
\begin{displaymath}
A_{x,y}\ot A_{z,u}\xrightarrow{\mathrm{mlt}_{xyzu}}A_{x,u} =
\begin{cases}
\mlt_{xyu}, & \textrm{if } y=z \\
0, & \text{else}
\end{cases}\textrm{ and } 
I\xrightarrow{\mathrm{uni}_{xy}}A_{x,y} =
\begin{cases}
\uni_{x}, & \textrm{if }x=y \\
0, & \textrm{else}
\end{cases}.
\end{displaymath}
The first one induces, for every $x,y,z,u$, a composite diagonal map as below ---
where the vertical arrows are the canonical injections --- hence the universal property of coproducts yields a (unique) map
$\mu\colon\hat{A}\ot\hat{A}\to\hat{A}$
\begin{displaymath}
\begin{tikzcd}[column sep=.7in]
\coprod\limits_{x,y,z,u} A_{x,y}\ot A_{z,u}
\ar[r,dashed,"\mu"] & \coprod\limits_{x,u}A_{x,u} \\
A_{x,y}\ot A_{z,u}\ar[u,hook]\ar[r,"\mathrm{mlt}_{xyzu}"']\ar[ur,dotted] &
A_{x,u}\ar[u,hook]
\end{tikzcd}
\end{displaymath}
which is easy to check that is associative.
If moreover the set of objects $X$ is finite and (finite) biproducts exist in $\Vv$, so 
$\coprod\limits_{x,y}A_{x,y}=\prod\limits_{x,y}A_{x,y}$,
then we also obtain a (unique) map 
\begin{displaymath}
\begin{tikzcd}[column sep=.7in]
& \prod\limits_{x,y}A_{x,y}\ar[d] \\
I\ar[r,"\mathrm{uni}_{xy}"']\ar[ur,dashed,"\eta"] & A_{x,y}
\end{tikzcd}
\end{displaymath}
which satisfies unity conditions. Therefore under these conditions for a $\Vv$-category $A$, $\hat{A}$ naturally obtains a monoid structure in $\Vv$.

Now suppose that $(H,\mlt,\uni,\lcomlt,\lcouni,\atpd)$ is a Hopf $\Vv$-category. Then $\hat H$ is also a
comonoid in $\Vv$, since comonoids are closed under colimits in any monoidal category
(see e.g.\ \cite{CDV}). Explicitly, the comultiplication and counit again follow
from the universal property of coproducts, induced by
\begin{displaymath}
H_{x,y}\xrightarrow{\lcomlt_{xy}}H_{x,y}\ot H_{x,y}\hookrightarrow
\coprod_{x,y} H_{x,y} \ot \coprod_{x,y} H_{x,y},\qquad H_{x,y}\xrightarrow{\lcouni_{xy}}I.
\end{displaymath}
It was shown in \cite{BCV} that with this structure, the packed form of a Hopf $\Mod_k$-category
with a finite set of objects is a weak Hopf algebra \cite{weakHopf}. 
Under the above conditions, this can also be proved for general Hopf $\Vv$-categories.
Explicitly, the compatibility between multiplication and comultiplication can be shown to hold due to the top of \cref{Hax1}, whereas
the other two axioms $\varepsilon\circ\mu\circ(\mu\ot1)=(\varepsilon\ot\varepsilon)\circ(\mu\ot\mu)\circ(1\ot\delta\ot1)=(\varepsilon\ot\varepsilon)\circ(\mu\ot\mu)\circ(1\ot\braid\ot1)\circ(1\ot\delta\ot1)$ and $(\delta\ot1)\circ\delta\circ\eta=(1\ot\mu\ot1)\circ(\delta\ot\delta)\circ(\eta\ot\eta)=$ can also be verified.
\end{proof}
In particular, applying the above proposition to
the previous example, one obtains the usual groupoid algebra $kG$ from $H$,
as a packed form: $kG = \bigoplus\limits_{\substack{x,y \in G}} H_{x,y}$. 

\begin{example}[\emph{Hopf opcategories}]\label{Hopfopcats}
If we replace $\Vv$ with $\Vv^\op$ at \cref{def:Hopfcat,def:sHopfcat},
we obtain the notion of a \emph{(semi) Hopf $\Vv$-opcategory}, called
\emph{dual Hopf category} in \cite{BCV}. Since $\Comon(\Vv^\op)\cong\Mon(\Vv)^\op$,
a semi-Hopf $\Vv$-opcategory $(C,\comlt,\couni,\lmlt,\luni)$
is precisely a $\Mon(\Vv)$-opcategory, i.e.\ it is equipped with cocomposition
and counit morphisms $\comlt_{xyz},\couni_{x}$ as in \cref{Vopcat}, together with local multiplication
and unit morphisms $\lmlt_{xy}\colon C_{x,y}\ot C_{x,y}\to C_{x,y}$, $\luni_{xy}\colon I\to C_{x,y}$
making each hom-object a monoid in $\Vv$, subject to compatibility conditions.
Moreover, a Hopf $\Vv$-opcategory comes with arrows $\atpd_{xy}\colon C_{y,x}\to C_{x,y}$ satisfying dual axioms
to \cref{HopfCatAntipodeEquations}:
\begin{displaymath}
\begin{tikzcd}
C_{x,y}\ot C_{y,x}\ar[rr,"1\ot\atpd_{xy}"] && C_{x,y}\ot C_{x,y}\ar[d,"\lmlt_{xy}"]  \\
C_{x,x}\ar[u,"\comlt_{xyx}"]\ar[r,"\couni_x"] & I\ar[r,"\luni_{xy}"] & C_{x,y}
\end{tikzcd}\quad
\begin{tikzcd}
C_{x,y}\ot C_{y,x}\ar[rr,"\atpd_{yx}\ot1"] && C_{y,x}\ot C_{y,x}\ar[d,"\lmlt_{yx}"]  \\
C_{x,x}\ar[u,"\comlt_{xyx}"]\ar[r,"\couni_x"] & I\ar[r,"\luni_{yx}"] & C_{y,x}
\end{tikzcd}
\end{displaymath}
\end{example}

\begin{proposition}\label{AstarHopfversions}
If $A$ is a (semi-)Hopf locally rigid $\Vv$-category, $A^\sop$ naturally obtains the structure of a (semi-)Hopf $\Vv$-opcategory.
\end{proposition}

\begin{proof}
By \cref{A*opcategory}, any $\Vv$-category $A$ gives rise to a
$\Vv$-opcategory $A^\sop$ given by $A^\sop_{x,y}=A^*_{y,x}$. Moreover, the local
comonoid structure turns into a local monoid structure under the strong antimonoidal $(-)^*\colon\Vv^\op\to\Vv$.
More explicitly, if $(A,\mlt,\uni,\lcomlt,\lcouni)$ is the semi-Hopf $\Vv$-category,
$(A^\sop,\phi\circ\mlt^*,\uni^*,\lcomlt^*\circ\phi,\lcouni^*)$ is an induced semi-Hopf $\Vv$-opcategory structure on
$A^\sop$ where the cocomposition and counit are given as in \cref{eq:A*opopcategory} and
multiplication and unit are
\begin{gather*}
\lmlt_{xy}\colon A^\sop_{x,y}\ot A^\sop_{x,y}=A^*_{y,x}\ot A^*_{y,x}\xrightarrow{\stackrel{\phi}{\cong}}(A_{y,x}\ot
A_{y,x})^*\xrightarrow{\lcomlt^*_{yx}}A^*_{y,x}=A^\sop_{x,y} \\
\luni_{xy}\colon \xrightarrow{\lcouni^*_{y,x}}A^*_{y,x}=A^\sop_{x,y}
\end{gather*}
Dually, if $C$ is a semi-Hopf $\Vv$-opcategory, its opposite dual $C^\sop$ is a semi-Hopf $\Vv$-category.
\end{proof}

\begin{remark}
In the $k$-linear setting, the above formulas take the following form, where $\sum{\cas{e}{xy}_i \ot \cas{f}{xy}_i} \in A_{x,y}\ot A^*_{x,y}$ denotes a finite dual basis of any $A_{x,y}$. The global comultiplication is given by
$$f\mapsto \sum_{i,j} f(\cas{e}{zy}_i\cas{e}{yx}_j) \cas{f}{yx}_j\ot \cas{f}{zy}_i$$ and the local multiplication is given by $f\ot g\mapsto f*g$ where $(f*g)(a)=g(a_1)f(a_2)$ for all $a\in A_{x,y}$.

Notice that in the one-object case of a $k$-bialgebra $A$, the induced $A^\sop$ as described above in fact coincides with the classical $(A^*)^{\op,\coop}$, namely the opposite-coopposite of the classical dual $k$-bialgebra $A^*$.
\end{remark}

The proposition below is the generalization of the classic `fusion map' formulation in this many-object setting.
\begin{proposition}\label{Galoisobjects}
Let $H$ be a Hopf $\Vv$-category. Then for any two objects $x,y$ in $H$, we have that
the above canonical map
\[\xymatrix{
H_{x,x}\ot H_{x,y} \ar[rr]^-{1\ot \delta_{xy}} && H_{x,x}\ot H_{x,y}\ot H_{x,y} \ar[rr]^-{m_{xxy}\ot 1} && H_{x,y}\ot H_{x,y}
}\]
is an isomorphism.
\end{proposition}

\begin{proof}
One can easily check that an inverse of the canonical map is given by
\[
H_{x,y}\ot H_{x,y}\xrightarrow{1\ot\delta_{xy}}H_{x,y}\ot H_{x,y}\ot H_{x,y} \xrightarrow{1\ot s_{xy}\ot 1}
H_{x,y}\ot H_{y,x}\ot H_{x,y} \xrightarrow{m_{xyx}\ot 1}H_{x,x}\ot H_{x,y}
\]
\end{proof}

Of course, the previous proposition also can be applied to Hopf opcategories.
Since it is well-known that, when working over a base field, i.e. $\Vv={\sf Vect_k}$, the bijectivity of the canonical map implies that the space of (co)invariants is trivial, we can deduce from the previous proposition the following interesting result.
\begin{theorem}
For a Hopf ${\sf Vect}_k$-category H, $H_{x,y}$ is an $H_{x,x}$-$H_{y,y}$ bi-Galois co-object for any pair of objects $x,y$ in $H$; for a Hopf ${\sf Vect}_k$-opcategory H, $H_{x,y}$ is an $H_{x,x}$-$H_{y,y}$ bi-Galois object for any pair of objects $x,y$ in $H$.
\end{theorem}
In particular, we find that each non-zero $H_{x,y}$ is isomorphic as a $k$-vector space to both $H_{x,x}$ and $H_{y,y}$. Moreover, $H_{x,x}$ and $H_{y,y}$ are isomorphic as $k$-vector spaces if $H_{x,y}$ or $H_{y,x}$ is non-zero.
This observation leads to some interesting examples of Hopf (op)categories as below. 
\begin{example}\label{Galoisobjectsremark}
Let $H$ be a Hopf algebra, and $A$ any (faithfully flat, right) Galois object of $H$. Then we know, see \cite{HopfbiGalois}, that one can construct a Hopf algebra $L=(A\ot A)^{coH}$ such that $A$ becomes an $L$-$H$ bi-Galois object. Hence we obtain a Hopf opcategory with two objects $x,y$ by putting $H_{x,x}=L$, $H_{y,y}=H$, $H_{x,y}=A$ and $H_{y,x}=A^{op}$. 
\end{example}

\subsection{The fundamental theorem of Hopf categories}

In this section, we recall the fundamental theorem for Hopf modules of Hopf $\Vv$-categories; details and proofs can be found in \cite{BCV}.
For $(A,\mlt,\uni,\lcomlt,\lcouni)$ a semi-Hopf $\Vv$-category, a right \emph{Hopf module} is a $\Vv$-graph $M$ over the same set of objects,
with a global $A$-action and a local $A$-coaction
\begin{displaymath}
\gact_{xyz}\colon M_{x,y}\ot A_{y,z}\to M_{x,z}, \quad \lcoact_{xy}\colon M_{x,y}\to M_{x,y}\ot A_{x,y} 
\end{displaymath}
making $M$ into an enriched $A$-module (\cref{Hopfcats}) and each $M_{x,y}$ into
an ordinary $A_{x,y}$-comodule, and furthermore satisfy
\begin{displaymath}
\begin{tikzcd}[column sep=.6in]
M_{x,y}\ot A_{y,z}\ar[r,"\lcoact_{xy}\ot\lcomlt_{yz}"]\ar[d,"\gact_{xyz}"']
& M_{x,y}\ot A_{x,y}\ot A_{y,z}\ot A_{y,z}\ar[r,"1\ot\braid\ot1"]
& M_{x,y}\ot A_{y,z}\ot A_{x,y}\ot A_{y,z}\ar[d,"\gact_{xyz}\ot\mlt_{xyz}"] \\
M_{x,z}\ar[rr,"\lcoact_{xz}"] && M_{x,z}\ot A_{x,z}
\end{tikzcd}
\end{displaymath}
There is a category with objects Hopf $A$-modules, and morphisms $\Vv$-graph maps that respect the global $A$-action and local $A$-coactions; it is denoted by $\Vv\mhyphen\Mod_A^A$. 

Dually, for a semi-Hopf $\Vv$-opcategory $(C,\comlt,\couni,\lmlt,\luni)$ as in \cref{Hopfopcats}, a right \emph{Hopf opmodule} is a $\Vv$-graph $N$ equipped with a global $C$-coaction and local $C$-action 
\begin{displaymath}
\gcoact_{xyz}\colon N_{x,z}\to N_{x,y}\ot C_{y,z}\quad \quad
\lact_{xy}\colon N_{x,y}\ot C_{x,y}\to N_{x,y}
\end{displaymath}
making $N$ into an enriched $C$-opmodule and each $N_{x,y}$ into an ordinary $C_{x,y}$-module, compatible in that 
$\gcoact_{xyz}\circ\lact_{xz}=(\lact_{xy}\ot \lmlt_{yz})\circ (1\ot\braid\ot 1)\circ (\gcoact_{xyz} \ot \comlt_{xyz})$.
The category of right Hopf opmodules over $C$ is denoted as $\Vv\mhyphen\opMod_C^C$.

\begin{example}\label{hopfmodulestructuredual}
Suppose that $(H,\mlt,\uni,\lcomlt,\lcouni,\atpd)$ is a locally rigid Hopf $\Vv$-category; recall by \cref{AstarHopfversions} that $H^\sop$ is a Hopf $\Vv$-opcategory.

\begin{enumerate}
\item\label{it:H1opmodule} $H$ is a right Hopf $H^\sop$-opmodule via the following coaction and action:
\begin{equation}\label{AisarightAstarop-opmodule} 
\gcoact_{xyz}\colon H_{x,z} \xrightarrow{1\cas{\coev}{zy}}H_{x,z}\ot H_{z,y}\ot H^*_{z,y} \xrightarrow{\mlt_{xzy}1}H_{x,y}\ot H^*_{z,y} 
\end{equation}
\begin{displaymath}
\lact_{xy}\colon H_{x,y}\ot H^*_{y,x}\xrightarrow{\braid}H^*_{yx}\ot H_{x,y} \xrightarrow{1\lcomlt_{xy}}H^*_{y,x}\ot H_{x,y}\ot H_{x,y}
\xrightarrow{1\atpd_{xy}1}H^*_{y,x}\ot H_{y,x}\ot H_{x,y} \xrightarrow{\cas{\ev}{yx}1}H_{x,y} 
\end{displaymath}
With this structure, H is called a \emph{type 1} Hopf $H^\sop$-opmodule, denoted $H_1$. 
\item\label{it:H2opmodule} $H^\op$ is a right $H^\sop$-opmodule, via the following action and coaction: 
\begin{displaymath}
\begin{tikzcd}[column sep=.25in]
\gcoact_{xyz}\colon H_{z,x}\ar[drrr,dashed] \ar[r,"1\cas{\coev}{zy}"] & H_{z,x}\ot H_{z,y}\ot H^*_{z,y} \ar[r,"1\atpd_{zy}1"] & H_{z,x}\ot H_{y,z}\ot H^*_{z,y} \ar[r,"\braid1"] & H_{y,z}\ot H_{z,x}\ot H^*_{z,y}\ar[d,"\mlt_{yzx}1"] \\
&&& H_{y,x}\ot H^*_{z,y}
\end{tikzcd}
\end{displaymath}
\begin{displaymath}
\lact_{xy}\colon H_{y,x}\ot H^*_{y,x}\xrightarrow{\lcomlt_{yx}1}H_{y,x}\ot H_{y,x}\ot H^*_{y,x} \xrightarrow{1\braid}H_{y,x}\ot H^*_{y,x}\ot H_{y,x} \xrightarrow{1\cas{\ev}{yx}}H_{y,x} 
\end{displaymath}
With this structure, $H^\op$ is called a \emph{type 2} Hopf $H^\sop$-opmodule, denoted $H_2$.

\item\label{it:H*1module} $H^*$ a right Hopf $H$-module 
via the action and coaction
\begin{equation}\label{AstarisarightA-module}
\begin{tikzcd}[column sep=.3in,row sep=.15in]
\gact_{xyz} \colon H^*_{x,y}\ot H_{y,z}\ar[ddrr,dashed]\ar[r,"1\atpd_{yz}\cas{\coev}{xz}"] & H^*_{x,y}\ot H_{z,y}\ot H_{x,z}\ot H^*_{x,z}\ar[r,"1\braid1"] 
& H^*_{x,y}\ot H_{x,z}\ot H_{z,y}\ot H^*_{x,z}\ar[d,"1\mlt_{xzy}1"] \\ 
&& H^*_{x,y}\ot H_{x,y}\ot H^*_{x,z}\ar[d,"{\cas{\ev}{xy}1}"] \\
&& H^*_{x,z} \\
\lcoact_{xy} \colon H^*_{x,y}\ar[dddrr,dashed] \ar[r, "1\cas{\coev}{xy}"] & 
H^*_{x,y}\ot H_{x,y}\ot H^*_{x,y} \ar[r, "1\lcomlt_{xy}1"]
& H^*_{x,y}\ot H_{x,y}\ot H_{x,y}\ot H^*_{x,y} \ar[d, "\braid11"]\\
&& H_{x,y}\ot H^*_{x,y}\ot H_{x,y}\ot H^*_{x,y} \ar[d, "1\cas{\ev}{xy}1"] \\
&& H_{x,y}\ot H^*_{x,y} \ar[d, "\braid"] \\ 
&& H^*_{x,y}\ot H_{x,y}
\end{tikzcd}
\end{equation}
With this structure, $H^*$ is called a \emph{type 1} Hopf $H$-module, denoted $H^*_1$.
\item\label{it:Astar2} $H^\sop$ is a right Hopf $H$-module via the following action and coaction 
\begin{displaymath}
\begin{tikzcd}[row sep=.15in,column sep=.35in]
\gact_{xyz}\colon H^*_{y,x}\ot H_{y,z}\ar[rrd,dashed] \ar[r, "11\cas{\coev}{zx}"] & H^*_{y,x}\ot H_{y,z}\ot H_{z,x}\ot H^*_{z,x} \ar[r, 
"1\mlt_{yzx}1"]
& H^*_{y,x}\ot H_{y,x} \ot H^*_{z,x} \ar[d, "\cas{\ev}{yx}1"] \\
&& H^*_{z,x} \\
\lcoact_{xy}\colon H^*_{y,x}\ar[rrdd,dashed] \ar[r, "1\cas{\coev}{yx}"] & 
H^*_{y,x}\ot H_{y,x}\ot H^*_{y,x} \ar[r, "1\lcomlt_{yx}1"] & 
H^*_{y,x}\ot H_{y,x}\ot H_{y,x}\ot H^*_{y,x} \ar[d, "\cas{\ev}{yx}\braid"]\\
& & H^*_{y,x}\ot H_{y,x} \ar[d, "1\atpd_{yx}"]\\
& & H^*_{y,x}\ot H_{x,y}
\end{tikzcd}
\end{displaymath}
With this structure, $H^\sop$ is called a \emph{type 2} Hopf $H$-module, denoted $H^*_2$. Notice that $\gact$ is precisely the induced action on the 
regular module from \cref{daggermodules}.
\end{enumerate}
\end{example}

\begin{example}\label{Csopopmodule}
Dually, if $(C,\comlt,\couni,\lmlt,\luni,\atpd)$ is a Hopf $\Vv$-opcategory then the following action and coaction make
$C^\sop$ a right Hopf $C$-opmodule.
\begin{displaymath}
\begin{tikzcd}[row sep=.15in, column sep=.35in]
\lact_{xy} \colon C^*_{y,x}\ot C_{x,y}\ar[ddrr,dashed]\ar[r,"1\atpd_{xy}\cas{\coev}{yx}"] &
C^*_{y,x}\ot C_{y,x}\ot C_{y,x}\ot C^*_{y,x}\ar[r,"1\braid1"] & C^*_{y,x}\ot C_{y,x}\ot C_{y,x}\ot C^*_{y,x}\ar[d,"1\lmlt_{yx}1"]\\
 && C^*_{y,x}\ot C_{y,x}\ot C^*_{y,x}\ar[d,"{\cas{\ev}{yx}1}"] \\
 && C^*_{y,x} \\
\gcoact_{xyz} \colon C^*_{z,x}\ar[r,"\cas{\coev}{yx}1"]\ar[ddrr,dashed] & C_{y,x}\ot C^*_ {y,x}\ot  C^*_{z,x} \ar[r,"\braid1"] 
& C^*_ {y,x}\ot C_{y,x}\ot  C^*_{z,x}\ar[d,"1\comlt_{yzx}1"] \\
&& C^*_ {y,x}\ot C_{y,z}\ot C_{z,x}\ot C^*_{z,x}\ar[d,"11\cas{\ev}{zx}"] \\ 
&& C^*_{y,x} \ot C_{y,z}
\end{tikzcd}
\end{displaymath}
Specifically $\gcoact$ is induced from the regular opmodule structure on $C$ as in \cref{daggermodules}.
\end{example}

We now recall the fundamental theorem for Hopf modules, \cite[Theorem 10.2]{BCV}. In its formulation, we denote by $\Vv$-$\dGrph$ the category of \emph{diagonal} $\Vv$-graphs, namely given by single-indexed families $(M_x)_{x\in X}$ of objects in $\Vv$. Notice that any $\Vv$-graph gives rise to a diagonal one, by considering only its endo-hom objects $M_{x,x}$.

\begin{theorem}\label{fundhopfmod}
Let $(A,\mlt,\uni,\lmlt,\lcouni)$ be a semi-Hopf $\Vv$-category and suppose that $\Vv$ has equalizers. The functor
\begin{equation}\label{eq:functorfundhopfmod}
 -\ot A\colon\Vv\mhyphen\dGrph \to \Vv\mhyphen\Mod_A^A
\end{equation}
that maps some $\{N_x\}_{x\in X}$ to $\{N_x\ot A_{x,y}\}_{x,y\in X}$ with $A$-action $1\ot\mlt_{xyz}$ and coaction $1\ot\lcomlt_{xy}$, has a right adjoint $(-)^{\co A}$ as in
$$
\begin{tikzcd}[column sep=1in]
\Vv\mhyphen\dGrph\ar[r,bend left=20,"-\ot A"]\ar[r,phantom,"\bot"description] &  \Vv\mhyphen\Mod_A^A\ar[l,bend left=20,dashed,"(-)^{\co A}"]
\end{tikzcd}
$$
defined on a Hopf $A$-module $(M,\tau,\rho)$ by the equalizer
\begin{equation}\label{eq:coinvariants}
\xymatrix{M^{\co A}_x\ar[rr]^{i_{x}} && M_{x,x} \ar@<.5ex>[rr]^-{\rho^M_{x,x}} \ar@<-.5ex>[rr]_-{1\ot \uni_{x}} && M_{x,x}\ot A_{x,x}}
\end{equation}
Moreover, $A$ is a Hopf $\Vv$-category if and only if the above functors establish an equivalence of categories; in particular, $M^{\co A}\ot A\cong M$ for any Hopf $A$-module $M$.
\end{theorem}

\begin{proof}(sketch)
The unit and counit are given respectively by $\alpha_{x}: N_{x} \to (N_{x} \ot A_{x,x})^{coA}_{x}$
such that $i_{x} \circ \alpha_{x} = N_{x} \ot \uni_{x}$
and $\beta_{xy} = \gact_{xxy} \circ (i_{x} \ot A_{x,y}): M^{coA}_{x} \ot A_{x,y} \to M_{x,y}$.
\end{proof}

Notice that $M^{\co A}_x$ is the \emph{space of coinvariants} for the local Hopf algebra $A_{x,x}$ in the $k$-linear case, see \cref{manyobjectHopfalgebra}, which can in that way be defined in any monoidal category $\Vv$ with equalizers.

The above theorem can also be deduced from viewing a Hopf category
as a special instance of a Hopf comonad on a naturally Frobenius map-monoidale~\cite{Bohm2017}, 
using the fundamental theorem of Hopf modules in that general setting~\cite{BohmLack}. 

Finally, the following result will be of use in later sections.
\begin{lemma}\label{antipodeisinvertible} 
If a Hopf $\Vv$-category $H$ is locally rigid, then its antipode is invertible.
\end{lemma}
\begin{proof}
The proof is essentially the same as \cite[Theorem 4.1]{Takeuchi1999}
and relies on the fundamental theorem for Hopf
modules for Hopf $\Vv$-categories. If $(H,\mlt,\uni,\lcomlt,\lcouni,\atpd)$ is the Hopf $\Vv$-category, we can apply the equivalence of \cref{fundhopfmod} to the right Hopf $H$-module $H^*_{1}$ described in
\cref{hopfmodulestructuredual}(\ref{it:H*1module}) to get a Hopf $H$-module isomorphism
$$\beta_{xy}: (H^*_{1})^{coH}_{x} \ot H_{x,y}
\to (H^{*}_{1})_{x,y}.$$ 
We now consider the following commutative diagram: since the braiding $\braid$ is invertible, and $\lcouni^*$ is split by $\uni^*$, the entire counter-clockwise composite that excludes $s_{x,y}$ constitutes a left inverse to the antipode.
\begin{displaymath}
\begin{tikzcd}
H_{xy} \ar[r, shift right=1.5, "\lcouni^*_{xx} 1"'] \ar[d, "\atpd_{xy}"'] & 
(H^*_{1})_{xx}  H_{xy} \ar[r, "\beta^{-1}_{xx} 1"] \ar[d, "1  \atpd_{xy}"] \ar[l, "\uni^*_{xx}\ 1"', shift right=1.5]& 
(H^*_{1})^{\co H}_{x}  H_{xx}  H_{xy} \ar[r, "\braid 1"] \ar[d, "1  1  \atpd_{xy}"] & 
H_{xx}  (H^*_{1})^{\co H}_{x}  H_{xy} \ar[dr, bend left =20, "1  \beta_{xy}"] \ar[d, "1  1  \atpd_{xy}"'] & \\
H_{yx} \ar[r, "\couni^*_{xx} 1"'] & 
(H^*_{1})_{xx}  H_{yx} \ar[r, "\beta^{-1}_{xx} 1"'] & 
(H^*_{1})^{\co H}_{x}  H_{xx}  H_{yx} \ar[r, "\braid 1"'] & 
H_{xx}  (H^*_{1})^{\co H}_{x}  H_{yx} \ar[r, "1 \zeta_{xy}"'] & 
H_{xx}  (H^*_{1})_{xy} 
\end{tikzcd}
\end{displaymath}
The middle diagrams commute by naturality and the right-most triangle commutes by definition of $\beta_{xy}$, where $\zeta_{xy} := (\cas{\ev}{xx} \ot 
1) \circ (1\ot\mlt_{xyx}\ot1) \circ (1\ot\braid\ot1) \circ (1\ot 1\ot \cas{\coev}{yx}) \circ (i_{x} \ot 1)$ is defined precisely in order to cancel 
$\beta$ introduced in \cref{fundhopfmod} for the appropriate action of $H_1^*$, \cref{AstarisarightA-module}.

Now a dual argument shows that $\atpd^*$ has a left inverse, hence $\atpd$ also has a right inverse because taking duals is a contravariant functor; therefore the antipode is invertible.
\end{proof}

\subsection{The fundamental theorem of Hopf opcategories}

\cref{fundhopfmod} can be appropriately dualized to produce a fundamental theorem for Hopf opmodules. However, due to some non-trivial subtle differences between the two cases, in this section we explicitly describe the basic constructions and proofs. In what follows, we fix $(C,\comlt,\couni,\lmlt,\luni)$ to be a semi-Hopf $\Vv$-opcategory as in \cref{Hopfopcats}, for $\Vv$ a braided monoidal category.

In order to specify a functor similarly to \cref{eq:functorfundhopfmod}, notice that for any diagonal $\Vv$-graph $\{N_x\}_{x\in X}$, the families $(N\ot C)_{x,y}:=N_x\ot C_{x,y}$ give a Hopf $C$-opmodule with $C$-action
$1\ot \lmlt_{xy}$ and coaction $1\ot \comlt_{xzy}$. This naturaly defines a functor
$$-\ot C:\Vv\mhyphen\dGrph\to \Vv\mhyphen\opMod^C_C\rlap{\ .}$$
On the other hand, for any Hopf $C$-opmodule $(M,\lact,\gcoact)$, we define the $x$-\emph{coinvariant space} $M^{\co C}_x$ of $M$ to be the limit in the following diagram
\begin{equation}\label{defconvariants}
\cd[@!C@C-3.5em]{
 & & & M^{\co C}_{x} \ar@{.>}[dl]|-{v_{xy}} \ar@{.>}[dr]|-{v_{xz}} \ar@{.>}[dlll]_-{v_{xw}} \ar@{.>}[drrr]^-{v_{xu}} & & & \\
M_{x,w} \ar[dr]|-{1\otimes \eta_{wy}} & & M_{x,y} \ar[dl]|-{\gcoact_{xwy}} \ar[dr]|-{1\ot \eta_{yz}} & & M_{x,z}
\ar[dl]|-{\gcoact_{xyz}} \ar[dr]|-{1\ot \eta_{zu}} & & M_{x,u} \ar[dl]|-{\gcoact_{xzu}} \\ 
\dots & M_{x,w}\otimes C_{w,y} & & M_{x,y}\otimes C_{y,z} & & M_{x,z}\otimes C_{z,u} & \dots 
}
\end{equation}
Explicitly, the object $M^{\co C}_x$ in $\Vv$ comes with maps $v_{xy}:M^{\co C}_x\to M_{xy}$ such that
$\gcoact_{xyz}\circ v_{xz} = (id\ot \luni_{yz})\circ v_{xy}$ for all $y,z$, and is universal with this property.
These spaces form a diagonal $\Vv$-graph $M^{\co C}=\{M^{\co C}_x\}_{x\in X}$ for any $x$, and this is set to be the mapping on objects of a functor $(-)^{\co C}:\Vv\mhyphen\opMod^C_C\to \Vv\mhyphen\dGrph$.

\begin{proposition}\label{fundopcat1}
For any complete category $\Vv$, there is an adjunction
$$\begin{tikzcd}\Vv\mhyphen\dGrph\ar[r,bend left=10,"-\ot C"]\ar[r,phantom,"\scriptstyle\bot"] & \Vv\mhyphen\opMod^C_C.\ar[l,bend left=10,"\;\;(-)^{\co C}"]\end{tikzcd}$$
\end{proposition}

\begin{proof}
For any $\{N_x\}_{x\in X}$ in $\Vv\mhyphen\dGrph$, the maps $1\ot\luni_{xy}:N_x\to N_x\ot C_{x,y}$ induce morphisms $\alpha_x\colon N_x\to (N\ot C)_x^{co H}$ by the universal property of the limit.
On the other hand, for any Hopf $C$-opmodule $M$ let $\beta_{xy}=\lact_{xy}\circ (v_{xy}\ot1)\colon M_x^{\co C}\ot C_{x,y}\to M_{x,y}$.
We can check that $\alpha$ and $\beta$ constitute a unit and counit for the proposed adjunction.
\end{proof}

\begin{theorem}\label{fundopcat_equiv}
Let $\Vv$ be a complete category.
If a semi-Hopf $\Vv$-opcategory $C$ is Hopf, the adjunction of \cref{fundopcat1} is an equivalence of categories. In particular, for any Hopf $C$-opmodule $M$,
\begin{displaymath}
M^{\co C}\otimes C \cong M
\end{displaymath}
\end{theorem}

\begin{proof}
It suffices to show that when $C$ is a Hopf $\Vv$-category, the adjunction $\textrm{-}\ot C\dashv (-)^{coC}$ is an adjoint equivalence, namely the unit and counit are isomorphisms.

An inverse for each $\alpha_x$ as defined in the previous proof is given by
$$\Gamma_{x} := \xymatrix{(N\ot C)^{\co C}_x \ar[rr]^-{v_{xx}}  &&N_x\ot C_{x,x} \ar[rr]^-{1\ot\couni_{xx}} && N_x}$$
It is clear that $\Gamma_{x} \circ \alpha_{x} = \id$ because of the commutativity of the following diagram:
\begin{equation*}
\begin{tikzcd}
N_{x} \ar[r,"\alpha_{x}"] \ar[rr,bend right=15,"1\ot \luni_{xx}"description]\ar[rrr,bend right=20,"1"description] & (N \ot C)^{\co C}_{x}\ar[r,"v_{xx}"]& N_{x}\ot C_{x,x} \ar[r,"1 \ot \couni_{xx}"] & N_{x}
\end{tikzcd}
\end{equation*}
For the other side composite, first note that there is only one endomorphism $f_{x}: (N\ot C)^{\co C}_{x} \to (N\ot C)^{\co C}_{x}$ such that $v_{xx} \circ f_{x} = v_{xx}$ by the universal property of limits; hence this is the identity. Moreover, $v_{xx} \circ (\alpha_{x} \circ \Gamma_{x}) = v_{xx}$ due to
\begin{equation*}
\begin{tikzcd}
(N\ot C)^{\co C}_{x} \ar[r,"v_{xx}"] \ar[drrr, bend right=10,"v_{xx}"']& N_{x}\ot C_{x,x} \ar[r,"1 \ot \couni_{xx}"] \ar[drr, bend right=10,"1"']& N_{x} \ar[r,"\alpha_{x}"]\ar[dr,"1\ot \luni_{xx}"']& (N\ot C)^{\co C}_{x}\ar[d,"v_{xx}"]\\
&&& N_{x} \ot C_{x,x}
\end{tikzcd}
\end{equation*}
so also $\alpha_{x} \circ \Gamma_{x} = \id$.

For each component $\beta_{xy}$ of the counit, an inverse is given by
$$\gamma_{xy} := \xymatrix{
M_{x,y} \ar[rr]^-{\lcoact_{xxy}} && M_{x,x}\ot C_{x,y} \ar[rr]^-{t_x\ot 1} && M_{x}^{\co C}\ot C_{x,y}
}$$
where $t_{x}$ is induced by the universal property of coinvariants and the family of maps $\nu_{xy} \circ 1 \ot \atpd_{yx}) \circ \gcoact_{xyx}$ which form a cone over the required diagram: indeed,
\begin{equation*}
\chi_{xyw} \circ \nu_{xw} \circ (M_{x,w} \ot \atpd_{wx}) \circ \chi_{xwx} = (M_{x,y} \ot \luni_{yw}) \circ \nu_{xy} \circ (M_{x,y} \ot \atpd_{xy}) \circ \chi_{xyx}.
\end{equation*} 
We can verify that this $\gamma_{xy}$ is a one-sided inverse of $\beta_{xy}$ by the following commutative square
\begin{equation*}
\begin{tikzcd}
M_{x,y}\ar[dd,bend right=60,"1"']\ar[r,"\gcoact_{xxy}"]\ar[d,"\gcoact_{xyy}"] & M_{x,x}\ot C_{x,y}\ar[d,"\gcoact_{xyx} \ot1"] \ar[r,"u_x \ot1"]\ar[ddr,phantom,"(**)"]& M^{\co C}_{x} \ot C_{x,y}\ar[dd,"v_{xy} \ot 1"]\\
M_{x,y} \ot C_{y,y}\ar[r,"1\ot \comlt_{yxy}"]\ar[d,"1\ot \couni_{yy}"]\ar[rd,phantom,"(*)"] & M_{x,y}\ot C_{y,x} \ot C_{x,y}\ar[d,"1\ot \atpd_{yx} \ot1"] & \\
M_{x,y}\ar[drr,bend right=30,"1"']\ar[dr,"1\ot \luni_{xy}"] & M_{x,y} \ot C_{x,y} \ot C_{x,y}\ar[r,"\nu_{xy} \ot1"] \ar[d,"1\ot \lmlt_{xy}"] & M_{x,y} \ot C_{x,y}\ar[d,"\nu_{xy}"]\\
 & M_{x,y} \ot C_{x,y} \ar[r,"\nu_{xy}"] & M_{x,y}
\end{tikzcd}
\end{equation*}
for any $x,y,z,w \in X$. The left and bottom triangle, the left upper square and the right lower square commute since $M$ is Hopf $C$-opmodule. The inner diagram $(*)$ follows from the Hopf opcategory axioms and $(**)$ from the universal property of the limit defining coinvariants. Hence
$\beta_{xy} \circ \gamma_{xy} = \id$.

For $\gamma_{xy} \circ \beta_{xy} = \id$, first note that 
\begin{equation}
\label{special1}
\chi_{xwy} \circ \nu_{xy} \circ (v_{xy} \ot C_{x,y}) = (\nu_{xw} \ot C_{w,y}) \circ (M_{x,w} \ot \comlt_{xwy}) \circ (v_{xw} \ot C_{x,y})
\end{equation}
and also it can be shown that 
\begin{equation}
\label{special2}
v_{xy} \circ t_{x} \circ \nu_{xx} \circ (v_{xx} \ot C_{x,x}) = v_{xy} \circ (M^{coC}_{x} \ot \couni_{xx})
\end{equation}
Since $((M^{coC}_{x}\ot C_{x,x})_{x}, (v_{xy} \circ (M^{coC}_{x}\ot \couni_{xx})_{xy}))$ is trivially a cone over the diagram \cref{defconvariants}, by the universal property there exists a unique morphism $h\colon M^{coC}\ot C \to M^{coC} \in \Vv\textrm{-}\dGrph$ such that $v_{xy} \circ h_{x} = v_{xy} \circ (M^{coC}_{x}\ot \couni_{xx})$ for every $x,y \in X$. By \cref{special2}, we know that $h_{x} = t_{x} \circ \nu_{xx} \circ (v_{xx} \ot C_{x,x})$. Since $(M^{coC}_{x} \ot \couni_{xx})$ satisfies this condition trivially, we can deduce by uniqueness of $h$ that they have to be equal:
\begin{equation}
\label{special3}
(M^{coC}_{x} \ot \couni_{xx}) = t_{x} \circ \nu_{xx} \circ (v_{xx} \ot C_{x,x})
\end{equation}
Finally, using the above data, we can compute
\begin{align*}
\gamma_{xy} \circ \beta_{xy} &= (t_{x} \ot C_{x,y}) \circ \chi_{xxy} \circ \nu_{xy} \circ (v_{xy} \ot C_{x,y})\\
&\overset{\cref{special1}}{=} (t_{x} \ot C_{x,y}) \circ (\nu_{xx} \ot C_{x,y}) \circ (M_{x,x} \ot \comlt_{xxy}) \circ (v_{xx} \ot C_{x,y})\\
&\overset{\cref{special3}}{=}  (M^{coC}_{x} \ot \couni_{xx} \ot C_{x,y}) \circ (M^{coC}_{x} \ot \comlt_{xxy})\\
&= M^{coC}_{x} \ot C_{x,y}
\end{align*}
where the last equality is due to $C$ being a $\Vv$-opcategory, hence the proof is complete.
\end{proof}

A `full' fundamental theorem for Hopf opmodules would include the converse of \cref{fundopcat1}; this may be readily proved by adapting the proof of 
the fundamental theorem of Hopf modules given in \cite{BCV}.  We omit it here since it is not required for our purposes.

\section{Frobenius \texorpdfstring{$\ca{V}$}{V}-categories}
\label{sec:FrobeniusVCats}

In \cite[\S 7]{Paper1a}, we introduced Frobenius $\Vv$-categories as Frobenius monoids inside the same monoidal bicategory where Hopf $\Vv$-categories arise as \emph{oplax Hopf monoids}. In this section, we provide characterizations of Frobenius categories in terms of Casimir elements, dual module structures, trace maps and Frobenius functors. These characterizations, necessary for our central results in \cref{sec:mainresults}, naturally generalise those for usual Frobenius algebras~\cite{Caenepeel2002} and are similar to those of Frobenius monads~\cite[Thm 1.6]{Frobeniusmonads}.

\begin{definition}\cite[{}7.1.1]{Paper1a}\label{def:Frobeniuscat}
A \emph{Frobenius $\Vv$-category} $A$ is a $\Vv$-category that is also a $\Vv$-opcategory, namely for every $x,y\in \Ob A$ there is an object $ A_{x,y}\in\Vv$ and maps
\begin{gather*}
\mlt_{xyz}\colon A_{x,y}\otimes  A_{y,z}\to A_{x,z}\qquad \uni_x\colon I\to A_{x,x} \\
\comlt_{xyz}\colon A_{x,z}\to A_{x,y}\otimes A_{y,z}\qquad \couni_x\colon A_{x,x}\to I 
\end{gather*}
satisfying the usual (co)associativity and (co)unitality axioms, and moreover the following diagrams commute:
\begin{equation}\label{frob1}
\begin{tikzcd}[column sep=.3in,row sep=.2in]
 A_{x,y}\ot A_{y,z}\ar[rr,"\comlt_{xwy}\ot1"]\ar[dd,"1\ot\comlt_{ywz}"']\ar[dr,"\mlt_{xyz}"description] &&
 A_{x,w}\ot A_{w,y}\ot A_{y,z}\ar[dd,"1\ot\mlt_{wyz}"] \\
&  A_{x,z}\ar[dr,"\comlt_{xwz}"description] & \\
 A_{x,y}\ot A_{y,w}\ot A_{w,z}\ar[rr,"\mlt_{xyw}\ot1"'] &&
 A_{x,w}\ot A_{w,z}
\end{tikzcd}
\end{equation}
\end{definition}
Along with \emph{Frobenius} $\Vv$-\emph{functors}, namely $\Vv$-graph morphisms which are functors and opfunctors, they form a category $\mathsf{Frob}$-$\Vv$-$\Cat$.

We will call the counit morphisms $\couni_{x}:A_{x,x}\to I$ of a Frobenius $\Vv$-category the {\em trace morphisms}.
Moreover, if $\Vv$ is braided, a Frobenius $\Vv$-category is \emph{symmetric} when
\begin{equation}\label{eq:symmetricFrob}
 \cd[@!C@R-2em]{
  A_{x,y} \ot A_{y,x} \ar[dd]_{\braid} \ar[r]^-{\mlt_{xyx}} & A_{x,x} \ar@/^.3em/[dr]^{\couni_{x}} &  \\
  & & I \\
  A_{y,x} \ot A_{x,y} \ar[r]_-{\mlt_{yxy}} & A_{y,y} \ar@/_.3em/[ur]_{\couni_{y}} &  \\
}
\end{equation}
which translates in the $k$-linear case to $\couni_{x}(ab)=\couni_{y}(ba)$.

A Frobenius monoid in any monoidal category $\Vv$ can be viewed as a one-object Frobenius $\Vv$-category, and in particular every diagonal hom-object $A_{x,x}\in\Vv$ is such. For more examples and discussion of related notions, see~\cite{Paper1a}.

The fact that Frobenius $\Vv$-categories properly generalize Frobenius monoids in $\Vv$ is also exhibited by the following result, which shows that the packed form of a Frobenius $\Vv$-category is a Frobenius monoid provided that the set of objects is finite. Notice that this result cannot be expected to hold in case the set of objects is infinite, since in such a case the packed form cannot be expected to be rigid in $\Vv$, whereas a Frobenius monoid is always rigid.

\begin{proposition}\label{compactification}
Suppose $A$ is a Frobenius $\Vv$-category with a finite object-set $X$, and that 
$\Vv$ has finite biproducts and $\otimes$ preserves them.
The packed form of $A$
$$\hat A=\coprod_{x,y\in X}A_{x,y}$$
is a Frobenius monoid in $\Vv$.
\end{proposition}

\begin{proof}
We already know by \cref{packed} that $(\hat A,\mu,\eta)$ is a monoid in $\Vv$. In a dual way to the multiplication defined therein, the maps
\begin{displaymath}
A_{x,y}\xrightarrow{\mathrm{comlt}_{xyzu}}A_{x,z}\ot A_{u,y}\equiv
\begin{cases}
\comlt_{xzu}, & \textrm{if } z=u \\
0, & \text{else}
\end{cases}\textrm{ and } 
A_{x,y}\xrightarrow{\mathrm{couni}_{xy}}I\equiv
\begin{cases}
\couni_{x}, & \textrm{if }x=y \\
0, & \textrm{else}
\end{cases}
\end{displaymath}
induce (uniquely) comultiplication and counit arrows
$\delta\colon\hat{A}\to\hat{A}\ot\hat{A}$, $e\colon\hat{A}\to I$
via the universal properties of (co)products.
It can then be verified that $\delta$ and $e$ make $\hat A$ into a comonoid, and moreover that $(\hat A,\mu,\eta,\delta,e)$
is a Frobenius monoid in $\Vv$.
\end{proof}

\begin{remark}
Let us note that the following converse of \cref{compactification} holds. For $X$ a finite set, consider the category of $X$-\emph{bigraded} $\Vv$-objects (namely packed forms of $\Vv$-graphs with set of objects $X$) with tensor product
$$(\coprod_{x,y\in X} A_{x,y}) \otimes (\coprod_{x,y\in X} B_{x,y}) = \coprod_{x,y\in X} (\coprod_{u\in X} A_{x,u}\ot A_{u,y})$$
Then a Frobenius monoid in this category is exactly the packed form of a Frobenius $\Vv$-category. 
\end{remark}

We introduce the following notation for any $\Vv$-category $(A,\mlt,\uni)$ which will be useful for the characterization of Frobenius $\Vv$-categories in what follows:
\begin{align}
V_1&=\{\couni=\{\couni_x\}_{x\in X}~|~\couni_x:\ A_{x,x}\to I\}\label{eq:comultiplications}\\
W_1&=\{\comlt=\{\comlt_{xyz}\}_{x,y,z\in X}~|~\comlt_{xyz}:A_{x,z}\to A_{x,y}\ot A_{y,z},\ \textrm{satisfying}\ \eqref{frob1}\}\nonumber
\end{align}
That is, the sets $V_1$ and $W_1$ consist respectively of candidate `trace morphisms' families and `comultiplication' families for a Frobenius structure on $A$.

\subsection{Characterization in terms of Casimir elements}

\begin{definition}\label{def:casimirfamily}
Let $(A,\mlt,\uni)$ be a $\Vv$-category with ${\Ob}A=X$. A \emph{Casimir family}
is a family $e$ of distinguished morphisms $\cas{e}{xy}\colon I\to A_{x,y}\ot A_{y,x}$ indexed by $(x,y)\in X^2$, satisfying the commutativity of the 
following diagram
\begin{equation}\label{casimirdiag}
 \begin{tikzcd}[column sep=.5in]
  A_{x,z}\ar[r,"\cas{e}{xy}\ot1"]\ar[d,"1\ot\cas{e}{zy}"'] & A_{x,y}\ot A_{y,x}\ot A_{x,z}
  \ar[d,"1\ot\mlt_{yxz}"] \\
  A_{x
  ,z}\ot A_{z,y}\ot A_{y,z}\ar[r,"\mlt_{xzy}\ot1"'] & A_{x,y}\ot A_{y,z}
 \end{tikzcd}
\end{equation}
for any triple $(x,y,z) \in X^3$. 
In the $k$-linear case, this gives an $X^2$-indexed family of elements
$\cas{e}{xy}=e^{1}_{x,y} \ot e^{2}_{y,x} \in A_{x,y} \ot A_{y,x}$
such that $ae^{1}_{z,y}\ot e^{2}_{y,z} = e^{1}_{x,y} \ot e^{2}_{y,x} a$
for all $a \in A_{x,z}$.

We denote by $W_2$ the set of all Casimir families for a given $\Vv$-category $A$.
\end{definition}

Using the above notation \cref{eq:comultiplications}, we obtain the following result.

\begin{lemma}\label{Frobcas}
For any $\Vv$-category $A$, we have a bijection between the sets $W_1$ and $W_2$.
\end{lemma}

\begin{proof}
For any family $\comlt=\{\comlt_{xyz}\}_{x,y,z\in X}\in W_1$, define 
morphisms $\cas{e}{xy}\colon I\xrightarrow{\uni_x}A_{x,x}\xrightarrow{\comlt_{xyx}}A_{x,y}\ot A_{y,x}$.
To check that these satisfy the Casimir property \cref{casimirdiag}, we examine the diagram
\begin{displaymath}
\begin{tikzcd}[row sep=.15in,column sep=.6in]
A_{x,z}\ar[dr,"\id"description]\ar[r,"\uni_x\ot1"]\ar[dd,"1\ot\uni_z"'] &
A_{x,x}\ot A_{x,z}\ar[d,"\mlt_{xxz}"]\ar[r,"\comlt_{xyx}\ot1"]
\ar[dr,phantom,"\scriptstyle\cref{frob1}"description] &
A_{x,y}\ot A_{y,x}\ot A_{x,z}\ar[dd,"1\ot\mlt_{yxz}"] \\
& A_{x,z}\ar[dr,"\comlt_{xyz}"]\ar[d,phantom,"\scriptstyle\cref{frob1}"description] & \phantom{A}\\
A_{x,z}\ot A_{z,z}\ar[ur,"\mlt_{xzz}"]\ar[r,"1\ot\comlt_{zyz}"'] &
A_{x,z}\ot A_{z,y}\ot A_{y,z} \ar[r,"\mlt_{xzy}\ot 1"'] & A_{x,y}\ot A_{y,z}
\end{tikzcd}
\end{displaymath}
where the left part is the unit axiom for any $\Vv$-category.
In $k$-linear language, we have $\cas{e}{xy}=\comlt_{xyx}(1_{x,x})$ and the diagram expresses that
for any $a\in A_{x,z}$,
\begin{displaymath}
a\comlt_{zyz}(1_{z,z}) = \comlt_{xyz}(a 1_{z,z})=
\comlt_{xyz}(a) = \comlt_{xyz}(1_{x,x}a) = \comlt_{xyx}(1_{x,x})a.
\end{displaymath}

Conversely, given a Casimir family $e=\{\cas{e}{xy}\}_{x,y\in X}\in W_2$, we 
define a family 
\begin{equation}\label{eq:inducedcomult}
\comlt_{xyz}\colon A_{x,z}\xrightarrow{1\ot\cas{e}{zy}}A_{x,z}\ot A_{z,y}\ot A_{y,z}\xrightarrow{\mlt_{xzy}\ot1}
A_{x,y}\ot A_{y,z} \;\;\stackrel{\cref{casimirdiag}}{=}\; (1\ot\mlt_{yzx})\circ(\cas{e}{xy}\ot 1)
\end{equation}
These indeed satisfy the Frobenius conditions \cref{frob1}: the first is verified by
\begin{displaymath}
 \begin{tikzcd}[column sep=.1in,row sep=.15in]
  A_{x,y}A_{y,z}\ar[rr,"{1\cas{e}{yw}1}"]\ar[dd,"{\mlt_{xyz}}"']\ar[dr,"11\cas{e}{zw}"] && 
  A_{x,y}A_{y,w}A_{w,y}A_{y,z}\ar[rr,"{\mlt_{xyw}11}"]\ar[dr,"{11\mlt_{wyz}}"]
  \ar[dl,phantom,"\scriptstyle\cref{casimirdiag}"description] &&
  A_{x,w}A_{w,y}A_{y,z}\ar[dd,"{1\mlt_{xyz}}"] \\
  & A_{x,y}A_{y,z}A_{w,z}\ar[rr,"{1\mlt_{yzw1}}"]\ar[dr,"{\mlt_{xyz}}"] &&
  A_{x,y}A_{y,w}A_{w,z}\ar[dr,"{\mlt_{xyw}1}"] & \\
  A_{x,z}\ar[rr,"{1\cas{e}{zw}}"'] && A_{x,z}A_{z,w}A_{w,z}\ar[rr,"{\mlt_{xzw}1}"'] &&
  A_{x,w}A_{w,z}
 \end{tikzcd}
\end{displaymath}
where the tensors have been omitted, and similarly for the second. In the $k$-linear case we get, for $a\in A_{x,y}$ and $b\in A_{y,z}$,
\begin{gather*}
\big((1\ot\mlt_{xyz})\circ(\comlt_{xwy}\ot1)\big)(a\ot b)= ae^{1}_{y,w}\ot e^{2}_{w,y}b
\stackrel{\cref{casimirdiag}}{=}abe^{1}_{z,w}\ot e^{2}_{w,z}=\\
\big(\mlt_{xyz}\circ\comlt_{xwz}\big)(a\ot b)=
\big((\mlt_{xyz}\ot1)\circ(1\ot\comlt_{ywz})\big)(a\ot b).
\end{gather*}

The above constructions provide well-defined functions $\alpha:W_1\leftrightarrows W_2:\beta$. Let us check that these constructions are mutual inverses. 
The identity $\beta\circ \alpha(\comlt)=\comlt$ follows form the following diagram:
$$\xymatrix{
A_{x,z} \ar[rr]^-{1\ot \uni_z} \ar@{=}[drr] && A_{x,z}\ot A_{z,z} \ar[rr]^-{1\ot \comlt_{zyz}}\ar@{}[drr]|{{\cref{frob1}}} \ar[d]^{\mlt_{xzz}}
&& A_{x,z}\ot A_{z,y}\ot A_{y,z} \ar[d]^-{\mlt_{xzy}\ot 1} \\
&& A_{x,z} \ar[rr]_-{\comlt_{xyz}}
&& A_{x,y}\ot A_{y,z}}
$$
Conversely, the identity $\alpha\circ \beta(e)=e$ follows from
$$\xymatrix{
I \ar[d]_-{\uni_x} \ar[rr]^{\cas{e}{xy}} &&  A_{x,y}\ot A_{y,x} \ar[d]|-{1\ot1\ot\uni_x} \ar@{=}[drr]\\
A_{x,x} \ar[rr]_-{1\ot \cas{e}{xy}} && A_{x,x}\ot A_{x,y}\ot A_{y,x} \ar[rr]_-{m_{xxy}\ot 1} && A_{x,y}\ot A_{y,x}
}$$
where the left square commutes by naturality and the triangle by the unitality condition of the $\Vv$-category $A$.
\end{proof}

We can now provide a first equivalent characterization of Frobenius categories.

\begin{proposition}\label{FrobeniusCasimir}
For any $\Vv$-category $A$, there is a bijective correspondence between 
\begin{itemize}
 \item comultiplication and counit families $(\comlt,\couni)$ that give $A$ the structure of a Frobenius $\Vv$-category;
\item Casimir families $e$ together with families $\nu=\{\nu_x\colon A_{x,x}\to I\}_{x\in X}$ such that the 
following triangles commute:
\begin{equation}\label{casimir}
\begin{tikzcd}[column sep=.5in]
A_{x,x}\ot A_{x,x}\ar[dr,"\nu_x\ot1"'] & I\ar[d,"\uni_x"]\ar[l,"\cas{e}{xx}"']\ar[r,"\cas{e}{xx}"] &
A_{x,x}\ot A_{x,x}\ar[dl,"1\ot\nu_x"] \\
& A_{x,x} &
\end{tikzcd}
\end{equation}
\end{itemize}
\end{proposition}

In the $k$-linear context, \cref{casimir} is expressed as $\nu_x(e^{1}_{x,x})\cdot e^{2}_{x,x} = e^{1}_{x,x}\cdot\nu_x(e^{2}_{x,x})= 1_{x,x}$.
The families of maps $(e, \nu)$ as above define a \emph{Frobenius system} for any $\Vv$-category $A$.

\begin{proof}
First, suppose that $(A,\mlt,\uni,\comlt,\couni)$ is a Frobenius $\Vv$-category. By \cref{Frobcas}, we know that $\comlt$ gives rise to a 
Casimir family $\cas{e}{xy}\colon I\xrightarrow{\uni_x}A_{x,x}\xrightarrow{\comlt_{xyx}}A_{x,y}\ot A_{y,x}$. 
If we define $\nu_x = \couni_x$, one easily verifies that \cref{casimir} is satisfied
using the counity axiom \cref{Vopcat1} for opcategories.

Conversely, suppose that $(A,\mlt,\uni)$ is a $\Vv$-category with
a Casimir family $e=\{\cas{e}{xy}\}_{x,y}$ and $\nu=\{\nu_x\}_x$ satisfying \cref{casimir}, namely a Frobenius system $(e,\nu)$.
By \cref{Frobcas}, the induced $\comlt_{xyz}=(1\ot\mlt_{yzx})\circ(\cas{e}{xy}\ot 1)$ already belong to $W_1$, and moreover we define coidentities by 
$\couni_x=\nu_x$. Then the coassociativity and counity conditions \cref{Vopcat1}
are satisfied by examining the following diagrams, where $\ot$ has been suppressed and separated subscripts have been concatenated for space purposes:
\begin{displaymath}
\adjustbox{scale=.95,center}{%
\begin{tikzcd}[column sep=.08in,row sep=.1in]
A_{xw}\ar[rr,"{1 \cas{e}{wy}}"]\ar[dd,"{1 \cas{e}{wz}}"'] &&
{A_{xw}  A_{wy}  A_{yw}}\ar[rr,"{\mlt_{xwy} 1}"]\ar[dl,"{111\cas{e}{wz}}"]
&& A_{xy}  A_{yw}\ar[dd,"{1 1 \cas{e}{wz}}"] \\
& A_{xw}A_{wy}A_{yw}A_{wz}A_{zw}\ar[drr,"{11\mlt_{ywz}1}"]
\ar[dd,phantom,"\scriptstyle\cref{casimirdiag}"description] &&& \\
A_{xw}  A_{wz}  A_{zw}\ar[ur,"{1 \cas{e}{wy} 1 1 1}"]
\ar[dr,"{1 1 \cas{e}{zy} 1}"']\ar[dd,"{\mlt_{xwz} 1}"']
&&& A_{xw}A_{wy}A_{yz}A_{zw}\ar[ddr,"{\mlt_{xwy}11}"']
& A_{xy}A_{yz}A_{wz}A_{zw}\ar[dd,"{1 \mlt_{ywz} 1}"] \\
& A_{x,w}  A_{w,z}  A_{zy}  A_{yz}  A_{zw}\ar[urr,"{1 \mlt_{xzy} 1 1}"']\ar[dr,"{\mlt_{xwz}111}"] &&&& \\
A_{xz}  A_{zw}\ar[rr,"{1 \cas{e}{zy} 1}"'] &&
A_{xz}  A_{zy}  A_{yz}  A_{zw}\ar[rr,"{\mlt_{xzy}11}"'] && A_{xy}  A_{yz}  A_{zw}
\end{tikzcd}}
\end{displaymath}
\begin{equation}\label{dualbasisdiag}
\adjustbox{scale=.95,center}{%
\begin{tikzcd}[column sep=.3in,row sep=.15in]
A_{xy}\ar[r,"{1\cas{e}{yx}}"]\ar[dr,"{\cas{e}{xx}1}"description]\ar[ddr,bend right=20,near end,"{\uni_x}1"']
\ar[ddrr,bend right=65,dashed,"\id"description] & 
A_{xy}A_{yx}A_{xy}\ar[r,"{\mlt_{xyx}1}"]\ar[d,phantom,"\scriptstyle\cref{casimirdiag}"] &
A_{xx}A_{xy}\ar[dd,"{\nu_x1}"] \\
\phantom{A}\ar[r,phantom,near end,"\scriptstyle\cref{casimir}"description] &
A_{xx}A_{xx}A_{xy}\ar[d,"{\nu_x11}"]\ar[ur,"{1\mlt_{xxy}}"'] & \\
& A_{xx}A_{xy}\ar[r,"{\mlt_{xxy}}"'] & A_{xy}
\end{tikzcd}\quad
\begin{tikzcd}[row sep=.5in,column sep=.3in]
A_{xy}\ar[drr,bend right=60,dashed,"\id"description]
\ar[dr,bend left=20,phantom,"\scriptstyle\cref{casimir}"description]
\ar[r,"1\cas{e}{yy}"]\ar[dr,"1\uni_x"'] & A_{xy}A_{yy}A_{yy}\ar[r,"\mlt_{xyy}1"]
\ar[d,"11\nu_y"] & A_{xy}A_{yy}\ar[d,"1\nu_y"] \\
& A_{xy}A_{yy}\ar[r,"\mlt_{xyy}"'] & A_{xy}
\end{tikzcd}}
\end{equation}
The unnamed sub-diagrams either commute trivially, or are $\Vv$-category axioms.
In the $k$-linear context, the above is established, for any $a \in A_{x,w}$
and $b\in A_{x,y}$, by
\begin{gather}
(1\ot\comlt_{yzw})\comlt_{xyw}(a)=ae^1_{wy}\ot e^2_{yw}e^1_{wz}\ot e^2_{zw}
\stackrel{\cref{casimirdiag}}{=}ae^1_{wz}e^1_{zy}\ot e^2_{yz}\ot e^2_{wz}=
(\comlt_{xyz}\ot1)\comlt_{xzw}(a)\nonumber \\
\nu_x(be^1_{y,x})\cdot e^{2}_{x,y}\stackrel{\cref{casimirdiag}}{=}
\nu_x(e^{1}_{x,x})\cdot e^{2}_{x,x}b\stackrel{\cref{casimir}}{=}1_{x,x}b=b,
\qquad
be^{1}_{y,y}\cdot\nu_y(e^{2}_{y,y})\stackrel{\cref{casimirdiag}}{=}b1_{y,y}=b
\label{dualbasis}
\end{gather}
Therefore $(A,\comlt,\couni)$ is a $\Vv$-opcategory which also satisfies \cref{frob1}, so it is indeed Frobenius.

The bijectivity of this correspondence follows directly from the bijectivity of the correspondence in \cref{Frobcas}.
\end{proof}

\begin{remark}\label{rem:counisarefunctionals}
Notice that as mere families of morphisms, traces (counits) $\couni$ from the opcategory structure and `functionals' $\nu$ from the Frobenius 
system are basically identical, belonging to $V_1$ as in \cref{eq:comultiplications}. Of course they ultimately satisfy different axioms, but 
as the previous proposition made clear, they are essentially the same hence can be used interchangeably.
\end{remark}

The following lemma establishes a very important property of Frobenius $\Vv$-categories, namely that they are locally rigid (i.e. each hom-object 
has a dual in $\Vv$) in a natural way.

\begin{lemma}\label{Frobisrigid}
Any Frobenius $\Vv$-category $A$ is locally rigid, with $A_{x,y}^*\cong A_{y,x}$ for any two objects $x,y$.
\end{lemma}

\begin{proof}
Since $A$ is equipped with a Frobenius system $(e,\nu)$, the evaluation and coevaluation maps can be defined as
\begin{gather*}
\cas{\ev}{xy}\colon A_{y,x}\ot A_{x,y}\xrightarrow{\mlt_{yxy}}A_{y,y}\xrightarrow{\nu_y}I \\
\cas{\coev}{xy}\colon I\xrightarrow{\cas{e}{xy}}A_{x,y}\ot A_{y,x}
\end{gather*}
and the two commutative diagrams \cref{dualbasisdiag} verify that $A_{y,x}$ is the dual of $A_{x,y}$.
\end{proof}

\begin{remark}\label{rem:dualbases}
In the $k$-linear context for a commutative ring $k$, we know that rigid objects are exactly finitely generated
and projective modules where the dual is given by all linear functionals. Hence
\cref{dualbasis}, which establishes \cref{Frobisrigid}, expresses the dual base property \cref{eq:dualcond} exhibiting 
$\{e^2_{x,y},\nu_x(-e^1_{y,x})\}$ as a finite dual basis for each $k$-module $A_{x,y}$.
Notice that 
$\{e^1_{x,y},\nu_x(e^2_{y,x}-)\}$ also constitutes a dual basis for $A_{x,y}$
since
\begin{displaymath}
a = a e^{1}_{y,y}\cdot\nu_y(e^{2}_{yy})= e^{1}_{x,y}\cdot\nu_{y}(e^{2}_{y,x}a).
\end{displaymath}
\end{remark}

\subsection{Characterization in terms of dual module structure}

One of the equivalent definitions of a (classical) Frobenius $k$-algebra $A$ is that $A$ is finite dimensional and isomorphic to its dual $A^{*}$ as a 
right $A$-module. In this section, we generalize this to Frobenius $\Vv$-categories.

For any locally rigid $\Vv$-category $A$, consider the right $A$-module $A^\dagger$ that is constructed as $A^\sop$ out of the left 
regular $A$-module as in \cref{daggermodules} and the left  $A$-module ${}^\dagger A$ that is constructed out of the right regular $A$-module. We 
denote
\begin{eqnarray*}
V_3&=&\Hom_A(A,A^\dagger), \qquad V'_3={_A\Hom}(A,{}^\dagger A)\\
W_3&=&\Hom_A(A^\dagger,A), \qquad W'_3={_A\Hom}({}^\dagger A,A)
\end{eqnarray*}
where for example $\Hom_A(A,A^\dagger)$ denotes the set of all right $A$-module morphisms from $A$ to $A^\dagger$: an element $f$ consists of maps
$f_{xy}\colon A_{x,y}\to A_{y,x}^*$ which satisfy \cref{VModmaps}, here
\begin{equation}\label{eq:Frobisomodule}
 \begin{tikzcd}[column sep=.4in]
 A_{x,y}\ot A_{y,z}\ar[d,"f_{xy}1"']\ar[rrr,"\mlt_{xyz}"] &&& A_{x,z}\ar[d,"f_{xz}"]\\
A^*_{y,x}\ot A_{y,z}\ar[r,"11\cas{\coev}{zx}"'] & A^*_{y,x}\ot A_{y,z}\ot A_{z,x}\ot A_{z,x}^*\ar[r,"1\mlt_{yzx}1"'] & A^*_{y,x}\ot 
A_{y,x}\ot A^*_{z,x}\ar[r,"\cas{\ev}{yx}1"'] & A^*_{z,x}
\end{tikzcd}
 \end{equation}

\begin{lemma}\label{VW3}
For any locally rigid $\Vv$-category $A$, there exist bijections
\begin{enumerate}
\item $V_1\cong V_3\cong V'_3$,
\item $W_1\cong W_3\cong W_3'$
\end{enumerate}
where $V_1$ and $W_1$ are as in \cref{eq:comultiplications}.
\end{lemma}

\begin{proof}
First notice that $V_3\cong V'_3$ and $W_3\cong W'_3$ follow by construction of the module structures on $A^\sop$.

\ul{(1)}. Given a family $\nu=\{\nu_x\}\in V_1$ (see \cref{rem:counisarefunctionals}), we can define
\begin{equation}\label{eq:rightAmodmap}
\psi_{xy}\colon A_{x,y}\xrightarrow{1\ot \cas{\coev}{yx}}A_{x,y}\ot A_{y,x}\ot A^{*}_{y,x}\xrightarrow{\mlt_{xyx}\ot1} 
A_{x,x}\ot A^{*}_{y,x}\xrightarrow{\nu_{x}\ot 1}A^{*}_{y,x} 
\end{equation}
Thes maps can be easily verified to satisfy \cref{eq:Frobisomodule}, so they form a right $A$-module morphism from $A$ to $A^\dagger$.
This defines a map $V_1\to V_3$. Conversely, given some $\psi\in V_3$, we define a family of morphisms
\begin{equation}\label{eq:inducednu}
\nu_x\colon A_{x,x}\xrightarrow{\uni_x\ot1}A_{x,x}\ot A_{x,x}\xrightarrow{\psi_{xx}\ot1}
A^*_{x,x}\ot A_{x,x}\xrightarrow{\cas{\ev}{xx}} I
\end{equation}
establishing a map $V_3\to V_1$. One can then verify that these two directions are mutual inverses.

\ul{(2)}. Recall that by \cref{Frobcas}, $W_1\cong W_2$, the set of Casimir families. Given a Casimir family $\cas{e}{xy}$, we can define an element $\phi\in W_3$ by means of the composition
\begin{displaymath}
\phi_{xy}\colon A^{*}_{y,x}\xrightarrow{1\ot\cas{e}{yx}}A^{*}_{y,x}\ot A_{y,x}\ot A_{x,y}\xrightarrow{\cas{\ev}{yx}\ot1}A_{x,y}
\end{displaymath}
A similar computation as in part (1) shows that $\phi$ is also a right $A$-module morphism using the Casimir property \cref{casimirdiag}.
Conversely, given $\phi\in W_3$, we claim that 
\begin{displaymath}
\cas{e}{xy}\colon I\xrightarrow{\cas{\coev}{xy}}A_{x,y}\ot A_{x,y}^*\xrightarrow{1\ot\phi_{yx}}A_{x,y}\ot A_{y,x}
\end{displaymath}
form a Casimir family. Indeed, the following commutativity
verifies \cref{casimirdiag}:
\begin{displaymath}
\begin{tikzcd}[row sep=.3in]
A_{xz}\ar[r,"{\cas{\coev}{xy}1}"]\ar[dd,"{1\cas{\coev}{zy}}"']
& A_{xy}A^*_{xy}A_{xz}\ar[d,"{111\cas{\coev}{zy}}"]\ar[r,"{1\phi_{yx}1}"] &
A_{xy}A_{yx}A_{xz}\ar[r,"{1\mlt_{yxz}}"]\ar[d,phantom,"{\scriptstyle(*)}"description] & A_{xy}A_{yz}\\
& A_{xy}A^*_{xy}A_{xz}A_{zy}A^*_{zy}\ar[r,"{11\mlt_{xzy}}"] &
A_{xy}A^*_{xy}A_{xy}A^*_{zy}\ar[r,"{1\cas{\ev}{xy}1}"] & A_{xy}A^*_{zy}\ar[u,"{1\phi_{yz}}"] \\
A_{xz}A_{zy}A^*_{zy}\ar[urr,phantom,"{\scriptstyle(**)}"description]
\ar[ur,"{\cas{\coev}{xy}111}"]\ar[rrr,"{11\phi_{yz}}"']
\ar[urrr,"{\mlt_{xzy}1}"'] &&& A_{xz}A_{zy}A_{yz}\ar[uu,bend right=40,"{\mlt_{xzy}1}"']
\end{tikzcd}
\end{displaymath}
where $(*)$ commutes since $\phi$ is a right $A$-module map between the regular $A$-module and its dual $A^\dagger$, and $(**)$ is the triangle 
equality for evaluation and coevaluation.
The above constructions provide well-defined maps between $W_1\cong W_2$ and $W_3$, which can be checked to be mutual inverses.
\end{proof}

\begin{remark}\label{rem:localrigiditynotessential}
One can observe that if $\Vv$ is a closed monoidal category, the bijections $V_1\cong V_3\cong V'_3$ are still valid without the assumption of 
local rigidity -- replacing $X^*$ by $[X,I]$. On the other hand, in this setting we still have well-defined maps $W_1\to W_3$ and $W_1\to W'_3$, 
but the proof for bijectivity is only valid under local rigidity of the $\Vv$-category $A$.

In the $k$-linear case, the Casimir family $e$ in function of a $\varphi\in W_3$ is explicitly given by
$\cas{u_i}{xy}\ot\varphi_{yx}(\cas{u_i}{\;xy_*})$,
where $\lbrace (\cas{u_i}{xy},\cas{u_i}{\;xy_*})\rbrace$ is a dual base for $A_{x,y}$.
The Casimir property in this case is explicitly checked as follows
\begin{align*}
a\cas{u_i}{zy}\ot\varphi_{yz}(\cas{u_i}{\;zy_*})& = \cas{u_i}{xy}\cdot\cas{u_i}{\;xy_*}(a\cas{u_i}{zy})
\ot\varphi_{yz}(\cas{u_i}{\;zy_*})=\cas{u_i}{xy}\ot\varphi_{yz}(\cas{u_i}{\;xy_*}(a\cas{u_i}{zy})\cdot\cas{u_i}{\;zy_*})= \\
& =\cas{u_i}{xy}\ot\varphi_{yz}(\cas{u_i}{\;xy_*}(a-))=\cas{u_i}{xy}\ot\varphi_{yx}(\cas{u_i}{\;xy_*})a,
\quad \forall a\in A_{x,z}.
\end{align*}
\end{remark}

The following result then gives the characterization of Frobenius $\Vv$-categories in terms of `Frobenius isomorphisms', namely $A$-module 
isomorphisms with $A^\sop$.

\begin{proposition}\label{frobmodules}
For a locally rigid $\Vv$-category $A$, there is a bijective correspondence between:
\begin{enumerate}
\item Frobenius systems on $A$;
\item isomorphisms between the right $A$-modules $A$ and $A^\dagger$; 
\item isomorphisms between the left $A$-modules $A$ and ${}^\dagger A$.
\end{enumerate}
In particular, a $\Vv$-category $A$ is Frobenius if and only if it is locally rigid and $A\cong A^\sop$ as right $A$-modules.
\end{proposition}
\begin{proof}
Recall that the right $A$-module structure on $A^\sop$ is the bottom of \cref{eq:Frobisomodule}.
We only prove the equivalence between (1) and (2) since the equivalence with (3) follows by symmetry. The last statement follows immediately from the 
stated correspondence in combination with \cref{Frobisrigid}.

Given a Frobenius system $(e, \nu)$ for $A$, by \cref{VW3} we can construct two right $A$-linear morphisms $\psi:A\leftrightarrows A^\dagger:\phi$. 
It can easily be checked that $\phi$ and $\psi$ are inverses using the (co)evaluation condition together with \cref{casimirdiag} and \cref{casimir}.

Conversely, given an isomorphism of right $A$-modules  $\psi:A\to A^\dagger$ with inverse $\phi$, we can obtain a Casimir element $E$ and a family $(\nu_x:A_{x,x}\to I)_{x\in X}$ from \cref{VW3}. Let us check that these make up a Frobenius system.
The left side of \cref{casimir} follows from
\begin{displaymath}
\begin{tikzcd}[row sep=.15in,column sep=.6in]
I\ar[r,"{\cas{\coev}{xx}}"]\ar[d,"{\uni_x}"'] &
A_{xx}A^*_{xx}\ar[d,"{\uni_x11}"]\ar[r,"{1\phi_{xx}}"] & A_{xx}A_{xx}\ar[r,"{\uni_x11}"]
& A_{xx}A_{xx}A_{xx}\ar[d,"{\psi_{xx}11}"] \\
A_{xx}\ar[d,"{\psi_{xx}}"']
& A_{xx}A_{xx}A^*_{xx}\ar[d,"{\psi_{xx}1}"]
&& A^*_{xx}A_{xx}A_{xx}\ar[d,"{\cas{\ev}{xx}1}"] \\
A^*_{xx}\ar[r,"{1\cas{\coev}{xx}}"]\ar[rr,dotted,bend right=10,"{\id}"description] &
A^*_{xx}A_{xx}A^*{xx}\ar[r,"{\cas{\ev}{xx}1}"]\ar[urr,bend left=10,"{11\phi_{xx}}"] &
A^*_{xx}\ar[r,"{\phi_{xx}}"] & A_{xx}
\end{tikzcd}
\end{displaymath}
using the evaluation-coevaluation property and $\phi$, $\psi$ being inverses.

For the right hand side of \cref{casimir}, we first notice that the induced $\nu_x$ constructed as \cref{eq:inducednu} are equivalently given by
\begin{displaymath}
A_{x,x}\xrightarrow{\psi_{xx}}A^*_{x,x}\xrightarrow{1\ot\uni_x}A^*_{x,x}\ot A_{x,x}
\xrightarrow{\cas{\ev}{xx}}I.
\end{displaymath}
This follows from $\psi$ being an $A$-module morphism and again the triangle equa\-li\-ties. The remaining
verification is now straightforward:
\begin{displaymath}
\begin{tikzcd}
I\ar[r,"\cas{\coev}{xx}"]\ar[d,"\uni_x"'] & A_{xx}A^*{xx}\ar[r,"{1\phi_{xx}}"]
\ar[rr,dotted,bend right=10,"\id"description]
& A_{xx}A_{xx}\ar[r,"{1\psi_{xx}}"] & A_{xx}A^*_{xx}\ar[r,"{11\uni_x}"] &
A_{xx}A^*_{xx}A_{xx}\ar[d,"{1\cas{\ev}{xx}}"] \\
A_{xx}\ar[urrrr,bend right=5,"{\cas{\coev}{xx}1}"]\ar[rrrr,"\id"description] &&&& A_{xx}
\end{tikzcd}
\end{displaymath}
\end{proof}

\begin{remark}
In the linear case, the equivalent formulations for $\nu_x$ are computed 
\begin{align*}
\nu_x(a)&=[\psi_{xx}(1_{xx})](a)=[\psi_{xx}(1_{xx})](a1_{xx}){=}[\psi_{xx})(1_{xx})\cdot a](1_{xx})
=\psi_{xx}(1_{xx}a)(1_{x,x})\\
&=\psi_{xx}(a)(1_{xx})
\end{align*}
from which it follows that
\begin{displaymath}
\cas{u_i}{xx}\cdot\nu_x(\varphi_{xx}(\cas{u_i}{\;xx*})){=}
\cas{u_i}{xx}\cdot[\psi_{x,x}(\varphi_{xx}(\cas{u_i}{\;xx*}))](1_{x,x})=
\cas{u_i}{xx}\cdot\cas{u_i}{\;xx*}(1_{x,x})=1_{xx}.
\end{displaymath}
\end{remark}

The following result shows that in fact, the Frobenius isomorphism $A\cong A^\sop$ is a $\Vv$-opcategory one rather than just an $A$-module one.

\begin{theorem}\label{Frobisoisopcatiso}
Let $(A,\mlt,\uni)$ be a Frobenius $\Vv$-category with Frobenius system $(e,\nu)$ and consider the dual $\Vv$-opcategory $(A^\sop,\comlt^{A^{\sop}}, \couni^{A^{\sop}})$ as in \cref{A*opcategory}.
Then the induced isomorphisms $\psi_{xy}\colon A_{x,y}\to A^*_{y,x}$ as in \cref{VW3} form an isomorphism of $\Vv$-opcategories $A\to A^\sop$.
\end{theorem}

\begin{proof}
Recall by \cref{FrobeniusCasimir} that the induced comultiplication on $A$ in terms of its Frobenius system $(e,\nu)$ is $\comlt^{A}_{xyz} = (\mlt_{xzy} \ot A_{y,z}) \circ (A_{x,z} \ot \cas{e}{zy})$ as in \cref{eq:inducedcomult} and the induced counit on $A$ is just $\couni=\nu$.
Furthermore, \cref{VW3} describes the morphisms $\psi_{xy}$ in terms of the Frobenius system as
$\psi_{xy} = (\nu_{x} \ot A^{*}_{y,x}) \circ (\mlt_{xyx} \ot A^{*}_{y,x}) \circ (A_{x,y} \ot \cas{\coev}{yx})$ as in \cref{eq:rightAmodmap}.

The following diagram shows that $\psi$ preserves the cocomposition
\begin{displaymath}
\begin{adjustbox}{max width=\textwidth}
\begin{tikzcd}
A_{xz}\arrow[r,"1\cas{\coev}{zx}"] \arrow[d,"1\cas{e}{zy}"']\arrow[dr,"1\cas{\coev}{zy}"]& A_{xz}A_{zx}A^{*}_{zx}\arrow[r,"\mlt_{xzx}1"] & A_{xx} A^{*}_{zx}\arrow[r,"\nu_{x}1"] & A^{*}_{zx} \arrow[r,"1\cas{\coev}{zy}"]& A^{*}_{zx}A_{zy}A^{*}_{zy}\arrow[d,"11\cas{\coev}{yx}1"]\\
A_{xz}A_{zy}A_{yz}\arrow[d,"\mlt_{xzy}1"'] & A_{xz}A_{zy}A^{*}_{zy}\arrow[d,"1\cas{e}{zy}11"]\arrow[dr,"11\cas{e}{yy}1"] & & & A^{*}_{zx}A_{zy}A_{yx}A^{*}_{yx}A^{*}_{zy}\arrow[dd,"1\mlt_{zyx}11"]\\
A_{xy}A_{yz}\arrow[d,"11\cas{\coev}{zy}"'] & A_{xz}A_{zy}A_{yz}A_{zy}A^{*}_{zy}\arrow[d,"11\mlt_{yzy}1"']\arrow[r,phantom,"\cref{casimirdiag}"] & A_{xz}A_{zy}A_{yy}A_{yy}A^{*}_{zy}\arrow[d,"111\nu_{y}1"]\arrow[dl,"1\mlt_{zyy}11"] & & \\
A_{xy}A_{yz}A_{zy}A^{*}_{zy}\arrow[d,"1\mlt_{yzy}1"'] & A_{xz}A_{zy}A_{yy}A^{*}_{zy}\arrow[d,"11\nu_{y}1"] & A_{xz}A_{zy}A_{yy}A^{*}_{zy}\arrow[dl,"1\mlt_{zyy}1"'] \arrow[rr,phantom,"(*)"]& & A^{*}_{yx}A^{*}_{zy}\\
A_{xy}A_{yy}A^{*}_{zy}\arrow[d,"1\nu_{y}1"'] & A_{xz}A_{zy}A^{*}_{zy}\arrow[dl,"\mlt_{xzy}1"]\arrow[dr,"11\cas{\coev}{yx}1"] &&&\\
A_{xy}A^{*}_{zy} \arrow[ddr,bend right,"1\cas{\coev}{yx}"']& & A_{xz}A_{zy}A_{yx}A^{*}_{yx}A^{*}_{zy}\arrow[d,"1\mlt_{zyx}11"] & & \\
& & A_{xz}A_{zx}A^{*}_{yx}A^{*}_{zy}\arrow[d,"1\mlt_{zyx}11"']  & & \\
 & A_{xy}A_{yx}A^{*}_{yx}A^{*}_{zy}\arrow[r,"\mlt_{xyx}11"'] & A_{xx}A^{*}_{yx}A^{*}_{zy} \arrow[rruuuu,bend right,"\nu_{x}11"'] & &
\end{tikzcd}
\end{adjustbox}
\end{displaymath}
where the inner diagram $(*)$ commutes because of \cref{casimir}, definition of $\Vv$-category and the evaluation-coevalution property. A similar diagram proves the counit condition.
\end{proof}

Notice that the previous characterizations of Frobenius $\Vv$-categories can be reformulated in terms of $\Vv$-opcategories. For 
example, since a $\Vv$-opcategory $C$ gives rise to a $\Vv$-category $C^\sop$ by \cref{A*opcategory}, \cref{frobmodules} would 
accordingly state that a $\Vv$-opcategory $C$ is Frobenius if and only if $C$ is locally rigid and isomorphic to $C^\sop$ as right 
$C$-opmodules, using the regular structure for $C$ and the one from \cref{Csopopmodule} for $C^\sop$.
In that case, one has the following corollary.

\begin{corollary}\label{dualfrob}
A $\Vv$-category $A$ is Frobenius if and only if the $\Vv$-opcategory $A^\sop$ is Frobenius.
\end{corollary}
\begin{proof}
This follows from the following equivalences.

\begin{center}
\begin{tabular}{rcl}
$\Vv$-category $A$ is Frobenius & $\Leftrightarrow$ & $A \cong A^\sop$ as right $A$-modules \\
& $\Leftrightarrow$ & $A^\sop \cong A$ as right $A^\sop$-opmodules \\
& $\Leftrightarrow$ & $\Vv$-opcategory $A^\sop$ is Frobenius
\end{tabular}
\end{center}
The first and last equivalences are \cref{frobmodules} and its dual statement, and the middle equivalence is \cref{prop:staropmodules}.
\end{proof}

Finally, the symmetry of the Frobenius definition is also expressed as follows.

\begin{proposition}\label{prop:isomodopmodFrob}
If $A$ is a Frobenius $\Vv$-category, then the categories $\Vv$-$\mathsf{opMod}_A$ and $\Vv$-$\Mod_A$ of $A$-modules and opmodules are isomorphic.
\end{proposition}

\begin{proof}
Suppose the $\Vv$-category $(A,\mlt,\uni)$ comes with a Frobenius system $(e,\nu)$, and $(N,\gcoact)$ is an $A$-opmodule and $(M,\gact)$ an 
$A$-module.
Define a functor $F\colon \Vv\textrm{-}\mathsf{opMod}_A\to\Vv\textrm{-}\Mod_A$ by $F(N_{x,y}) = N_{x,y}$ with action 
\begin{displaymath}
 N_{x,y}\ot A_{y,z}\xrightarrow{\gcoact_{xzy}\ot1}N_{x,z}\ot A_{z,y}\ot A_{y,z}\xrightarrow{1\ot\mlt_{zyz}}N_{x,y}\ot 
A_{z,z}\xrightarrow{1\ot\nu_z}N_{x,y}
\end{displaymath}
mapping an $A$-opmodule map to the same morphism in $A$ which can be shown to commute with the above defined actions.
Furthermore, define $G\colon\Vv\textrm{-}\Mod_A \to \Vv\textrm{-}{\sf opMod}_A$ by $G(M_{x,y}) = M_{x,y}$ with coaction 
\begin{displaymath}
 M_{x,y}\xrightarrow{1\ot\cas{e}{yz}}M_{x,y}\ot A_{y,z}\ot A_{z,y}\xrightarrow{\gact_{xyz}\ot1}M_{x,z}\ot A_{z,y}
\end{displaymath}
Those two functors are inverse to one another, and the proof is complete.
\end{proof}

\subsection{Characterization in terms of trace morphisms}

In this subsection, we provide yet another characterization of Frobenius $\Vv$-categories, generalizing the classical characterization of Frobenius 
algebras in terms of properties of the trace morphisms. We moreover show how this is related to so-called `Calabi-Yau' categories. 

\begin{definition}\label{def:bilinearform}
For any $\Vv$-graph $G$, a \emph{bilinear form} $\Gamma$ is a collection of morphisms $\Gamma_{xy} \colon
G_{x,y}\otimes G_{y,x} \to I$ in $\Vv$. 
If $\Vv$ is braided, a bilinear form is said to be \emph{symmetric} when
\begin{equation*}
\cd[@!C@R-2em]{
  G_{x,y} \ot G_{y,x} \ar[dd]_{\braid} \ar@/^.6em/[dr]^-{\Gamma_{xy}} &  \\
  & I \\
  G_{y,x} \ot G_{x,y} \ar@/_.6em/[ur]_-{\Gamma_{yx}} &  \\
}
\end{equation*}
commutes for all $x,y$. 

If $A$ is a $\Vv$-category and $\Gamma$ is a bilinear form on $A$, we say that $\Gamma$ is {\em balanced} when the following commutes
\begin{equation}\label{eq:balancedGamma}
\xymatrix{
A_{x,y}\ot A_{y,z}\ot A_{z,x} \ar[rr]^-{\mlt_{xyz}\ot 1} \ar[d]_{1\ot \mlt_{yzx}} && A_{x,z}\ot A_{z,x} \ar[d]^{\Gamma_{xz}}\\
A_{x,y}\ot A_{y,x} \ar[rr]_-{\Gamma_{xy}} && I
}
\end{equation}
\end{definition}

We denote the set of all balanced bilinear forms on a $\Vv$-category $A$ by
$$V_4=\{\Gamma=\{\Gamma_{xy}\}_{x,y\in X}\;|\;\Gamma_{xy}\colon A_{x,y}\ot A_{y,x}\to I \textrm{ satisfying }\cref{eq:balancedGamma}\}$$
With notation as in \cref{eq:comultiplications}, we obtain a correspondence to the traces or functionals on $A$, see \cref{rem:counisarefunctionals}.

\begin{lemma}\label{V4}
For any $\Vv$-category $A$, there is a bijection $V_1\cong V_4$.
\end{lemma}

\begin{proof}
For any $\nu=\{\nu_x\colon A_{x,x}\to I\}_{x}\in V_1$, define a bilinear form by
\begin{equation}\label{eq:inducedGamma}
\Gamma_{xy}\colon A_{x,y}\ot A_{y,x}\xrightarrow{\mlt_{xyx}}A_{x,x}\xrightarrow{\nu_x}I.
\end{equation}
From the associativity of $A$, we immediately obtain that $\Gamma$ is balanced. 

Conversely, given a balanced bilinear form $\Gamma$ on $A$, we define 
$$\nu_x\colon A_{x,x}\xrightarrow{\uni_x \ot 1}A_{x,x}\ot A_{x,x}\xrightarrow{\Gamma_{xx}}I$$
Since $\Gamma$ is balanced, we also have that $\nu_x=\Gamma_{xx}\circ (1\ot\uni_x)$. 

It can be easily verified that these two constructions are inverses.
\end{proof}

For any locally rigid $\Vv$-category $A$, \cref{VW3,V4} establish that 
\begin{equation}\label{eq:V134}
V_1\cong V_3\cong V'_3\cong V_4.
\end{equation}
A balanced bilinear form $\Gamma\in V_4$ or a corresponding family of trace morphisms $\nu\in V_1$ will be called left (respectively right) {\em 
non-degenerate} if the corresponding element in $V_3$ (resp. $V'_3$) is a split monomorphism. Spelled out in the more general setting of a monoidal closed 
category (see \cref{rem:classicalduals}) where $\eta\colon Y\to[X,Y\ot X]$ is the tensor-hom adjuction unit, this leads to the following definition.

\begin{definition}\label{def:leftrightnondegenerate}
Suppose $\ca{V}$ is monoidal closed. A bilinear form $\Gamma$ on a $\Vv$-category $A$ is {\em left non-degenerate} when all maps
\begin{displaymath}
\Gamma_{xy}^1\colon A_{x,y} \xrightarrow{\eta}[A_{y,x},A_{x,y}\ot A_{y,x}]\xrightarrow{[1,\Gamma_{xy}]}[A_{y,x},I]=A_{y,x}^*
\end{displaymath}
are split monomorphisms in $\Vv$. If $A$ is locally rigid, this says that 
\begin{displaymath}
\Gamma_{xy}^1\colon A_{x,y} \xrightarrow{1\ot\cas{\coev}{yx}}A_{x,y}\ot A_{y,x}\ot A_{y,x}^*\xrightarrow{\Gamma_{xy}\ot1}A_{y,x}^*
\end{displaymath}
are split monomorphisms.
Symmetrically, $\Gamma$ is {\em right non-degenerate} when the maps
\begin{displaymath}
\Gamma_{xy}^2\colon A_{x,y} \xrightarrow{\eta}[A_{y,x},A_{x,y}\ot A_{y,x}]\xrightarrow{[1,\invbraid]}[A_{y,x},A_{y,x}\ot A_{x,y}]\xrightarrow{[1,\Gamma_{yx}]}[A_{y,x},I]=A_{y,x}^*
\end{displaymath}
are split monomorphisms in $\Vv$. We say that $\Gamma$ is {\em non-degenerate} if and only if it is both left and right non-degenerate.

Equivalently, a family $\nu=\{\nu_x\}_x$ of traces is {\em non-degenerate} if and only if the corresponding bilinear form \cref{eq:inducedGamma} is; 
in the locally rigid setting, this is the case when \cref{eq:rightAmodmap} and its `switched' (using braiding) are split monomorphisms in $\Vv$.
\end{definition}

In the $k$-linear case (where $k$ is a field) we find that $\Gamma^1_{xy}(a) \in A^*_{x,y}$ is defined by the formula $\Gamma^1_{xy}(a)(b)=\Gamma_{xy}(a\ot b)$. So 
$\Gamma^1_{xy}$ being a split monomorphism (which is in this case equivalent to being a mononomorphism or to being injective) means that if $\Gamma_{xy}(a\ot b)=0$ for all $b$, then $a=0$ as usually stated.

The next lemma gives several sufficient conditions to obtain non-degeneracy of a bilinear form.

\begin{lemma}\label{nondegenerate}\hfill
\begin{enumerate}
\item If $A$ is a Frobenius $\Vv$-category, there exists a non-degenerate balanced bilinear form on $A$.
\item If a bilinear form is symmetric, then left non-degeneracy is equivalent to right non-degeneracy (hence to non-degeneracy).
\item If $\Vv={\sf Vect}_k$ and $A$ is locally rigid then for a bilinear form, left non-degeneracy is equivalent to right non-degeneracy (hence to 
non-degeneracy).
\end{enumerate}
\end{lemma}

\begin{proof}
\ul{(1)}. If $A$ is Frobenius, then the counit
of the opcategory structure corresponds by \cref{V4} to a balanced bilinear form $\Gamma$ on $A$. On the 
other hand, we know by \cref{frobmodules} that $A$ is isomorphic as a right $A$-module with $A^\dagger$ and as a left $A$-module with ${}^\dagger A$. 
As observed above, these isomorphisms are exactly given by $\Gamma^1$ and $\Gamma^2$, hence in particular $\Gamma$ is already 
non-degenerate.\\
\ul{(2)}. If a bilinear form $\Gamma$ is symmetric, it follows directly from the definitions that $\Gamma^1=\Gamma^2$.\\
\ul{(3)}. In this case, local rigidity means that $\dim A_{x,y}=\dim A^*_{y,x}$ and injective maps between vector spaces of the same finite dimension 
are automatically bijective.
\end{proof}

We now proceed to the equivalent characterization of Frobenius structures in terms of bilinear forms, and moreover for symmetric ones as in 
\cref{eq:symmetricFrob}.

\begin{proposition}\label{CalabiYauEquiv}
Let $A$ be a locally rigid $\Vv$-category. There is a bijective 
correspondence between the following:
\begin{enumerate}
\item (symmetric) Frobenius structures on $A$;
\item (symmetric) non-degenerate balanced bilinear forms on $A$;
\item (symmetric) non-degenerate families of trace morphisms on $A$.
\end{enumerate}
\end{proposition}

\begin{proof}
We already know by \cref{V4} that for a given (locally rigid) $\Vv$-category, there is a bijective correspondence between trace families (candidate 
counits) for a Frobenius structure on $A$ and balanced bilinear forms on $A$, which carries on to non-degenerate ones by \cref{def:leftrightnondegenerate}.

From the construction of this correspondence, it is clear that a given 
Frobenius structure is symmetric if and only if the associated bilinear form is so. Furthermore, we know from \cref{nondegenerate} that the bilinear 
form $\Gamma\in V_4$ corresponding to the counits $\nu\in V_1$ of a Frobenius system is non-degenerate.

It only remains to prove that a non-degenerate balanced bilinear form on $A$ endows $A$ with a (unique) Frobenius system. Let $\Gamma$ be such a 
bilinear form: by \cref{eq:V134} we know that it corresponds to a right $A$-linear morphism $\Gamma^1:A\to A^\dagger$ 
and a left $A$-linear morphism $\Gamma^2:A\to {}^\dagger A$. Moreover, non-degeneracy of $\Gamma$ implies that all $\Gamma^1_{xy}$ and 
$\Gamma^2_{xy}$ are split monomorphisms in $\Vv$ and so they have left inverses that we will denote respectively 
by $\Lambda^1_{xy}$ and $\Lambda^2_{xy}$. 
Define
\begin{displaymath}
\ol \Lambda^2_{xy}\colon A^*_{y,x}\xrightarrow{1\ot\cas{\coev}{xy}}A^*_{y,x}\ot A_{x,y}\ot A^*_{x,y}\xrightarrow{\braid\ot\Lambda^2_{yx}}A_{x,y}\ot A^*_{y,x}\ot A_{y,x}\xrightarrow{1\ot\cas{\ev}{yx}}A_{x,y}
\end{displaymath}
Using the fact that $\Lambda^2$ is inverse to $\Gamma^2$ and the evaluation/coevaluation condition, we can verify that $\ol\Lambda^2$ is a 
right inverse of $\Gamma^1_{xy}$. Therefore, $\Gamma^1$ has both a left inverse $\Lambda^1$ and a right inverse $\ol\Lambda^2$ so is an 
isomorphism, and as a result $A$ is Frobenius by \cref{frobmodules}.
\end{proof}

Now recall from \cite[p.\ 176]{Costello2007} that a \emph{Calabi-Yau category} is a $k$-linear category for a field $k$, equipped with 
a family of trace maps $\Tr_x \colon A_{x,x} \to I$ with the property that for all $x,y\in X$ the associated pairing 
$\Gamma_{xy}=\Tr_x\circ\mlt_{xyx}$ is non-degenerate and symmetric.
In our notation, if we substitute $\Tr_x$ by $\nu_x$ and recall 
\cref{eq:inducedGamma}, we immediately obtain the 
following characterization of Calabi-Yau categories as a corollary of the previous proposition.

\begin{corollary}
A locally rigid $k$-linear category is Calabi-Yau if and only if it is symmetric Frobenius.
\end{corollary}

\subsection{Characterization in terms of adjoint functors}

Another classical characterization of Frobenius algebras says that a $k$-algebra $A$ is Frobenius if and only if the forgetful functor $U:\Mod_A\to \Mod_k$ that forgets the $A$-action is a \emph{left} adjoint of the free functor $-\ot A\colon\Mod_k\to\Mod_A$. This is the origin of the alternative `Frobenius functor' terminology that refers to adjoints which are both left and right to the same functor. Before we generalize this result, let us first observe that for any $\Vv$-category we have the usual free-forgetful adjunction.

\begin{proposition}
Let $A$ be a $\Vv$-category. Then the forgetful functor $$U\colon\Vv\textrm{-}\Mod_A\to \Vv\textrm{-}{\sf dGrph},$$ defined on objects by $U(M,\gact)\mapsto \{M_{x,x}\}_{x\in X}$ has a left adjoint $F=-\ot A$.
\end{proposition}

\begin{proof} 
Similarly to \cref{fundhopfmod}, the functor $F=-\ot A\colon\Vv\textrm{-}{\sf dGrph}\to \Vv\textrm{-}\Mod_A$ maps a diagonal graph $\{N_x\}_{x}$ to the right $A$-module $\{N_x\ot A_{x,y}\}_{x,y}$ with action $1\ot \mlt_{xyz}$.

For any $N\in \Vv\textrm{-}{\sf dGrph}$, define $\alpha_N^x=1\ot \uni_x\colon N_x\to N_x\ot A_{x,x}$ and for any $(M,\gact)\in \Vv\textrm{-}\Mod_A$ define $\beta_M^{xy}=\gact_{xyz}\colon M_{x,x}\ot A_{x,y}\to M_{x,y}$. It can verified that $\alpha$ and $\beta$ are the unit and counit of the desired adjunction.
\end{proof}

The following lemma gives the essence of the characterization in terms of the adjunction. 

\begin{lemma}\label{VWnatural}
Let $(\Vv,\otimes, I)$ be a category such that the monoidal unit $I$ is a (regular) generator and all endofunctors $-\ot V:\Vv\to \Vv$ preserve 
(regular) epimorphisms.
For a $\Vv$-category $A$, the sets $V_1$ of traces \cref{eq:comultiplications} and $W_2$ of Casimir families \cref{casimirdiag} are moreover in bijection to the following sets of natural transformations
\begin{enumerate}
\item $V_1\cong \Nat(U\circ F,\id);$
\item $W_2\cong \Nat(\id,F\circ U).$
\end{enumerate}
\end{lemma}

\begin{proof}
(1) Consider a family $\couni=\{\couni_x\}_{x\in X}\in V_1$. We define a natural transformation $\beta\colon U\circ F\to \id$ for any object $N\in \Vv\textrm{-}{\sf dGrph}$ by
$$\beta_N^x=1\ot\couni_x:N_x\ot A_{x,x}\to N_x.$$
Conversely, given $\beta: U\circ F\to \id$, for the unit diagonal graph $D=\{I\}_x$ we obtain $\{\beta_D^x\colon A_{x,x}\to I\}_{x\in X}\in V_1$. 
Now for any $\Vv$-${\sf dGrph}$-morphism $n:D\to N$ with components $n_x:I\to N_x$, we find by naturality of $\beta$ that $\beta_N^x\circ (n_x\ot 1)=n_x\circ \beta_D^x$. On the other hand, $I$ is a generator so the morphisms $n:I\to N$ are jointly epimorphic;  since endofunctors $-\ot V$ preserve (regular) epimorphisms, it follows that $(n\ot 1)$ are jointly epimorphic as well. Therefore, we obtain that $\beta_N^x=1\ot \beta_D^x$ for all $N$ and $x$, from which we deduce that the above defined maps between $V_1$ and $\Nat(U\circ F,\id)$ are bijections.

(2) Let $e\in W_2$ be a Casimir family. Define a natural transformation $\id\to F\circ U$ by
\begin{equation}\label{defalphafromcas}
\alpha_M^{xy}\colon M_{x,y}\xrightarrow{1\ot \cas{e}{yx}}M_{x,y}\ot A_{y,x}\ot A_{x,y}\xrightarrow{\mlt_{xyx}\ot 1}M_{x,x}\ot A_{x,y}
\end{equation}
for any right $A$-module $M$ and any $x,y\in X$. The fact that $\alpha_M^{xy}$ is a right $A$-module morphism follows from the Casimir property of 
$e$, and naturality of $\alpha$ is immediate.

Now let $\alpha:\id\to F\circ U$ be a natural transformations and fix $y\in X$. If $h:X\to X$ is defined by $h(z)=y$ for all objects $z\in X$, recall that the shuffle $A^h$ of the regular right $A$-module, as defined in \cref{Hopfcats}, is given by $A^h_{z,u}=A_{y,u}$ for all $z,u\in X$ with a right action given by $\gact^h_{zuv}=\mlt_{yuv}$. Then we define
\begin{equation}\label{defcasfromalpha}
\cas{e}{yx}\colon I\xrightarrow{\uni_y}A_{y,y}=A^h_{x,y} \xrightarrow{\alpha_{A^h}^{xy}}A^h_{x,x}\ot A_{x,y} = A_{y,x}\ot A_{x,y}
\end{equation}
which is verified to satisfy the Casimir property in the end of the current proof.

In case the natural transformation $\alpha$ arises from a Casimir family as in \eqref{defalphafromcas}, we see by the commutativity of the following diagram that the morphisms defined in \eqref{defcasfromalpha} are exactly the initial ones:
\[
\xymatrix{
I \ar[rr]^-{\uni_y} \ar[d]_{\cas{e}{yx}} && A_{y,y}{=}A^h_{x,y} \ar[d]_{\cas{e}{yx}1} \ar[rr]^-{1e^{yx}} 
&& A^{h}_{x,y}A_{y,x}A_{x,y}=A_{y,y}A_{y,x}A_{x,y} \ar[d]^-{\mlt_{yyx}1} \\
A_{y,x}A_{x,y}\ar[rr]^-{11\uni_y} && A_{y,x}A_{x,y}A_{y,y} \ar[rr]^-{1\mlt_{xyy}} && 
A^\tau_{x,x}A_{x,y}=A_{y,x}A_{x,y}
}
\]
Conversely, we will show that for a natural transformation $\alpha:\id\to F\circ U$, the maps defined as in \eqref{defalphafromcas} where $\cas{e}{yx}$ is constructed as in \eqref{defcasfromalpha} are exactly the components of the $\alpha$. To this end, for a fixed $y\in X$, any right $A$-module $M$ gives rise to a right $A$-module $MA$ whose components and action are given by
$$(MA)_{u,v}=M_{u,y}\ot A_{y,v}, \qquad \mu^{MA}_{uvw}=1\ot \mlt_{yvw} :(MA)_{u,v}\ot A_{v,w} \to (MA)_{u,w}$$
for all $u,v,w\in X$. Then for any family of morphisms $n_u:I\to M_{u,y}$ with $u\in X$, we obtain a right $A$-module morphism 
$n:A^h\to MA$ given by
$$n_{uv}:\xymatrix{A^h_{u,v}=A_{y,v} \ar[rr]^-{n_u\ot 1} && M_{u,y}\ot A_{y,v}=(MA)_{u,v}}$$
for all $u,v\in X$.
As a consequence, we find that the following diagram commutes
\[
\xymatrix@R=.2in{
(MA)_{x,y}{=}M_{x,y}A_{y,y} \ar[rr]^-{\alpha_{MA}^{xy}} && (MA)_{x,x}A_{x,y}{=}M_{x,y}A_{y,x}A_{x,y}\\
& A^\tau_{x,x}A_{x,y}{=}A_{y,x}A_{x,y} \ar[ur]^-{n_x11}\\
A^\tau_{x,y}{=}A_{y,y} \ar[uu]^-{n_{x}1} \ar[ur]^{\alpha^{xy}_{A^h}} \ar[rr]_-{n_x1} 
&& M_{x,y}A^\tau_{x,y}{=}M_{x,y}A_{y,y} \ar[uu]_-{1\alpha_{A^h}^{xy}} 
}
\]
The commutativity of the upper triangle follows by the naturality of $\alpha$ applied to the right $A$-module morphism $n$, and the lower triangle commutes by naturality of the tensor product. Hence we find that $\alpha_{MA}^{xy}\circ (n_x\ot 1)=(1\ot \alpha_{A^h}^{xy})\circ (n_x\ot 1)$.
Since this hold for all choices of the morphisms $n_x:I\to M_{x,y}$, we obtain that 
\begin{equation}\label{equalityalpha}
\alpha_{MA}^{xy}=1\ot \alpha_{A^h}^{xy}
\end{equation}
by the generator condition on $I$.
Thanks to this identity, the following diagram commutes by naturality of $\alpha$
\[\xymatrix{
M_{x,y} \ar[rr]^-{1\uni_y} \ar@{=}[drr] && M_{x,y}A_{y,y}{=}M_{x,y}A^h_{x,y} \ar[d]_{\mu_{xyy}} 
\ar[rr]^-{\alpha_{MA}^{xy}{=}1\ot \alpha_{A^\tau}^{xy}} && M_{x,y}A^h_{x,x}A_{x,y}{=}M_{x,y}A_{y,x}A_{x,y}
\ar[d]^{\mu_{xyx}1}\\
&& M_{x,y} \ar[rr]_-{\alpha^{xy}_M} && M_{x,x}A_{x,y}
}\]
Hence $\alpha^{x,y}_M=(\mu_{xyx}\ot 1)\circ(1\ot \cas{e}{yx})$, where the morphisms $\cas{e}{yx}$ are defined as in \eqref{defcasfromalpha}.

The proof is complete, subject to the verification of the Casimir property \cref{casimirdiag} of $\cas{e}{yx}$ \eqref{defcasfromalpha} by the following diagram:
\[\xymatrix{
A_{x,y} \ar[r]^-{1\uni_y} \ar[d]_{\uni_x1} \ar@{=}[dr] & A_{x,y}A_{y,y}{=}A_{x,y}A^{h'}_{z,y} 
\ar[rr]^-{1\alpha_{A^{h'}}^{zy} =\alpha_{A^h A}^{zy} } \ar[d]^{\mlt_{xyy}}
&& A_{x,y}A_{y,z}A_{z,y}{=}A_{x,y}A^{h'}_{z,z}A_{z,y} \ar[d]^{\mlt_{xyz}1}\\
A^h_{z,x}A_{x,y}{=}A_{x,x}A_{x,y} \ar[r]_{\mlt_{xxy}} \ar[d]_{\alpha_{A^h}^{zx}1} 
& A_{x,y} {=} A^{h}_{z,y} \ar[rr]^{\alpha^{zy}_{A^h}} && A^h_{z,z}A_{z,y}{=}A_{x,z}A_{z,y} \ar@{=}[d] \\
A^h_{z,z}A_{z,y}A_{x,y}{=}A_{x,z}A_{z,y}A_{x,y} \ar[rrr]^-{1\mlt_{zxy}} &&& A_{x,z}A_{z,y}
}
\]
The functions $h,h'\colon X\to X$ are defined by $h(u)=x$ and $h'(u)=y$ for all $u\in X$,
the right $A$-module $A^h A$ has $(A^h A)_{u,v}=A_{x,y}\ot A_{y,v}$ and $1\ot \alpha_{A^{h'}}^{zy} =\alpha_{A^h A}^{zy}$ is just \eqref{equalityalpha} applied to the right $A$-module $A^h$.
The right upper square commutes by naturality of $\alpha$
and the lower square commutes by right $A$-linearity of $\alpha_{A^h}$. 
\end{proof}

Recall (e.g. \cite{Caenepeel2002}) that a functor $F:\Cc\to \Dd$ is called {\em Frobenius} if it has isomorphic left and right adjoints. If we denote the (left or right) adjoint of $F$ by $U$, then we also say that $(F,U)$ is a \emph{Frobenius pair of functors}. Hence we obtain the following characterization of Frobenius $\Vv$-categories.

\begin{theorem}
Under the same assumptions as \cref{VWnatural},
a $\Vv$-category is Frobenius if and only if the functor $-\ot A\colon\Vv$-${\sf dGrph}\to \Vv$-$\Mod_A$ is Frobenius.
\end{theorem}
Precisely, there is a bijective correspondence between Frobenius structures on $A$ and pairs of natural transformations $U\circ F\to\id$ and $\id\to F\circ U$ making $(F,U)$ into a Frobenius pair of functors.
\begin{proof}
It can be verified that $F$ is a right adjoint of $U$ with unit $\alpha$ and counit $\beta$ if and only if the corresponding Casimir family $e$ and family of maps $\epsilon=\nu$ by means of \cref{VWnatural} satisfy the Frobenius system conditions \eqref{casimir}.
\end{proof}

\section{The Larson-Sweedler theorem}
\label{sec:mainresults}

In this section, having introduced all the required structures, we proceed to the main goal of this work: a
generalization of the Larson-Sweedler theorem for Hopf $\Vv$-categories. We first briefly recall the original setting for $k$-algebras over a field or principal ideal domain, then we generalize integral theory for Hopf $\Vv$-categories and finally we prove the main results relating Hopf and Frobenius structures on a $\Vv$-category.

\subsection{Classical Larson-Sweedler Theorem}\label{sec:classicLS}

Let us recall
the original statement found in~\cite{LS}.

\begin{theorem*}[Larson-Sweedler]\label{classicalLarsonSweedler}
Let $H$ be a finite dimensional bialgebra over the principal ideal domain $R$. Then the following conditions are equivalent:
\begin{enumerate}
\item there exists an antipode for $H$;
\item there exists a non-singular left integral in $H$.
\end{enumerate}
If $\Lambda$ is a non-singular left integral in $H$, and $\Lambda_1$ is any left integral in $H$, there exists $a \in R$ such that $\Lambda_1 = a\Lambda$.
\end{theorem*}
By a left integral in $H$, one means an element $t\in H$ satisfying $ht=\epsilon(h)t$ for all $h\in H$. An integral is moreover {\em non-singular} if the linear maps 
\begin{align}
p: H^*\to H,& \qquad p(f)=f(t_{(1)})t_{(2)}\label{classicalnonsingular}\\
q: H^*\to H,& \qquad q(f)=t_{(1)}f(t_{(2)})\nonumber
\end{align}
are bijective. Notice that when working over a field, the existence of an isomorphism $H\cong H^*$ implies the finite dimensionality of $H$ -- so that this assumption can be dropped in part (2) of the above statement.
Moreover, one can verify that the composition of $p$ with the antipode of the Hopf algebra $H$ yields a right $H$-linear map $H^*\to H$. Since the antipode of a finite-dimensional Hopf algebra is automatically bijective (see one-object case of \cref{antipodeisinvertible}) this right $H$-linear map is an isomorphism, giving a Frobenius structure on $H$.
Because of this, and since any Frobenius algebra is finite dimensional (see one-object case of \cref{frobmodules}), the Larson-Sweedler theorem is sometimes rephrased by saying that a Hopf algebra is finite dimensional if and only if it is Frobenius. 

However, an important remark should be made here: it is of course well-known that a bialgebra can only have one unique antipode. Moreover, the 
Larson-Sweedler theorem says that a non-singular integral in a finite dimensional Hopf algebra exists, and such an integral is also unique (up to 
scalar multiplication). On the other hand, there can exist many Frobenius structures on the same finite dimensional Hopf algebra. However, only one 
of these Frobenius structures will correspond exactly to a non-singular integral!

Let us clarify the above by the following example. Of course, any group algebra $kG$ over a finite group $G$ has the structure of a finite dimensional Hopf algebra, with multiplication and unit extending those of the group and comultiplication, counit and antipode extending $g\mapsto g\ot g$, $g\mapsto e_G$ and $g\mapsto g^{\textrm{-}1}$. The Larson-Sweedler theorem then ensures that it has a non-singular integral, which in this case is given by the element $\sum_{g\in G} g$, and the Casimir element of the corresponding Frobenius structure on such a Hopf algebra is given by $\sum_{g\in G} g\ot g^{-1}$. 

In particular, for $G=C_4$ the cyclic group with four elements generated by $g$, we find that a Frobenius system for $kG$ is given by the Casimir element
$$e\ot e+ g\ot g^3 + g^2\ot g^2 + g^3\ot g$$
and linear functional $\delta_e$, the dual base vector of $e$.
On the other hand, one can easily verify that there is another Frobenius system on $kC_4$ given by the Casimir element
$$e\ot g + g\ot e + g^2\ot g^3+ g^3\ot g^2$$
and linear functional $\delta_g$, the dual base vector of $g$.
It is known (see \cref{prop:integralCasimir} below) that any Casimir element $e^1\ot e^2$ on a Hopf algebra $H$ leads to a left integral $e^1\epsilon(e^2)$ and on the other hand any left integral $t$ leads to a Casimir element $t_{(1)}\ot s(t_{(2)})$. However, this correspondence is not bijective: in general, there are more Casimir elements than integrals. In case of the example $C_4$, we see that both Casimir elements for $kC_4$ above give the same integral $t=e+g+g^2+g^3$. On the other hand, starting with this integral, the above construction gives back only the first Casimir element.

Finally, let us also remark that a bialgebra can be Frobenius without being Hopf: indeed, this is the case exactly when the Frobenius structure does not correspond to an integral. Take for example the monoid algebra $kM$ where $M=\{e,g~|~g^2=g\}$. Then $kM$ is Frobenius via the system with Casimir element
$$e\ot e+g\ot g$$ 
and linear functional $\delta_e$. However, since $M$ is a monoid and not a group, the bialgebra $kM$ is not a Hopf algebra. One can also see that the only integral in $kM$ is $g$, and this integral is singular since $p(g)(\delta_e)=\delta_e(g)g=0$.

\subsection{Integral theory for semi-Hopf $\Vv$-categories}
In this section, we generalize the theory of integrals in the many-object setting, which is necessary for the expression and proof of the main \cref{oppositeLarsonSweedler2} as well as intermediate results. 
For what follows, fix $(A,\mlt,\uni,\lcomlt,\lcouni)$ to be a semi-Hopf $\Vv$-category for a braided monoidal category $\Vv$.

\begin{definition}\label{integralfamilies}
A {\em left integral family} for $A$ is a collection $t=\{\cas{t}{xy}\colon I\to A_{x,y}\}_{x,y\in X}$ of morphisms in $\Vv$
that satisfy the commutativity of
\begin{equation}\label{integraldiag}
\begin{tikzcd}[row sep=.2in,column sep=.7in]
A_{z,x} \ar[rr,"1\ot\cas{t}{xy}"]\ar[dd,"1\ot\cas{t}{zy}"']\ar[dr,"\lcouni_{zx}"description]
&& A_{z,x}\ot A_{x,y}\ar[dd,"\mlt_{zxy}"]\\
& I\ar[dr,"\cas{t}{zy}"description]& \\
A_{z,x}\ot A_{z,y}\ar[rr,"\lcouni_{zx}\ot1"'] && A_{z,y}
\end{tikzcd} 
\end{equation}
where the bottom triangle commutes trivially.
In fact, if we consider the unit $\Vv$-graph $\Ii$ given by $\Ii_{x,y}=I$ viewed as a left $A$-module via $\lcouni_{xy}\ot1\colon A_{x,y}\ot I\to I$,
a left integral family can equivalently be viewed as an identity-on-objects
left $A$-module morphism $t\colon\Ii\to A$.

A {\em right integral family} for $A$ is defined symmetrically via the property 
\begin{equation}\label{rightintegraldiag}
\mlt_{xyz}\circ(\cas{t}{xy}\ot1)=\cas{t}{xz}\circ\lcouni_{yz}. 
\end{equation}
\end{definition}

In the $k$-linear case, morphisms $k\to A_{x,y}$ can be identified with elements of the vector space, and a left integral family is then can be written in the form
$$\{\cas{t}{xy}\in A_{x,y}~|~ a\cas{t}{xy}=\lcouni_{zx}(a)\cdot\cas{t}{zy}, \forall z\in X \ {\textrm{and\ }} \forall a\in A_{z,x}\}.$$

The following result establishes the close relationship between integral and Casimir families from \cref{def:casimirfamily}, for (semi-) Hopf categories.

\begin{proposition}\label{prop:integralCasimir}
Every Casimir family $e=\{\cas{e}{xy}\}_{x,y\in X}$ for $A$ gives rise to a left integral family $t_e$ via
\begin{equation}\label{t(e)}
\cas{t}{xy}:= I\xrightarrow{\cas{e}{xy}}A_{x,y}\ot A_{y,x}
\xrightarrow{1\ot\lcouni_{yx}}A_{x,y}
\end{equation}
If $A$ is moreover Hopf, every left integral family $t=\{\cas{t}{xy}\}_{x,y\in X}$ gives rise to
a Casimir family $e_t$ via
\begin{equation}\label{e(t)}
\cas{e}{xy}:=I\xrightarrow{\cas{t}{xy}}A_{x,y}
\xrightarrow{\lcomlt_{xy}}A_{x,y}\ot A_{x,y}\xrightarrow{1\ot\atpd_{xy}}
A_{x,y}\ot A_{y,x}
\end{equation}
In fact, for every integral family $t$ in a Hopf category it holds that $t_{e_t}=t$. 
\end{proposition}

\begin{proof}
 In order to verify \cref{integraldiag}, we examine the following commutative diagram
\begin{displaymath}
\begin{tikzcd}[column sep=.7in]
A_{z,x}\ar[r,"1\ot\cas{e}{xy}"]\ar[d,"\cas{e}{zy}\ot1"']
\ar[dr,phantom,"\scriptstyle\cref{casimirdiag}"description] &
A_{z,x}\ot A_{x,y}\ot A_{y,x}\ar[r,"1\ot1\ot\lcouni_{yx}"]\ar[d,"\mlt_{zxy}\ot1"] &
A_{z,x}\ot A_{x,y}\ar[d,"\mlt_{zxy}"] \\
A_{z,y}\ot A_{y,z}\ot A_{z,x}\ar[r,"1\ot\mlt_{yzx}"']
\ar[rr,bend right=15,"{1\ot\lcouni_{yz}\ot\lcouni_{zx}}"']
\ar[rr,bend right=7,phantom,"\scriptstyle\cref{Hax1}"description]
& A_{z,y}\ot A_{y,x}\ar[r,"1\ot\lcouni_{yx}"'] & A_{z,y}
\end{tikzcd}
\end{displaymath}
where the bottom composite is precisely $\cas{t}{zy}\ot\lcouni_{zx}$.

Conversely, the Casimir condition \cref{casimirdiag} can be verified via the following calculation
\begin{align*}
\mlt_{xzy}1 &\circ 11\atpd_{zy} \circ 1\lcomlt_{zy} \circ 1\cas{t}{zy} \\
&\stackrel{(*)}{=} \mlt_{xzy}1 \circ 11\atpd_{zy} \circ 1\lcomlt_{zy} \circ 1\cas{t}{zy}\circ 1\lcouni_{xz} \circ \lcomlt_{xz}\\
&= \mlt_{xzy}1 \circ 111\lcouni_{xz} \circ 11\atpd_{zy}1 \circ 1\lcomlt_{zy}1 \circ 1\cas{t}{zy}1 \circ \lcomlt_{xz}\\
&\stackrel{(**)}{=}\mlt_{xzy}1 \circ 11\mlt_{yzz} \circ 111\uni_{z} \circ 111\lcouni_{xz} \circ 11\atpd_{zy}1 \circ 1\lcomlt_{zy}1 \circ 1\cas{t}{zy}1 \circ \lcomlt_{xz}\\
&\stackrel{\cref{HopfCatAntipodeEquations}}{=}  \mlt_{xzy}1 \circ 11\mlt_{yzz} \circ 111\mlt_{zxz} \circ 111\atpd_{xz}1 \circ 111\lcomlt_{xz} \circ 11\atpd_{zy}1 \circ 1\lcomlt_{zy}1 \circ 1\cas{t}{zy}1 \circ \lcomlt_{xz}\\
&= \mlt_{xzy}1 \circ 11\mlt_{yzz} \circ 111\mlt_{zxz} \circ 111\atpd_{xz}1 \circ 11\atpd_{zy}11 \circ 111\comlt_{xz} \circ 1\lcomlt_{zy}1 \circ 1\cas{t}{zy}1 \circ \lcomlt_{xz}\\
&= \mlt_{xzy}1 \circ 11\mlt_{yzz} \circ 111\mlt_{zxz} \circ 111\atpd_{xz}1 \circ 11\atpd_{zy}11 \circ 1\comlt_{zy}11 \circ 11\comlt_{xz}\circ 1\cas{t}{zy}1 \circ \lcomlt_{xz}\\
&\stackrel{\ref{antipodeproperties}}{=} 1\mlt_{yxz} \circ \mlt_{xzy}11 \circ 11\atpd_{xy}1\circ 11\mlt_{xzy}1 \circ 11 \braid^{-1}1 \circ 1 \lcomlt_{zy}11 \circ 1 \braid^{-1}1 \circ 11\cas{t}{zy}1 \circ \comlt_{xz}1\\
&= 1\mlt_{yxz} \circ 1\atpd_{xy}1 \circ 1\mlt_{xzy}1 \circ 1 \braid^{-1}1 \circ \mlt_{xzy}111 \circ 1 \lcomlt_{zy}11 \circ 1 \braid^{-1}1 \circ 11\cas{t}{zy}1 \circ \comlt_{xz}1\\
&\stackrel{(***)}{=}1\mlt_{yxz} \circ 1\atpd_{xy}1 \circ \mlt_{xzy}\mlt_{xzy}1 \circ 1 \braid 11 \circ \lcomlt_{xz}\lcomlt_{zy}1\circ 1 \cas{t}{zy}1 \circ \lcomlt_{xz}\\
&\stackrel{\cref{Hax1}}{=} 1\mlt_{yxz} \circ 1\atpd_{xy}1 \circ \lcomlt_{xy}1 \circ \mlt_{xzy}1 \circ 1 \cas{t}{zy}1 \circ \lcomlt_{xz}\\
&\stackrel{\cref{integraldiag}}{=} 1\mlt_{yxz} \circ 1\atpd_{xy}1 \circ \lcomlt_{xy}1 \circ \lcouni_{xz}11 \circ 1 \cas{t}{xy}1 \circ \lcomlt_{xz}\\
&\stackrel{(*)}{=} 1\mlt_{yxz} \circ 1\atpd_{xy}1 \circ \lcomlt_{xy}1 \circ \cas{t}{xy}1
\end{align*}
Explicitly, $(*)$ uses the local comultiplication, $(**)$ the $\Vv$-category structure and $(***)$ the naturality of the braiding with $\braid_{I,A} = \braid^{-1}_{I,A}\cong\id_A$. Finally, taking the integral family of an integral-induced Casimir family returns the initial one as follows:
\begin{displaymath}
\begin{tikzcd}[row sep=.2in]
I\ar[r,"\cas{t}{xy}"] & A_{x,y}\ar[r,"\lcomlt_{xy}"]\ar[drr,bend right=10,"1"description]
& A_{x,y}\ot A_{x,y}\ar[r,"1\ot\atpd_{xy}"]\ar[dr,"1\ot\lcouni_{xy}"description]
& A_{x,y}\ot A_{y,x}\ar[d,"1\ot\lcouni_{yx}","\scriptstyle(\ref{antipodeproperties})"'] \\
&&& A_{x,y}
\end{tikzcd}
\end{displaymath}
\end{proof}

\begin{remark}
As we already discussed in \cref{classicalLarsonSweedler}, in general a Hopf $\Vv$-category (or just a Hopf algebra) has more Casimir families than integral families. For example, if the semi-Hopf category $A$ has an op-antipode $\ol\atpd$ (e.g.\ $A$ is locally rigid and Hopf), given a left integral $t$, one can easily check that the following family is also Casimir
$$\xymatrix{I\ar[r]^-{\cas{t}{xy}} & A_{x,y} \ar[r]^-{\lcomlt_{xy}}  &A_{x,y}\ot A_{x,y} \ar[r]^-\sigma  &A_{x,y}\ot A_{x,y} \ar[r]^-{1\ot\ol\atpd_{yx}} & A_{x,y}\ot A_{y,x}}$$
In case $A$ is commutative or cocommutative, then this Casimir element is the same as the one from \eqref{e(t)}.
In the same way, starting with a right integral family $t$, we find that the following composites
\begin{gather}\label{e(t)right}
\xymatrix{I\ar[r]^-{\cas{t}{xy}} & A_{x,y} \ar[r]^-\delta  &A_{x,y}\ot A_{x,y}  \ar[r]^-{\atpd_{xy}\ot 1} & A_{y,x}\ot A_{x,y}} \\
\xymatrix{I\ar[r]^-{\cas{t}{xy}} & A_{x,y} \ar[r]^-\delta  &A_{x,y}\ot A_{x,y} \ar[r]^-\sigma  &A_{x,y}\ot A_{x,y} \ar[r]^-{\ol\atpd_{yx}\ot 1} & A_{x,y}\ot A_{y,x}}\nonumber
\end{gather}
both form Casimir families.
\end{remark}

Notice that a Casimir family for a $\Vv$-category $A$ is in particular a bilinear form (\cref{def:bilinearform}) on $A$, viewed as a $\Vv^\op$-graph.
Since Casimir elements can be constructured from integrals by the above \cref{prop:integralCasimir}, non-degeneracy of bilinear forms as in \cref{def:leftrightnondegenerate} can be traced back to integrals as well. This lies at the origin of the notion of non-singularity for an integral family as expressed below, where the split monomorphism condition of non-degeneracy here corresponds to a split epimorphism condition due to $\ca{V}^\op$. 

Suppose the enriching base $\Vv$ is monoidal closed, and denote $A^*_{x,y}=[A_{x,y},I]$ which is not necessarily the categorical dual of $A_{x,y}$ as discussed in \cref{rem:classicalduals}. For any (left or right) integral family $t=\{\cas{t}{xy}\}$ of $A$, we define two families of 
morphisms 
\begin{equation}\label{eq:pxy}
p_{xy}\colon A^*_{x,y}\xrightarrow{1\ot \cas{t}{xy}}A^*_{x,y}\ot A_{x,y}\xrightarrow{1\ot \lcomlt_{xy}}
A^*_{x,y}\ot A_{x,y}\ot A_{x,y} \xrightarrow{\cas{\ev}{xy}\ot 1}A_{x,y}
\end{equation}
\begin{equation}\label{eq:qxy}
q_{xy}\colon
A^*_{x,y} \xrightarrow{\cas{t}{xy}\ot 1}A_{x,y}\ot A^*_{x,y} \xrightarrow{\lcomlt_{xy}\ot 1}
A_{x,y}\ot A_{x,y}\ot A^*_{x,y}\xrightarrow{1\ot\braid}A_{x,y}\ot A^*_{x,y}\ot A_{x,y}\xrightarrow{1\ot\cas{\ev}{xy}}A_{x,y}
\end{equation}

\begin{definition}\label{def_non_singular}
Suppose $A$ is a semi-Hopf $\Vv$-category, where $\Vv$ is braided monoidal closed.
A (left or right) integral family $t=\{\cas{t}{xy}\}$ is called \emph{left non-singular} if all maps $p_{xx}$ for $x\in X$
are split epimorphisms.

Similarly, the (left or right) integral family $t$ is called {\em right non-singular} if all maps $q_{xx}$ for $x\in X$
are split epimorphisms. If the integral family $t$ is both left and right non-singular, we say that it is {\em non-singular}.
\end{definition}

In the $k$-linear case for any commutative ring, the composites of \cref{def_non_singular} are
\begin{align*}
p_{xy}&:A^*_{x,y}\to A_{x,y},\;\; p_{xy}(f)=f(\cas{t}{xy}_{(1)})\cdot \cas{t}{xy}_{(2)}\\
q_{xy}&:A^*_{x,y}\to A_{x,y},\;\; q_{xy}(f)=\cas{t}{xy}_{(1)}\cdot f(\cas{t}{xy}_{(2)})
\end{align*}
If $k$ is a field and $A$ is locally finite dimensional, then non-singularity implies that all $p_{xx}$ and $q_{xx}$ are isomorphisms.
In other words, if a $k$-linear semi-Hopf category (with $k$ a field) has a non-singular integral then all diagonal bialgebras $H_{xx}$ have a non-singular integral in the classical sense, as in \cref{classicalnonsingular}.
Moreover, we will later prove \cref{oppositeLarsonSweedler}, (\ref{itemfour}) that if a semi-Hopf $\Vv$-category $A$ has a non-singular left and right integral, then all maps $p_{xy}$ and $q_{xy}$ are isomorphisms. 

The following result shows how to construct a right integral family from a left one, using an invertible antipode.

\begin{proposition}\label{leftrightintegral}
Suppose $(H,\mlt,\uni,\lcomlt,\lcouni,\atpd)$ is a Hopf $\Vv$-category with invertible antipode. 
If $t=\{\cas{t}{xy}\colon I\to H_{x,y}\}_{x,y}$ is a left integral family, then 
$$\atpd\circ t:=\{ I\xrightarrow{\cas{t}{yx}}H_{y,x}\xrightarrow{\atpd_{yx}}H_{x,y}\}_{x,y}$$ is a  right integral family for $H$. Moreover, if $t$ is left (right) non-singular, then $\atpd\circ t$ is right (left) non-singular.
\end{proposition}

\begin{proof}
The right integral property \cref{rightintegraldiag} of $\{\atpd_{yx}\circ\cas{t}{yx}\}_{x,y}$ can be verified by the following computation, where e.g. $\atpd_{yx}H_{y,z}$ denotes $\atpd_{yx}\ot 1_{H_{y,z}}$:
\begin{align*}
\mlt_{xyz} \circ \atpd_{yx}H_{y,z} \circ \cas{t}{yx}H_{y,z}
&= \mlt_{xyz} \circ H_{x,y}\atpd_{zy} \circ H_{x,y}\atpd^{-1}_{zy}\circ \atpd_{yx}H_{y,z} \circ \cas{t}{yx}H_{y,z}\\
&= \mlt_{xyz} \circ H_{x,y}\atpd_{zy} \circ  \atpd_{yx}H_{z,y} \circ H_{y,x}\atpd^{-1}_{zy}\circ  \cas{t}{yx}H_{y,z}\\
&= \atpd_{zx} \circ \mlt_{zyx} \circ \braid^{-1} \circ H_{y,x}\atpd^{-1}_{zy}\circ \cas{t}{yx}H_{y,z}\\
&= \atpd_{zx} \circ \mlt_{zyx} \circ H_{z,y}\cas{t}{yx} \circ \atpd^{-1}_{zy}\\
&= \atpd_{zx} \circ \cas{t}{zx} \circ \lcouni_{xy} \circ \atpd^{-1}_{zy}\\
&= \atpd_{zx} \circ \cas{t}{zx} \circ \lcouni_{yz} \\
\end{align*}
where we used \cref{antipodeproperties} and \cref{integraldiag}, namely the condition that left integrals satisfy.

Furthermore, to show that $\atpd\circ t$ is right non-singular when $t$ is left non-singular, as per \cref{def_non_singular} we need to find a right-sided inverse $\ol{q}_{x}$ to the composite 
$$H^*_{x,x}\xrightarrow{\cas{t}{xx}1}H_{x,x}\ot H^*_{x,x}\xrightarrow{\atpd_{xx}1}H_{x,x}\ot H^*_{x,x}\xrightarrow{\lcomlt_{xx}1}H_{x,x}\ot H_{x,x}\ot H^*_{x,x}\xrightarrow{1\braid}H_{x,x}\ot H^*_{x,x}\ot H_{x,x}\xrightarrow{1\cas{\ev}{xx}} H_{x,x}
$$

For $t$, we know that there exists a $\ol{p}_{x}$ such that $p_{xx} \circ \ol{p}_{x} = H^{*}_{x,x}$ for $p_{xx}$ as in \cref{eq:pxy}. It can now be verified, using that $\lcomlt_{yx}\circ\atpd_{xy}=\braid\circ(\atpd_{xy}\ot\atpd_{xy})\circ\lcomlt_{xy}$ from \cref{antipodeproperties}, that 
$$\ol{q}_{x}:=H_{x,x}\xrightarrow{\atpd^{-1}_{xx}}H_{x,x}\xrightarrow{\ol{p}_{x}}H^*_{x,x}\xrightarrow{(\atpd_{xx}^{-1})^*}H^*_{x,x}
$$
is the required splitting, where in the monoidal closed setting $f^*$ just means $[f,1]$.
\end{proof}

Remark that in the previous proof, we only used that the antipode is an anti-Hopf category morphism. Hence one proves in the same way that if $t$ is a right integral family, then $\atpd\circ t$ is a left integral family.

The following technical lemma will be needed in the proof of our main \cref{oppositeLarsonSweedler2} and relates the inverses of the split epimorphisms in the non-singularity condition to one another via an invertible antipode.

\begin{lemma}\label{leftrightnonsingular}
Suppose $t$ is a non-singular left integral family for a Hopf $\Vv$-category $H$ with invertible antipode.
For any two maps $f_x,g_x \colon H_{x,x} \to I$, 
if the composites
\begin{gather*}
I \xrightarrow{\cas{t}{xx}}H_{x,x} \xrightarrow{\lcomlt_{xx}}H_{x,x} \ot H_{x,x}\xrightarrow{1\ot f_{x}}H_{x,x} \\
I \xrightarrow{\cas{t}{xx}} H_{x,x} \xrightarrow{\lcomlt_{xx}}H_{x,x} \ot H_{x,x}\xrightarrow{g_{x}\ot 1}H_{x,x} 
\end{gather*}
are both equal to $\uni_{x}$, then $g_x=f_x\circ \atpd_{xx}^{-1}$.

In particular, if $\ol{q}_{x}$ is a right inverse of the split epimorphism $q_{xx}$ \cref{eq:pxy} and $\ol{p}_{x}$ a right inverse of $p_{xx}$ accordingly, then 
$$\cas{\ev}{xx} \circ (\ol{p}_{x}\ot1) \circ (\uni_{x}\ot1) = \cas{\ev}{xx}\circ \braid \circ (1\ot \ol{q}_{x}) \circ (1\ot \uni_{x}) \circ \atpd^{-1}_{xx}\colon H_{x,x}\longrightarrow I.$$
\end{lemma}

\begin{proof}
The result follows from the following computation:
\begin{align*}
g_{x} &= g_{x} \circ \mlt_{xxx} \circ H_{x,x}\uni_{x}\\
&= g_{x} \circ \mlt_{xxx} \circ H_{x,x}H_{x,x}f_{x} \circ H_{x,x}\lcomlt_{xx} \circ H_{x,x}\cas{t}{xx}\\
&= g_{x} \circ \mlt_{xxx} \circ H_{x,x}H_{x,x}f_{x} \circ H_{x,x}H_{x,x}\atpd^{-1}_{xx} \circ H_{x,x}H_{x,x}\atpd_{xx}  \circ H_{x,x}\lcomlt_{xx} \circ H_{x,x}\cas{t}{xx}\\
&= g_{x} \circ H_{x,x}f_{x} \circ H_{x,x}\atpd^{-1}_{xx}\circ \mlt_{xxx}H_{x,x} \circ H_{x,x}H_{x,x}\atpd_{xx}  \circ H_{x,x}\lcomlt_{xx} \circ H_{x,x}\cas{t}{xx}\\
&\leftstackrel{(*)}{=} g_{x} \circ H_{x,x}f_{x} \circ H_{x,x}\atpd^{-1}_{xx}\circ H_{x,x}\mlt_{xx} \circ H_{x,x}\atpd_{xx} H_{x,x} \circ \lcomlt_{xx}H_{x,x} \circ \cas{t}{xx}H_{x,x}\\
&\leftstackrel{(**)}{=} g_{x} \circ H_{x,x}f_{x} \circ H_{x,x}\mlt_{xx}\circ H_{x,x}\braid^{-1}\circ H_{x,x}\atpd^{-1}_{xx}\atpd^{-1}_{xx}  \circ H_{x,x}\atpd_{xx}H_{x,x} \circ \lcomlt_{xx}H_{x,x} \circ \cas{t}{xx}H_{x,x}\\
&= f_{x} \circ \mlt_{xxx} \circ \braid^{-1} \circ H_{x,x}\atpd^{-1}_{xx} \circ g_{x}H_{x,x} \circ \lcomlt_{x,x}H_{x,x} \circ \cas{t}{xx}H_{x,x}\\
&= f_{x} \circ \mlt_{xxx} \circ \braid^{-1} \circ H_{x,x}\atpd^{-1}_{xx} \circ \uni_{x}H_{x,x}\\
&= f_{x} \circ \mlt_{xxx} \circ H_{x,x}\uni_{x} \circ \atpd^{-1}_{xx}\\
&= f_{x} \circ \atpd^{-1}_{xx}
\end{align*}
In $(*)$ we used \cref{prop:integralCasimir} and in $(**)$ we used \cref{antipodeproperties}.
For the second part, notice that both composites below equal $\uni_{x}$, by definition of $\ol{p}$ and $\ol{q}$:
\begin{equation*}
I\xrightarrow{\uni_x}H_{x,x} \xrightarrow{\ol{p}_{x}}H^*_{x,x}\xrightarrow{1\cas{t}{xx}}H^*_{x,x}\ot H_{x,x}\xrightarrow{1\lcomlt_{xx}}H^*_{x,x}\ot H_{x,x}\ot H_{x,x}\xrightarrow{\cas{\ev}{xx}1}H
\end{equation*}
\begin{equation*}
I \xrightarrow{\uni_x}H_{x,x}\xrightarrow{\ol{q}_{x}}H^*_{x,x}\xrightarrow{\cas{t}{xx}1}H_{x,x}\ot H^*_{x,x}\xrightarrow{\lcomlt_{xx}1}H_{x,x}\ot H_{x,x}\ot H^*_{x,x}\xrightarrow{1\braid}H_{x,x}\ot H^*_{x,x}\ot H_{x,x}\xrightarrow{1\cas{\ev}{xx}}H
\end{equation*}
Therefore by choosing $g_{x} = \cas{\ev}{xx} \circ (\ol{p}_{x}\ot H_{x,x}) \circ (\uni_{x}\ot H_{x,x})$ and $f_{x} = \cas{\ev}{xx} \circ \braid \circ (H_{x,x}\ot \ol{q}_{x}) \circ (H_{x,x}\ot \uni_{x})$, the result follows.
\end{proof}

Integral families of \cref{integralfamilies} can be expressed in any semi-Hopf $\Vv$-category $A$. 
In what follows, we require some extra assumptions in order to define an `integral space' as an object in the enriching category $\Vv$, whose generalized elements are precisely those families. Notice that due to standard conventions, left integrals are constructed using right internal homs and vice versa; since $\Vv$ is braided, that subtlety can be ignored.

\begin{definition}\label{def:int}
Suppose $A$ is a semi-Hopf $\Vv$-category, where $\Vv$ is monoidal closed with all limits. 
The {\em left integral space} of $A$ is the diagonal graph $\int^\ell_A= \left\{\left(\int^\ell_A\right)_z\right\}_{z\in X}$ where each object $\left(\int^\ell_A\right)_z$, denoted henceforth $\int^\ell_{A,z}$, is the limit of a diagram in $\Vv$ as below
\begin{equation}\label{defintegralW}
\cd[@R-1em]{
 & & \int^\ell_{A,z} \ar@{-->}[dl]_{t_{xz}} \ar@{-->}[dr]^-{t_{yz}} & & \\
 \dots & A_{x,z} \ar[dl]_-{\ol \mlt_{wxz}} \ar[dr]^-{\ol{\lcouni_{xy}\ot1}} & & A_{y,z} \ar[dl]_-{\ol \mlt_{xyz}}
\ar[dr]^-{\ol{\lcouni_{yu}\ot1}} & \dots \\ 
[A_{w,x},A_{w,z}] && [A_{x,y},A_{x,z}] && [A_{y,u},A_{y,z}]
}
\end{equation}
The morphisms $\ol {\lcouni_{xy}\ot1}$ and $\ol \mlt_{xyz}$ are the adjuncts of $\lcouni_{xy}\ot1\colon A_{x,y}\ot
A_{x,z}\to A_{x,z}$ and $\mlt_{xyz}\colon A_{x,y}\ot A_{y,z}\to A_{x,z}$ under the tensor-hom adjunction for
right-closure.
The limiting cone under $\int^\ell_{A,z}$ is determined by the dashed maps $t_{yz}:\int^\ell_{A,z} \to A_{y,z}$.

The {\em right integral space} $\int^r_A$ of $A$ is computed similarly using left closure, by taking the limit of the diagram
\begin{equation}\label{eq:rightintegral}
\begin{tikzcd}
&& \int^r_{A,x}\ar[dr,dashed]\ar[dl,dashed] && \\
\ldots & A_{x,y}\ar[dr,"{\ol{1\otimes\lcouni_{zy}}}"'] && A_{x,z}\ar[dl,"{\ol{\mlt}_{xzy}}"] & \ldots\\
&& {[}A_{z,y},A_{x,y}{]} &&
\end{tikzcd}
\end{equation}
\end{definition}

As one can expect, in the $k$-linear case an integral space $\int^\ell_{A,x}$ is given exactly by the k-linear space of all integral families of 
the form $\cas{t}{yx}$ for arbitrary $y\in X$. Since it is important 
for what follows, let us spell out the exact connection between integral families of \cref{integralfamilies} and integral spaces 
of \cref{def:int} for general $\Vv$-categories.

The key observation that connects the above definition with integral families is that under the tensor-hom adjunction, each commuting square in the limit diagram \eqref{defintegralW} for example corresponds to a commuting
\begin{equation}\label{eq:squareunderadjunction}
\xymatrix{
A_{z,x}\ot\int^\ell_{A,y} \ar[rr]^-{1\ot t_{xy}} \ar[d]_{1\ot t_{zy}} && A_{z,x}\ot A_{x,y} \ar[d]^-{\mlt_{zxy}}\\
A_{z,x}\ot A_{z,y}\ar[rr]^-{{\lcouni_{zx}\ot1}} && A_{z,y}
}
\end{equation}
Then any identity-on-object diagonal $\Vv$-graph morphism $u\colon \Ii\to \int^\ell_A$ in $\Vv$ consists of morphisms $u_{z}\colon I\to\int^\ell_{A,z}$ which in turn correspond to maps $\cas{t}{xz}$ as below
\begin{equation}\label{universalintegral}
\begin{tikzcd}[row sep=.2in]
&& I\ar[d,"u_z"]\ar[ddl,bend right,"\cas{t}{xz}"']\ar[ddr,bend left,"\cas{t}{yz}"] && \\
&& \int^\ell_{A,z}\ar[dr,dashed,"t_{yz}"]\ar[dl,dashed,"t_{xz}"'] && \\
\ldots & A_{x,z}\ar[dr,"{\ol{\lcouni_{xy}\ot1}}"'] && A_{y,z}\ar[dl,"{\ol{\mlt}_{xyz}}"] & \ldots\\
&& {[}A_{x,y},A_{x,z}{]} &&
\end{tikzcd}
\end{equation}
which all together form an $X^2$-family $\{\cas{t}{xy}\colon I\to A_{x,y}\}$ satisfying precisely \cref{integraldiag}.

In the one-object case, in any monoidal closed category $\Vv$ with limits, both spaces reduce to equalizers due to the shape of the limiting diagrams: the left integral space of a bimonoid $A$ is given by 
\begin{equation}\label{eq:intspacebimonoid}
\begin{tikzcd}
\int^\ell_{A}\ar[r] & A\ar[r,shift left,"{\ol{\lcouni\ot1}}"]\ar[r,shift right,"{\ol{\mlt}}"'] & {[}{A,A}{]}
\end{tikzcd}
\end{equation}
and similarly the right integral space is the equalizer of the adjuncts of the multiplication $m\colon A\ot A\to A$ and $1\ot\lcouni\colon A\ot A\to A$.

On the other hand, in the many-object setting again, one can wonder how the definition of integral space dualizes to a semi-Hopf $\Vv$-opcategory  $(C,\comlt,\couni,\lmlt,\luni)$. The construction of a similar limit will now use the (global) counit $\couni$ and the (local) multiplication $\lmlt$, and 
the switch between these local and global structures makes the limit in this case into a sheer equalizer. Hence we obtain the following definition.

\begin{definition}
For a semi-Hopf $\Vv$-opcategory  $(C,\comlt,\couni,\lmlt,\luni)$, the left integral space of $C$ is the diagonal graph which consists of the equalizers
\begin{equation}\label{eq:equalizerintegral}
\begin{tikzcd}
  \int^\ell_{C,z}\ar[r] & C_{z,z}\ar[r,shift left,"{\ol{\couni_z\ot1}}"]\ar[r,shift right,"{\ol{\lmlt_{zz}}}"'] &
{[}{C_{z,z}},{C_{z,z}}{]}
\end{tikzcd}
\end{equation}
for all $z\in X$.
Symmetrically, we will denote by $\int^r_C$ the right integral space of $C$.
\end{definition}

\cref{leftrightintegral} draws a correspondence between left and right integral families for Hopf categories; this is established as an isomorphism between the integral spaces below.

\begin{proposition}\label{integralspacesareisomorphic}
If $H$ is a Hopf $\Vv$-category with invertible antipode, then $\int^r_H \cong \int^\ell_H$.
\end{proposition}

\begin{proof}
It suffices to show that there is a natural isomorphism between the diagrams \cref{defintegralW,eq:rightintegral} over which the limits are computed, which is easier to see if we translate them under the tensor-hom adjunction (since the one uses right and the other left closure):
$$
\begin{tikzcd}[column sep=.2in]
A_{x,y}\ot A_{x,z}\ar[dr,"\lcouni_{xy}\ot1"'] && A_{x,y}\ot A_{y,z}\ar[dl,"\mlt_{xyz}"] & A_{z,x}\ot A_{y,x}\ar[dr,"1\ot\lcouni_{yx}"'] && A_{z,y}\ot A_{y,x}\ar[dl,"\mlt_{zyx}"] \\
& A_{x,z} &&& A_{z,x}&
\end{tikzcd}
$$
The following commutative squares give isomorphisms between the left and right legs of the above diagrams respectively, where the bottom isomorphism between the common targets is the same:
$$
\begin{tikzcd}[column sep=.3in]
A_{x,y}A_{x,z}\ar[d,"\lcouni_{xy}1"']\ar[r,"\braid"] & A_{x,z}A_{x,y}\ar[r,"\atpd_{xz}\atpd_{xy}"] & A_{z,x}A_{y,x}\ar[d,"1\lcouni_{yx}"] & A_{x,y}A_{y,z}\ar[d,"\mlt_{xyz}"']\ar[r,"\braid"] & A_{y,z}A_{x,y}\ar[r,"\atpd_{yz}\atpd_{xy}"] & A_{z,y}A_{y,x}\ar[d,"\mlt_{zyx}"] \\
A_{x,z}\ar[rr,"\atpd_{xz}"'] && A_{z,x} & A_{x,z}\ar[rr,"\atpd_{xz}"'] && A_{z,x} &
\end{tikzcd}
$$
These are verified using the standard properties relating antipodes with counits and multiplications from \cref{antipodeproperties}, and the 
braidings imply the passage between left and right closure in braided monoidal categories.
\end{proof}

 The following result relates the integral spaces to the coinvariant spaces as in \cref{eq:coinvariants,defconvariants} of specific regular Hopf (op)modules, namely $H_1$, the Hopf $H^\sop$-opmodule $H$ described in \cref{hopfmodulestructuredual}(\ref{it:H1opmodule}) and $H^*_1$, the Hopf $H$-module $H^*$ described in  \cref{hopfmodulestructuredual}(\ref{it:H*1module}). 
As these modules only exist if $H$ is locally rigid, we restrict to this setting now.
In the $k$-linear case, the second part below was shown in \cite[Prop.~10.5]{BCV}.

\begin{proposition}\label{rk_coinv_integrals}
If $H$ is a locally rigid Hopf $\Vv$-category then
$$\int^r_{H,x} \cong (H_1)^{\co H^\sop}_x
\quad\text{and}\quad
\int^r_{H^\sop,x} \cong (H^*_1)^{\co H}_x
$$
\end{proposition}

\begin{proof}
Since $H$ is locally rigid, the (left) internal hom $[H_{y,z},H_{x,z}]$ in $\Vv$ is given, up to isomorphism, by
$H_{x,z} \ot H_{y,z}^* $; in that case, the maps $\ol m_{xzy}$ and $\ol{1 \ot \lcouni_{zy}}$ of \cref{eq:rightintegral}
are precisely the global $H^\sop$-coaction $\gcoact_{xyz}$ from \cref{AisarightAstarop-opmodule} and $1\otimes \lcouni^*_{zy}$.
Thus, regardless of the exact specification of closure, there is a natural isomorphism between the diagrams over which each limit is computed and so the integral space and coinvariant space are themselves isomorphic.
It is worth mentioning that the coinvariant of $H_1$ is computed using only its coaction and therefore the antipode does
not play any role in this part.

Similarly for the second part, recall the Hopf $\Vv$-opcategory structure of $H^\sop$
given in \cref{AstarHopfversions}: the right integral space is given by the equalizer \cref{eq:equalizerintegral},
whereas the equalizer \cref{eq:coinvariants} gives the coinvariant space of the Hopf $H$-module $H^*$. It remains to
compare the local coaction \cref{AstarisarightA-module} to the adjunct of the induced local multiplication
$(\lcomlt^*_{xx}\circ\phi)$ for $H^\sop$, and also $1\ot\uni_x$ to the adjunct of $1\ot\uni_x^*$, the induced
global counit which end up being isomorphic.
\end{proof}

\begin{remark}\label{leftHsopmodule}
The previous result has of course also a version for \emph{left} integrals. For example, one can see that the space of left integrals in $H$ is isomorphic to the space of coinvariants for the left Hopf $H^\sop$-opmodule structure on the $\Vv$-graph $H$ with following global coaction and local action
$$
H_{x,y} \xrightarrow{\coev_{zx}1}H_{z,x}\ot H^*_{z,x} \ot H_{x,y}\xrightarrow{\braid1}H_{z,x}\ot H^*_{z,x} \ot H_{x,y} \xrightarrow{1\mlt_{zxy}}  
H^*_{z,x}\ot H_{z,y}
$$
$$
H^*_{y,x}\ot H_{x,y} \xrightarrow{1\lcomlt_{xy}} H^*_{y,x}\ot H_{x,y}\ot H_{x,y} \xrightarrow{1\braid} H^*_{y,x}\ot H_{x,y}\ot H_{x,y} \xrightarrow{1\atpd_{xy}1}H^*_{y,x}\ot H_{y,x}\ot H_{x,y} \xrightarrow{\cas{\ev}{yx}1} H_{x,y}$$
similarly to the right module structure described in \cref{hopfmodulestructuredual}(\ref{it:H1opmodule}).
\end{remark}

\subsection{Main theorems}

Let $H$ be a Hopf $\Vv$-category. In case the underlying $\Vv$-category is Frobenius, we just say that $H$ is Frobenius, or that $H$ is a {\em Frobenius Hopf} $\Vv$-category. The next results proves the `uniqueness of integrals' for Frobenius semi-Hopf categories.

\begin{proposition}\label{integralresults}
If $A$ is a Frobenius semi-Hopf $\Vv$-category, then its integrals are non-trivial and unique in the sense that $\int^\ell_{A,x}\cong I$ for all $x\in X$.
\end{proposition}

\begin{proof}
If $A$ is Frobenius, by \cref{FrobeniusCasimir} it comes equipped with a Frobenius system
$(e,\nu)$ whose Casimir family gives rise to an integral family $t_e$ denoted $\{\cas{t_e}{xy}\}$ by \cref{prop:integralCasimir}. As explained by \cref{universalintegral},
the integral family $t_e$ is in bijection with a unique family of morphisms $u_y:I\to \int^\ell_{A,y}$ such that $t_{xy}\circ u_y=\cas{t_e}{xy}$.
We will show that each such $u_x$ is an isomorphism, with inverse $\nu_x\circ t_{xx}$ where $\nu$ are the Frobenius functionals.

The following diagram establishes that $\nu_x\circ t_{xx}$ is a right inverse of $u_x$:
\begin{displaymath}
\begin{tikzcd}[column sep=.7in,row sep=.4in]
& \int^\ell_{A,x}\ar[r,"t_{xx }"] & A_{x,x}\ar[dr,"\nu_x"] & \\
I \ar[ur,"u_x"]\ar[urr,"\cas{t_e}{xx}"description]\ar[urr,phantom,bend right=9,"\scriptscriptstyle\cref{t(e)}"]
\ar[r,"\cas{e}{xx}"description]\ar[rr,dotted,bend right=15,"\uni_x"description]
\ar[rr,phantom,bend right=9,"\scriptscriptstyle\cref{casimir}"description]
\ar[rrr,phantom,bend right=12,near end,"\scriptscriptstyle\cref{Hax1}"description]
\ar[rrr,dotted,bend right=20,"1"description]
& A_{x,x}\ot A_{x,x}\ar[ur,"1\ot\lcouni_{xx}"']\ar[r,"\nu_x\ot1"description]
& A_{x,x}\ar[r,"\lcouni_{xx}"'] & I
\end{tikzcd}
\end{displaymath}

The following computation shows that $\nu_x\circ t_{xx}$ is a left inverse of $u_x$: for all $y\in X$, the following diagram commutes
\begin{displaymath}
 \begin{tikzcd}[row sep=.3in]
  && I\ar[r,"u_x"]\ar[dr,"\cas{e}{yx}"']\ar[drr,"\cas{t_e}{yx}"description] & \int^\ell_{A,x}\ar[r,"t_{yx}"]\ar[dr,phantom,"{\scriptstyle\cref{universalintegral}}"] & A_{yx}\ar[d,equal] \\  
  & A_{xx}\ar[ur,"\nu_x"]\ar[r,"\cas{e}{yx}1"] & A_{yx}\ot A_{xy}\ot A_{xx}\ar[rd,"1\lcouni_{xy}1"]\ar[r,"11\nu_x"]\ar[d,phantom,"{\scriptstyle\cref{eq:squareunderadjunction}}"] & \stackrel{\cref{t(e)}}{A_{yx}\ot A_{xy}}\ar[r,"1\lcouni_{xy}"'] & A_{yx}\ar[d,equal] \\
  \int^\ell_{A,x}\ar[ddd,equal]\ar[dr,"t_{yx}"description]\ar[ur,"t_{xx}"]\ar[r,"\cas{e}{yx}1"] & A_{yx}\ot A_{xy}\ot\int^\ell_{A,x}\ar[ur,"11 t_{xx}"]\ar[r,"11t_{yx}"] & A_{yx}\ot A_{xy}\ot A_{yx}\ar[dd,phantom,"{\scriptstyle\cref{casimirdiag}}"]\ar[r,"1\mlt_{xyx}"] &  A_{yx}\ot A_{xx}\ar[r,"1\nu_x"] & A_{yx}\ar[ddd,equal]\\
  & A_{yx}\ar[ur,"{\cas{e}{yx}1}"description]\ar[dr,"1\cas{e}{xx}"] &  & \\
  & \int^\ell_{A,x}\ot A_{xx}\ot A_{xx}\ar[d,bend right=50, phantom, "{\scriptstyle\cref{casimir}}"]\ar[d,"11\nu_{x}"]\ar[r,"t_{yx}11"] & A_{yx}\ot A_{xx}\ot A_{xx}\ar[uur,"\mlt_{yxx}1"description] \\
  \int^\ell_{A,x}\ar[ur,"1\cas{e}{xx}"]\ar[r,"1\uni_x"]\ar[dr,"t_{yx}"']& \int^\ell_{A,x}\ot A_{xx}\ar[r,"t_{yx}1"] & A_{yx}\ot A_{xx} \ar[rr,"\mlt_{yxx}"] && A_{yx} \\
  & A_{yx}\ar[ur,"1\uni_x"description]\ar[urrr,"1"description]
  \end{tikzcd}
\end{displaymath}
Notice that the family $\{t_{yx}\}_y$ is jointly monic since $\int^\ell_{A,x}$ is defined as a limit, therefore $u_x\circ \nu_x\circ t_{xx}=1$ and the proof is complete.
\end{proof}

We refer to the above theorem as ``uniqueness of integrals'', since it shows that for a Frobenius semi-Hopf $\Vv$-category, 
two integral families differ only up to automorphisms of the monoidal unit $I$.
In the one-object case, this implies the folklore result that the integral space \cref{eq:intspacebimonoid} of a Frobenius Hopf monoid in $\Vv$ is isomorphic to the monoidal unit.
In the $k$-linear case for a field $k$, the automorphisms $k\to k$ are just scalars; we recover the classical uniqueness of integrals of Frobenius Hopf algebras up to a scalar, which is part of the classical Larson-Sweedler theorem as discussed in \cref{classicalLarsonSweedler}.

The next result shows that for a Hopf category, being Frobenius and being \emph{locally Frobenius} (namely all local comonoids $H_{x,y}$ are Frobenius in $\Vv$) are two equivalent properties. In particular, this generalizes the classical result that the underlying algebra of a Hopf algebra $H$ is Frobenius if and only if the underlying coalgebra of $H$ is Frobenius.

\begin{theorem}\label{generalLarsonSweedler}
Suppose $H$ is a Hopf $\Vv$-category. The following are equivalent:
\begin{enumerate}[(i)]
\item $H$ is Frobenius (i.e. $H$ is a Frobenius as a $\Vv$-category).
\item $H$ is locally rigid and $\int^\ell_{H,x}\cong I$ for all $x\in X$.
\item \label{it:locFrob} $H$ is locally Frobenius (i.e. all comonoids $H_{x,y}$ are Frobenius in $\Vv$).
\item $H$ is locally rigid and $\int^\ell_{H^*,x}\cong I$ for all $x\in X$.
\end{enumerate}
\end{theorem}

\begin{proof}

$(i)\Rightarrow {(ii)}$. If $H$ is Frobenius, then it is locally rigid by \cref{Frobisrigid} and the left integral space is isomorphic to the monoidal unit by \cref{integralresults}.

$(ii)\Rightarrow {(iii)}$.
The fundamental theorem of Hopf $\Vv$-opcategories (\cref{fundopcat_equiv}) applied to the $H^{\sop}$-Hopf opmodule $H_1$ of \cref{hopfmodulestructuredual}(\ref{it:H1opmodule}) yields an isomorphism of $H^{\sop}$-opmodules 
\begin{equation*}
(H_{1})^{\co H^{\sop}} \ot H^{\sop} \cong H_{1}
\end{equation*}
which using \cref{rk_coinv_integrals} results in 
\begin{equation*}
\int^r_H \ot H^{\sop} \cong H_{1}.
\end{equation*}
Since $H$ is locally rigid, the antipode is invertible by \cref{antipodeisinvertible} so we can apply \cref{integralspacesareisomorphic} to get $\int^r_H\cong \int^\ell_H\cong I$ and hence
we can conclude that $H^{\sop} \cong H_{1}$ as Hopf $H^{\sop}$-opmodules. This implies that for all $x,y \in X$
$$H^{*}_{y,x} \cong (H_{1})_{x,y}=H_{x,y}$$
as right $H^{*}_{y,x}$-modules, where we regard $H^{*}_{y,x}$ as a (local) monoid in $\Vv$. 
It is easy to check that $\atpd_{xy}\colon(H_{1})_{x,y} \to (H_{2})_{x,y}$ for $H_2$ of \cref{hopfmodulestructuredual}(\ref{it:H2opmodule}) is a right $H^{*}_{y,x}$-module (iso)morphism and combining these, we obtain a right $H^{*}_{y,x}$-module isomorphism $H^{*}_{y,x}\cong (H_{2})_{x,y}=H_{y,x}$.
Now using \cref{frobmodules} for the $1$-object case (i.e. any monoid in $\Vv$) we find that every $H^*_{y,x}$ is a Frobenius monoid. It is well-known for any Frobenius monoid in $\Vv$ that its dual is also Frobenius, hence this proves $(iii)$.

Although this is in principle superfluous, let us also prove how $(ii)$ implies $(i)$. Above, we already showed that $H^{\sop} \cong H_{1}$ as Hopf 
$H^{\sop}$-opmodules which means exactly that the $\Vv$-opcategory $H^{\sop}$ is Frobenius by the dual statement of \cref{frobmodules}.  Hence, it 
follows from \cref{dualfrob} that $H$ is a Frobenius $\Vv$-category, or equivalently $H \cong H^\dagger$.

$(iii)\Rightarrow (iv)$. If each comonoid $H_{x,y}$ is Frobenius in $\Vv$, then it is dualizable and $H^*_{x,y}$ is a Frobenius monoid 
(\cref{dualfrob} in the one-object case). As a result of the 1-object case of the direction $(i)\Rightarrow(ii)$ earlier, it is ensured that each 
$\int^\ell_{H^*,x} \cong I$.

$(iv)\Rightarrow (i)$.
If $H$ is locally rigid then the antipode is invertible by \cref{antipodeisinvertible} and therefore
\begin{displaymath}
 (H^*_1)^{\co H} \cong \int^r_{H^{*,\op}} \cong \int^\ell_{H^{*,\op}} \cong I
\end{displaymath}
by \cref{rk_coinv_integrals,integralspacesareisomorphic}.
Then by the fundamental theorem of Hopf modules, \cref{fundhopfmod}, $H^*_1 \cong (H^*_1)^{\co H} \otimes H \cong H$ as Hopf $H$-modules where 
$H^*_1$ is as in \cref{hopfmodulestructuredual}, (\ref{it:H*1module}). In particular, they are isomorphic as plain $H$-modules. Since $s^*\colon 
H^*_2\to H^*_1$ is an $H$-module isomorphism, $H\cong H_2^*$ of \cref{hopfmodulestructuredual}, (\ref{it:Astar2}) which by \cref{frobmodules} implies 
that $H$ is a Frobenius $\Vv$-category.
\end{proof}

\begin{remark}\label{remarkFrobstructuresLS}
Let $H$ be a locally rigid Hopf $\Vv$-category with a `unique' right integral family $t=\{\cas{t}{xy}\}$. Using the explicit formula for the 
isomorphism in the fundamental theorem for Hopf modules, \cref{fundhopfmod}, one can also obtain  an explicit formula for the Frobenius isomorphisms, 
following the proof $(ii) \Rightarrow (i)$. In particular, we find that the Frobenius isomorphism $\phi:H\cong H^\sop$ as right $H$-modules is given 
by 
\begin{displaymath}
\phi_{xy}\colon H^{*}_{y,x}\xrightarrow{1\cas{t}{xy}}H^{*}_{y,x} \ot H_{x,y} \xrightarrow{1\lcomlt_{xy}}H^{*}_{y,x} \ot H_{x,y}  \ot H_{x,y} \xrightarrow{1\atpd_{xy}1}H^{*}_{y,x} \ot H_{y,x} \ot H_{x,y} \xrightarrow{\cas{\ev}{yx}1}H_{x,y}
\end{displaymath}
In the k-linear case, this gives the formula $\varphi_{xy}(f) = f(\atpd_{xy}(\cas{t}{xy}_{(1)}))\cas{t}{xy}_{(2)}$. 

The corresponding Casimir element is exactly Casimir element associated to the right integral $t$ as in \eqref{e(t)right}. To complete this to 
a full Frobenius system one also needs a family of trace morphisms which are then given by
$$\nu_x: \xymatrix{H_{x,x} \ar[r]^-{\uni_x\ot 1} & H_{x,x}\ot H_{x,x} \ar[r]^-{\phi_{xx}^{-1}\ot 1} & H^*_{x,x}\ot H_{x,x} \ar[r]^-{\ev} & I }$$
Since the explicit form of $\phi_{xx}^{-1}$ depends on the explicit form of the isomorphism $\int_H^r\cong I$, we unfortunately don't have a more explicit form for the trace morphisms.

Using symmetric arguments based on the left Hopf $H^\sop$-opmodule from \cref{leftHsopmodule}, we find an (in general different) Frobenius structure on $H$, where the Casimir element in this case is the one from \eqref{e(t)}, induced by a left integral. If we denote the left integral family as $\cas{u}{xy}$ then we obtain this way a Frobenius isomorphism of left $H$-modules of the form
\begin{equation}\label{eq:FrobIsoFromHopfLeft}
\phi_{xy}\colon H^{*}_{y,x}\xrightarrow{1\cas{u}{xy}}H^{*}_{y,x} \ot H_{x,y} \xrightarrow{1(\braid\circ\lcomlt_{xy})}H^{*}_{y,x} \ot H_{x,y}  \ot H_{x,y} \xrightarrow{1\atpd_{xy}1}H^{*}_{y,x} \ot H_{y,x} \ot H_{x,y} \xrightarrow{\cas{\ev}{yx}1}H_{x,y}
\end{equation}

In a similar way, starting again with a unique right integral family $\cas{t}{xy}$ and following the proof $(ii)\Rightarrow(iii)$ above, we can also find an explicit formula form the local Frobenius isomorphisms. In this case, these come out as
\begin{equation}\label{LocalFrobIsoFromHopf}
\psi_{yx}
\colon H^{*}_{y,x}\xrightarrow{1\cas{t}{xy}}H^{*}_{y,x} \ot H_{x,y} \xrightarrow{1\lcomlt_{xy}}H^{*}_{y,x} \ot H_{x,y}  \ot H_{x,y} \xrightarrow{1\atpd_{xy}\atpd_{xy}}H^{*}_{y,x} \ot H_{y,x} \ot H_{y,x} \xrightarrow{\cas{\ev}{yx}1}H_{y,x}
\end{equation}
or in the $k$-linear case $\psi_{y,x}(f)=f(S(\cas{t}{xy}_{(1)}))S(\cas{t}{xy}_{(2)})$, which can explictly be checked to be $H^*_{xy}$-linear. This 
formula does not give us a Casimir element in $H^{*}_{y,x}\ot H^*_{y,x}$ for this structure (for this one would now need the explicit form 
$\psi_{yx}^{-1}$), but we do obtain an explicit form for the (non-degenerate) trace morphism:
$$
\nu_{yx}:\xymatrix{H^*_{y,x}\ar[r]^-{1\cas{t}{xy}} & H^{*}_{y,x} \ot H_{x,y} \ar[r]^-{1\ot \atpd_{xy}} & H^{*}_{y,x} \ot H_{y,x} \ar[r]^-\ev & I
}
$$

Finally, let us remark that the Frobenius structures that one obtains from the proof $(iv) \Rightarrow (i)$ could be different, since it makes us of the fundamental theorem for Hopf $H$-modules and integrals in $H^\sop$, while the previous one comes from the fundamental theorem for $H^{\sop}$-Hopf opmodules and integral in $H$. The previous theorem should therefore not be understood as a theorem stating a bijective correspondence between structures, but it is only an equivalence on existence of certain structures!
\end{remark}

The above \cref{generalLarsonSweedler} gives quite a lot of information regarding Hopf Frobenius $\Vv$-categories. We now gather it all together, 
connecting the different structures on the appropriate level and concluding to \cref{hopffrobconnection}.

Suppose $H$ is a locally rigid Frobenius Hopf $\Vv$-category. Apart from its $\Vv$-category and local comonoid structure as a Hopf category, 
$H$ also has a $\Vv$-opcategory structure from being Frobenius, as well as a local monoid structure by \cref{generalLarsonSweedler}, 
(\ref{it:locFrob}). As was remarked earlier, several Frobenius structures might exist on the same Hopf 
$\Vv$-category: our next aim is to show that when constructed properly, these four (category, opcategory, local algebra, local 
coalgebra) structures on $H$ can be combined in different ways to constitute Hopf and Frobenius structures.

\begin{theorem}\label{twoopcategoriesarethesame}
Let $H$ be a locally rigid Frobenius Hopf $\Vv$-category. Then the Frobenius structure induces a $\Vv$-opcategory structure on $H$, 
such that $H\cong H^\sop$ as Hopf 
$\Vv$-categories via \eqref{eq:FrobIsoFromHopfLeft}.
\end{theorem}

\begin{proof}
Let us denote by $t=\{\cas{t}{xy}\}$ a (unique) right integral family in $H$. From the proof of \cref{generalLarsonSweedler} and 
\cref{remarkFrobstructuresLS}, we know that $H$ is locally Frobenius with Frobenius isomorphisms given by \eqref{LocalFrobIsoFromHopf}. 
If we denote by $u=\atpd\circ t=\{\atpd_{yx}\circ\cas{t}{xy}\}$ the left integral family obtained from $t$ by \cref{leftrightintegral}, we also 
know that it endows $H$ with a Frobenius structure, such that the Frobenius isomorphism is of the form 
\eqref{eq:FrobIsoFromHopfLeft}.

Therefore, in particular, $H$ is a $\Vv$-opcategory (not necessarily with the initial opcategory structure!).
By \cref{Frobisoisopcatiso} we immediately obtain that the family of morphisms $\phi_{xy}\colon H^\sop_{x,y}\to H_{x,y}$ form an isomorphism of 
$\Vv$-opcategories for the considered structures. 

Moreover $H$ has a local monoid structure which is the opposite (using the braiding) of the one induced by the right integral $t$. From the 
one-object 
dual of the same result, \cref{Frobisoisopcatiso}, we know that the morphisms $\psi_{xy}$ from \eqref{LocalFrobIsoFromHopf} are monoid morphisms. Now 
one can easily observe that 
$$\phi_{x,y}=\psi_{y,x}\circ \ol\atpd^*_{x,y}$$
where $\ol\atpd$ is the inverse op-antipode as in \cref{atpdbijopatpd}.
Since $\ol\atpd^*$ is a local anti-monoid automorphism of $H^\sop$, we find that the morphisms $\phi_{x,y}$ are also anti-monoid morphisms, and hence 
become monoid morphisms when we consider the opposite monoid structures on $H_{x,y}$.

We can therefore conclude that the isomorphisms $\phi_{xy}$ are both morphisms of $\Vv$-opcategories and local monoids. Since $H^\sop$ is a Hopf 
$\Vv$-category by \cref{AstarHopfversions}, the given structures on $H$ also satisfy the axioms of a Hopf $\Vv$-opcategory and $\phi$ 
is a Hopf category isomorphism.
\end{proof}

The connection between the four structures on a locally rigid Frobenius Hopf category and their combinations is summarized in the following table, 
which was already observed by Street \cite{streetFAMC} for group-algebras $kG$ over a field $k$.

\begin{center}
\begin{table}[!htbp]
\caption{}
\begin{tabular}{|c|c|c|}
\hline
 & Hopf category  $A$ & Hopf opcategory  $A\cong A^{\sop}$\\
 \hline
 Frobenius category  $A$ & $\mlt_{xyz}$ & $\comlt_{xyz}$ \\
 \hline
 Local Frobenius $A$ & $\lcomlt_{xy}$ & $\lmlt_{xy}$ \\
 \hline
\end{tabular}
\label{hopffrobconnection}
\end{table}
\end{center}

We now investigate, starting from a semi-Hopf rather than Hopf category, how non-singularity of integrals can give rise to antipodes and as a result also Frobenius structures on a $\Vv$-category. Recall the definitions of a right/left antipode and op-antipode, \cref{def:Hopfcat,def:opantipode}, as well as non-singularity of integrals, \cref{def_non_singular}.
\begin{theorem}\label{oppositeLarsonSweedler}
Suppose $(A,\mlt,\uni,\lcomlt,\lcouni)$ is a locally rigid semi-Hopf $\Vv$-category. 
\begin{enumerate}[(i)] 
\item If it has a right non-singular left integral family, then it has right antipode;
if it also has a left non-singular right integral family, then it is Hopf.
\item If it has a left non-singular left integral family, then it has a left op-antipode; if it also has has a right non-singular right integral 
family, then it has an op-antipode.
\item If it has a non-singular left and right integral family, then its antipode is invertible.
\item If it has a non-singular left and right integral family,
then it is Frobenius. \label{itemfour}
\item If $A$ is Hopf with a right integral family $t=\{\cas{t}{xy}\}$ such that it is Frobenius via the induced Casimir family 
$e_t=\{(1\ot\atpd_{xy})\circ\lcomlt_{xy}\circ\cas{t}{xy}\}$, then $t$ is non-singular. \label{itemfive}
\end{enumerate}
\end{theorem}

\begin{proof}\hfill

(i) Denoting the right inverse of $q_{xx}$  \cref{eq:qxy} by $\ol{q}_{x}$, define a composite
\begin{equation}\label{eq:fx}
f_x:A_{x,x}\xrightarrow{1\ot\uni_x}A_{x,x}\ot A_{x,x}\xrightarrow{1\ot \ol q_x}A_{x,x}\ot A^*_{x,x}\xrightarrow{\braid}A^*_{x,x}\ot 
A_{x,x}\xrightarrow{\cas{\ev}{xx}}I
\end{equation}
We then define:
\begin{equation*}
\begin{tikzcd}
\atpd_{xy}\colon A_{x,y}\arrow[rrrdd,dashed]\arrow[r,"1\ot\cas{t}{yx}"]&A_{x,y}\ot A_{y,x}\arrow[r,"1\ot\lcomlt_{yx}"]&
A_{x,y}\ot A_{y,x}\ot A_{y,x}\arrow[r,"\braid\ot1"]& A_{y,x}\ot A_{x,y}\ot A_{y,x}\arrow[d,"1\ot\mlt_{xyx}"]\\
&&& A_{y,x}\ot A_{x,x}\arrow[d,"1\ot f_x"]\\
&&& A_{y,x}
\end{tikzcd}
\end{equation*}
In the $k$-linear case, this is defined by the formula $\atpd_{xy}(a)=\cas{t}{yx}_{(1)}\cdot f_x(a\cas{t}{yx}_{(2)})$.
One now checks easily that the top of \cref{HopfCatAntipodeEquations} follows directly from the left integral condition \cref{integraldiag} so that 
$\atpd_{xy}$ form a right antipode.

Now given a left non-singular right integral family, let $\ol p_y$ be the right inverse of $p_{yy}$ \cref{eq:pxy} and put $g_y=\ev\circ (\ol p_y \ot 
1)\circ(\uni_y \ot 1)\colon A_{y,y}\to I$. We then define a left antipode by
\begin{displaymath}
\begin{tikzcd}
\atpd'_{xy} \colon A_{x,y} \arrow[rrrdd,dashed]\arrow[r,"{\cas{t}{yx}\ot1}"] & A_{y,x}\ot A_{x,y}\arrow[r,"\delta_{yx}\ot1"]& A_{y,x}\ot A_{y,x}\ot A_{x,y} \arrow[r,"1\braid"] & A_{y,x}\ot A_{x,y}\ot A_{y,x}\arrow[d,"\mlt_{xyx}1"]\\
&&& A_{y,y}\ot A_{y,x}\arrow[d,"g_y\ot1"]\\
&&& A_{y,x} 
\end{tikzcd}
\end{displaymath}
In any semi-Hopf $\Vv$-category, if both left and right antipodes exist then they are equal and thus an antipode.

\noindent
(ii) This is essentially dual to the previous part. 
In the first case define 
\begin{displaymath}
\atpd_{xy} \colon A_{x,y}\xrightarrow{1\ot\cas{t}{xy}}A_{y,x}\ot A_{x,y}\xrightarrow{1\ot\lcomlt_{xy}}
A_{y,x}\ot A_{x,y}\ot A_{x,y}\xrightarrow{\mlt_{yxy}\ot1}
A_{y,y}\ot A_{x,x}\xrightarrow{g_y\ot1}A_{y,x}
\end{displaymath}
with $g_y$ as above; and in the second case define
\begin{displaymath}
\atpd'_{xy} \colon A_{x,y}\xrightarrow{\cas{t}{yx}\ot1}A_{y,x}\ot A_{x,y}\xrightarrow{\lcomlt_{yx}\ot1}
A_{y,x}\ot A_{y,x}\ot A_{x,y}\xrightarrow{1\ot\mlt_{yxy}}
A_{y,x}\ot A_{x,x}\xrightarrow{1\ot f_x}A_{y,x}
\end{displaymath}
with $f_x$ as above \cref{eq:fx}.

\noindent
(iii) This is immediate when we combine parts (i) and (ii) with \cref{atpdbijopatpd}.

\noindent
(iv) By part (iii) the semi-Hopf $\Vv$-category $A$ is Hopf and the antipode is invertible. 
Therefore, as in \cref{e(t)} the integral family $t$ gives rise to a Casimir family $e_t$, given by 
$(1\ot\atpd_{xy})\circ\lcomlt_{xy}\circ\cas{t}{xy}$.
By setting $\nu_x\coloneqq A_{x,x}\xrightarrow{\atpd^{-1}_{xx}} A_{x,x}\xrightarrow{f_x}I$
where $f_x$ is as in \cref{eq:fx}, we can show that $(e_t,\nu)$ is a Frobenius system for $A$.
Indeed, \cref{casimir} is satisfied as follows
\begin{align*}
A_{x,x}\nu_{x} \circ \cas{e}{xx} &= A_{x,x}\cas{\ev}{xx}\circ A_{x,x}\braid \circ A_{x,x}A_{x,x}\ol{q}_{x}\circ A_{x,x}A_{x,x}\uni_{x} \circ A_{x,x}\atpd^{-1}_{xx} \circ A_{x,x}\atpd_{xx}\circ \lcomlt_{xx} \circ \cas{t}{xx} \\
&= A_{x,x}\cas{\ev}{xx}\circ A_{x,x}\braid \circ \lcomlt_{xx}A^*_{x,x}\circ \cas{t}{xx}A^*_{x,x} \circ \ol{q}_{x} \circ \uni_{x}\\
&= q_{xx}\circ \ol{q}_{xx} \circ \uni_{x}\\
&= \uni_{x}
\end{align*}

\begin{align*}
\nu_{x}A_{x,x} \circ \cas{e}{xx} &= \cas{\ev}{xx}A_{x,x}\circ \braid A_{x,x}\circ A_{x,x}\ol{q}_{xx}A_{x,x}\circ A_{x,x}\uni_{x}A_{x,x} \circ \atpd^{-1}_{xx}A_{x,x} \circ A_{x,x}\atpd_{xx}\circ \lcomlt_{xx} \circ \cas{t}{xx}\\
&\leftstackrel{(*)}{=} \cas{\ev}{xx}A_{x,x}\circ \ol{p}_{x}A_{x,x}A_{x,x}\circ \uni_{x}A_{x,x}A_{x,x} \circ A_{x,x}\atpd_{xx}\circ \lcomlt_{xx} \circ \cas{t}{xx}\\
&= \atpd_{xx}\circ \cas{\ev}{xx}A_{x,x}\circ \ol{p}_{xx}A_{x,x}A_{x,x}\circ \uni_{x}A_{x,x}A_{x,x} \circ  \lcomlt_{xx} \circ \cas{t}{xx}\\
&=  \atpd_{xx}\circ \cas{\ev}{xx}A_{x,x}\circ A^*_{x,x}\lcomlt_{xx} \circ A^*_{x,x}\cas{t}{xx} \circ \ol p_{xx}\circ \uni_x\\
&= \atpd_{xx}\circ p_{xx}\circ \ol p_{xx}\circ  \uni_{x}=\atpd_{xx}\circ \uni_{x}\\
&= \uni_{x}
\end{align*}
where $(*)$ follows from \cref{leftrightnonsingular}.

\noindent
(v) By \cref{frobmodules}, we know that $A$ is locally rigid -- thus the antipode is invertible by \cref{antipodeisinvertible} -- and also we have an 
isomorphism of $A$ and $A^{\sop}$ as left and right $A$-modules. Let us denote $\phi_{xy}\colon A^{\sop}_{x,y} \to A_{x,y}$ and 
$\phi^{'}_{xy}\colon A^{\sop}_{x,y} \to A_{x,y}$ for respectively the right and left $A$-module isomorphisms as in \cref{VW3}, namely $\varphi_{xy} = 
(\cas{\ev}{yx} \ot A_{x,y}) \circ (A^{*}_{y,x} \ot \cas{e}{yx})$ and $\varphi^{'}_{xy} = (A_{x,y} \ot \cas{\ev}{yx}) \circ (A_{x,y} \ot \braid) \circ 
(\cas{e}{xy} \ot A^{*}_{y,x})$. The invertibility of $p_{xx}$ and $q_{xx}$ as in \cref{def_non_singular} now follows from the following 
factorisations, which in fact shows that all arbitrary-indices composites are not only split epimorphisms but isomorphisms. 
\begin{equation*}
\begin{tikzcd}[column sep=.4in,row sep=.4in]
A^{*}_{x,y} \ar[r,"\atpd^{*}_{xy}"] \ar[rr, bend right=15, "\varphi_{yx}"']& A^{*}_{y,x}\ar[r,"p_{yx}"] & A_{y,x}
\end{tikzcd}
\qquad
\begin{tikzcd}[column sep=.4in,row sep=.4in]
A^{*}_{y,x} \ar[r,"q_{yx}"]\ar[rr, bend right=15, "\varphi^{'}_{xy}"']  & A_{y,x} \ar[r,"\atpd_{yx}"] & A_{x,y} 
\end{tikzcd}
\end{equation*}
\end{proof}

We are now ready to formulate and prove the main result of this paper, 
which we call the Larson-Sweedler theorem\index{Larson-Sweedler} for Hopf $\Vv$-categories.

\begin{theorem}\label{oppositeLarsonSweedler2}
Suppose $A$ is a locally rigid semi-Hopf $\Vv$-category. The following are equivalent:
\begin{enumerate}[(i)]
\item $A$ is Hopf and has a non-singular right integral family;
\item $A$ has both a non-singular right integral family, and a non-singular left integral family;
\item $A$ is Hopf and Frobenius;
\item $A$ is Hopf and $\int^\ell_{A,x} \cong I$;
\item any statement dual the those above, for the dual semi-Hopf $\Vv$-opcategory $A^*$.
\item interchanging left and right in statements (i) and (iv).
\end{enumerate}
\end{theorem}

\begin{proof}\hfill

\noindent
$(i)\Rightarrow(ii)$. Since $A$ is locally rigid and Hopf, we know by \cref{antipodeisinvertible} that the antipode is invertible. Hence by \cref{leftrightintegral}, $A$ also admits a non-singular left integral family.

\noindent
$(ii)\Rightarrow(iii)$.
It follows from (i) and (iv) of \cref{oppositeLarsonSweedler}.

\noindent
$(iii)\Rightarrow(iv)$. It follows from \cref{generalLarsonSweedler}.

\noindent
$(iv)\Rightarrow (i)$.
If $A$ has a Hopf $\Vv$-category structure and $\int^{l}_{A,x} \cong I$, we know from \cref{generalLarsonSweedler} that $A$ has a Frobenius structure. Moreover, as explained in \cref{remarkFrobstructuresLS}, one can choose a Frobenius structure such that the Casimir family of this Frobenius structure is exactly the one as stated in \cref{oppositeLarsonSweedler} (v). Hence this last mentioned theorem tells us that there is a non-singular right integral family. 

\noindent
$(v)$. This is obvious by duality. We know that $A$ is Frobenius and Hopf if and only if the opcategory $A^*$ is Frobenius and Hopf, see \cref{dualfrob,twoopcategoriesarethesame}.

\noindent
$(vi)$. Obvious since e.g.\ (ii) is left-right symmetric.
\end{proof}

The above version of the Larson Sweedler might look slightly weaker than the classical theorem. Indeed, for example item (i) assumes both non-singular 
left and right integrals, whereas in the classical theorem one needs only one. Similarly, the Hopf condition is always combined with 
an additional assumption, such as Frobenius or uniqueness of integrals. However, remark as well that the assumptions made in the classical 
Larson-Sweedler theorem are much stronger than the ones imposed herein. Indeed, the classical Larson-Sweedler theorem considers only Hopf algebras 
that are free of finite rank over a principal ideal domain. Our theorem concerns any locally rigid Hopf categories and applies in particular to any 
finitely generated and projective Hopf algebra over an arbitrary commutative ring.

Nevertheless, our last aim is to show that \cref{oppositeLarsonSweedler2} subsumes the classical Larson-Sweedler theorem for Hopf algebras 
when these stronger conditions are imposed. First recall that by working over a PID, every projective module is free and therefore one can use 
dimension arguments when dealing with free modules over a PID. Other (commutative) rings that have this property are local rings (Kaplansky) and 
polynomial rings over a field (Quillen–Suslin). Therefore, for the rest of this section we will consider modules over a 
commutative base ring $k$ for which every finitely generated and projective module is free. As in \cite{LS}, we will say that a $k$-module is {\em finite dimensional} if it 
is projective (hence free) of finite (and constant since $k$ is commutative) rank over $k$. Let us first generalize \cite[Lemma 1]{LS} to the 
multi-object case.

\begin{lemma}\label{Lemma1LS}
Let $A$ be a locally rigid $k$-linear semi-Hopf category, such that for any two objects $x,y$, $A_{x,y}$ and $A_{y,y}$ have the same dimension if $A_{x,y}$ is non-zero. If $\atpd$ is a right antipode of $A$, then it is also a left antipode of $A$.
\end{lemma}

\begin{proof}
Since $\atpd$ is a right antipode we know that $A_{x,y} * \atpd_{xy} = \uni_{x} \circ \lcouni_{xy}$ for every $x,y \in A$, where $*$ denotes the 
convolution product. Define $$\Gamma_{xy}: \Hom(A_{x,y},A_{y,y}) \to \Hom(A_{x,y}, A_{y,x})$$ by $\Gamma_{x,y}(f) = f* \atpd_{xy} $ for every $f: 
A_{x,y} \to A_{y,y}$. This map is clearly surjective, since we can write every $g \in \Hom(A_{x,y}, A_{y,x})$ as $\Gamma_{xy}(g * A_{x,y})$. Using the 
fact that all $A_{x,y}$ and $A_{y,y}$ have the same dimension, we can conclude that every $\Gamma_{xy}$ is bijective.

Clearly the map $\Gamma'_{xy}: \Hom(A_{x,y}, A_{y,x}) \to \Hom(A_{x,y}, A_{y,y}): r \mapsto r * A_{x,y}$ is a right inverse of $\Gamma_{xy}$. And by 
bijectivity of $\Gamma_{xy}$, any one-sided inverse is a two-sided inverse and hence $\Gamma'_{xy}$ is also bijective. This implies that there exists 
a morphism $u_{x,y} \in \Hom(A_{x,y}, A_{y,x})$ such that $\uni_{y} \circ \lcouni_{xy} = u_{x,y}* A_{x,y}$ for every $x,y \in A$.

Moreover $u_{x,y} = u_{x,y} * (A_{x,y} * \atpd_{xy}) =  (u_{x,y} * A_{x,y}) * \atpd_{xy} = \atpd_{xy}$. We can conlude that $\atpd$ is indeed an antipode for $A$.
\end{proof}

\begin{corollary}\label{classicalLScat}
In case $\Vv=\Mod_k$, where $k$ is a ring such that all projective modules are free, the equivalent statements of \cref{oppositeLarsonSweedler2} are furthermore equivalent to
\begin{enumerate}[(i)]
\setcounter{enumi}{6}
\item $A$ has a right non-singular left integral family and for any two objects $x,y$, $A_{x,y}$ and $A_{y,y}$ have the same dimension if $A_{x,y}$ is non-zero. 
\item $A$ is Hopf.
\end{enumerate}
\end{corollary}

\begin{proof}
$\ul{(i)/(vi)\Rightarrow(vii)}$. The left version (i) tells in particular that $A$ has a right non-singular left integral family. Since $A$ is Hopf, we know by \cref{Galoisobjects} that $A_{x,y}\ot A_{y,y}\cong A_{x,y}\ot A_{x,y}$ in $\Mod_k$ and therefore $A_{x,y}$ and $A_{y,y}$ have the same dimension if $A_{x,y}$ is non-zero. 

$\ul{(vii)\Rightarrow(viii)}$. By \cref{oppositeLarsonSweedler}(i), the existence of a right non-singular left integral family implies that $A$ has a right antipode and therefore by \cref{Lemma1LS}, $A$ also has a two-sided antipode.

$\ul{(viii)\Rightarrow(iv)}$. As in the proof of \cref{generalLarsonSweedler}, by the fundamental theorem for Hopf modules we obtain $A_{x,x}^*\cong \int^\ell_{A,x}\ot A_{x,x}$.
Since $A_{x,x}$ and $A_{x,x}^*$ have the same dimension, we find that the dimension of $\int^{\ell}_{A,x}$ must be one, i.e.\ 
$\int^{\ell}_{A,x}$ is free of rank one.
\end{proof}

\section{Examples}
\label{sec:applications}

In this section, we gather a few important examples that are obtained as results of the generalization of the Larson-Sweedler Theorem, and we provide some directions for further research.

\subsubsection*{Hopf algebras in a monoidal category \texorpdfstring{$\Vv$}{V}}

For its one-object case, \cref{generalLarsonSweedler} gives a version of the Larson-Sweedler theorem for Frobenius and Hopf algebras in any braided 
monoidal category $\ca{V}$. 
In particular, if $H$ is a Hopf monoid in a braided monoidal category $\Vv$, then it is Frobenius if and only if it is dualizable and its 
integral 
\cref{eq:intspacebimonoid} is isomorphic to the monoidal unit $I$. 

In particular, by regarding  the 1-object case of \cref{classicalLScat}, we recover the ’classical’ Larson-Sweedler theorem (for Frobenius and Hopf 
$k$-algebras). In the same way, by considering the 1-object case for the monoidal Hom-category $\tilde{\mathcal H}({\mathcal C})$ associated to a 
braided monoidal category $\mathcal C$ as constructed in \cite{stefhomhopf}, we obtain a version of the Larson-Sweedler for monoidal Hom-Hopf 
algebras. In the same way, by choosing suitable braided monoidal categories, one can derive the Larson-Sweedler theorem for graded Hopf algebras and 
Yetter-Drinfel'd Hopf algebras \cite{S2002}.

Let us remark that there is a subtle difference between the classical case, where the considered monoidal category is $\Mod_k$, and the other mentioned cases, where the considered monoidal category is a category of graded vector spaces or Yetter-Drinfel'd modules. As the monoidal unit $k$ of $\Mod_k$ is also a generator in this category, any element $t$ of the Hopf algebra $H$ can be understood as a morphism $t:k\to H$ in this category. This is no longer true for graded Hopf algebras and Yetter-Drinfel'd modules. For example, in the case of graded modules, morphisms from the monoidal unit to $H$ correspond to homogeneous elements of degree $0$ (where $0$ is the unit of the grading group). Consequently, it follows from \cref{oppositeLarsonSweedler2} that a graded Hopf algebra is Frobenius if and only if it is finite dimensional in each degree and has a non-singular integral of degree $0$. Nevertheless, there exist graded Hopf algebras with (non-singular) integrals of arbitrary degree, but these are not Frobenius.

\subsubsection*{Turaev's Hopf group-algebras}
Recall from \cite{Z2004} the definition of a \emph{Hopf $G$-algebra}\index{Hopf $G$-algebra}. Let $G$ be a group. A Hopf $G$-algebra $H$ consist of a $G$-indexed family of $k$-coalgebras $\left(H_{g}, \Delta_{g}, \couni_{g}\right)_{g \in G}$ endowed with the following data.
\begin{enumerate}
\item A family of coalgebra morphisms $\mu = \left(\mu_{g,h}: H_{g} \ot H_{h} \to H_{gh}\right)_{g,h \in G}$, called the multiplication such that
\begin{equation*}
\mu_{gh,l} \circ (\mu_{g,h} \ot H_{l}) = \mu_{g,hl} \circ (H_{h} \ot \mu_{h,l})
\end{equation*}
for every $g,h,l \in G$
\item A coalgebra morphism $\eta: k \to H_{1}$, called the unit such that
\begin{equation*}
\mu_{g,1} \circ (H_{g} \ot \eta) = H_{g} = \mu_{1,g} \circ (\eta \ot H_{g})
\end{equation*}
for every $g \in G$
\item A family of coalgebra isomorphisms $\psi = \left(\psi^{g}_{h}: H_{g} \to H_{hgh^{-1}}\right)$, which need to satisfy:
\begin{equation*}
\psi^{hgh^{-1}}_{l} \circ \psi^{g}_{h}  = \psi_{lh}
\end{equation*}
\begin{equation*}
\psi^{gh}_{l} \circ \mu_{g,h} = \mu_{lgl^{-1},lgl^{-1}} \circ (\psi^{g}_{l}\ot \psi^{h}_{l})
\end{equation*}
\begin{equation*}
\psi^{1}_{h}\circ \eta = \eta
\end{equation*}
for every $g,h,l \in G$
\item A family of maps $\atpd = \left(\atpd_{g}: H_{g} \to H_{g^{-1}}\right)_{g \in G}$
such that
\begin{equation*}
\mu_{g^{-1},g}\circ (\atpd_{g} \ot H_{g}) \circ \Delta_{g} = \mu_{g,g^{-1}} \circ (H_{g} \ot \atpd_{g}) \circ \Delta_{g} = \eta \circ \couni_{g}
\end{equation*}
for every $g \in G$
\end{enumerate} 

Dually one has the notion of a \emph{Hopf $G$-coalgebra}\index{Hopf $G$-coalgebra}. Every Hopf G-(co)algebra can be turned into a k-linear Hopf (op)category. We provide here the construction for a Hopf G-algebra, for a Hopf $G$-algebra $\left((H_{g})_{g \in G}, \mu, \eta, \Delta, \couni\right)$ we define the $k$-linear Hopf category $(\tilde{H}_{x,y})_{x,y \in G} $ by $\tilde{H}_{x,y} := H_{x^{-1}y}$, $\mlt_{xyz}: \tilde{H}_{x,y} \ot \tilde{H}_{y,z} = H_{x^{-1}y} \ot H_{y^{-1}z} \xrightarrow{\mu_{x^{-1}y,y^{-1}z}} H_{x^{-1}z} = \tilde{H}_{x,z}$, $\uni_{x} = \eta$, $\lcomlt_{xy} = \Delta_{x^{-1}y}$ and $\lcouni_{xy} = \couni_{x^{-1}y}$, see \cite[Proposition 6.2]{BCV}.
In case $H$ is a Hopf $G$-algebra such that all $H_g$ are finite dimensional, we can apply the Larson-Sweedler theorem for Hopf categories (see \cref{oppositeLarsonSweedler2}) to the Hopf category associated to this Hopf $G$-algebra and obtain in this way a Frobenius $k$-linear category with $k$-linear morphisms $\comlt_{xyz}: \tilde{H}_{x,y} = H_{x^{-1}y} \to H_{x^{-1}y} \ot H_{x^{-1}y} = \tilde{H}_{x,y} \ot \tilde{H}_{x,y}$ and $\couni_{x}: \tilde{H}_{x,x} = H_{1} \to k$ satisfying conditions \cref{frob1}.

A natural question is whether there exists already a version of the Larson-Sweedler theorem for Hopf $G$-algebras, without using the passage to Hopf categories as described above. A first naive approach would be to use the result from \cite{CD2006}, which states that a Hopf $G$-algebra is an Hopf algebra in a suitably constructed monoidal category of families of $k$-vector spaces called Turaev/Zunino category, and to apply the Larson-Sweedler theorem for Hopf algebras in this monoidal category. However, this will not lead to the desired result, as for this we should require that the Hopf $G$-algebra $(G,H_g)$ is a rigid object in the Zunino category. As this category is equiped with a strict monoidal forgetful functor to $\Set$, sending the indexing group $G$ to its underlying set, the rigidity of $(G,H_g)$ in the Zunino category implies that the set $G$ is a rigid object in $\Set$, which means that it is a singleton, and hence this can only be applied to the classical case of a usual finite dimensional Hopf algebra.
On the other hand, the notion of a Frobenius $G$-algebra already appeared in \cite{T2010}: 
a $G$-algebra $A$ together with a symmetric $k$-bilinear form $\rho: A \ot A \to k$ such that
\begin{enumerate}
\item $\rho(A_{g} \ot A_{h}) = 0$ if $h \neq g^{-1}$
\item The restriction of $\rho$ to  $A_{g} \ot A_{g^{-1}}$ is non-degenerate for every $g \in G$
\item $\rho \circ (\mu_{g,h} \ot A_{l}) = \rho \circ (A_{g}\ot\mu_{h,l})$
\end{enumerate}

In a similar way the construction to obtain a Hopf category out of a Hopf $G$-algebra, one can construct a Frobenius category out of a Frobenius $G$-algebra. If $(A, \rho)$ is a Frobenius $G$-algebra as decribed above, then $\tilde{A}_{x,y}:= A_{x^{-1}y}$ is indeed a $k$-linear Frobenius category. The $k$-linear category structure is obtained in the exact same way as for the Hopf case, $\mlt_{xyz}:= \mu_{x^{-1}y,y^{-1}z}$ and $\uni_{x}  := \eta$. To see it is Frobenius we use the characterization given in \cref{CalabiYauEquiv}. The bilinear form $\Gamma_{xy}: \tilde{A}_{x,y} \ot \tilde{A}_{y,x} \to k$ can be defined as the restriction of $\rho$ to $A_{x^{-1}y} \ot A_{y^{-1}x}$, which is non-degenerate by definition of a Frobenius $G$-algebra.

No other equivalent definitions were given in this reference. Based on our work, one could provide equivalent characterizations for Frobenius $G$-algebras as those described in \cref{sec:FrobeniusVCats} in case of Frobenius categories can be obtained; observe the strong similarity between the definition of a Frobenius $G$-algebra and the definition of a $k$-linear Frobenius algebra as described in \cref{CalabiYauEquiv}. Furthermore, we believe  that a Larson-Sweedler type theorem in this setting can also be obtained in such a way that the following diagram commutes:

\begin{center}
\begin{tikzcd}
\text{Hopf }~ G\text{-algebra} \arrow[rr]\arrow[d,"L-S"']&& \text{Hopf}~\text{ category}\arrow[d,"L-S"]\\
\text{Frobenius }~ G\text{-algebra} \arrow[rr]&& \text{Frobenius}~\text{ category}
\end{tikzcd}
\end{center}

Following this idea, the notion of Hopf $G$-algebra and Hopf category could be unified by means of a more general version of Hopf categories, where the indexing set $X \times X$ is replaced by any groupoid, since both definitions essentially rely on the groupoid structures of $X \times X$ and $G$. Moreover, such further work would also investigate a unified Larson-Sweedler theorem in this setting.

\subsubsection*{Weak (multiplier) Hopf algebras}

In \cite{BCV} (see also \cref{packed}) it is shown that for a $k$-linear Hopf category $A$ with a finite set of objects $X$, $\bigoplus_{x,y\in X} A_{x,y}$ is a weak Hopf algebra. 
If each $A_{xy}$ is in fact finite-dimensional, then
\cref{classicalLScat} in combination with \cref{compactification}
ensures that $\bigoplus_{x,y\in X} A_{x,y}$ is a weak Hopf algebra which is also Frobenius.
This could also be deduced from \cite{iovanov}, since the base of a weak Hopf algebra associated to a Hopf category with a finite number of objects is the cartesian product $k^n$ where $n$ is the finite number of objects in the category.

In case the set of objects $X$ is not finite, the same construction of the `packed' algebra $\oplus_{x,y\in X} A_{x,y}$ will lead to a weak multiplier Hopf algebra, which is Frobenius as an algebra. This can be compared to the Larson-Sweedler theorem for weak multiplier Hopf algebras as proven in \cite{KVD2018}.

\subsubsection*{Groupoid algebra}

Consider a groupoid $G$, a field $k$ and let $G_{x,y}$ be the set of maps from $y$ to $x$. Put $A_{x,y} = kG_{x,y}$. As explained in \cite{BCV}, $A$ has the structure of a $k$-linear Hopf category. We briefly recall the structure:
The multiplication is the one from the groupoid and extended linearly. Every $kG_{x,y}$ has the structure of a coalgebra: $\lcomlt_{xy}(g) = g \ot g$ and $\lcouni_{xy}(g) = 1$. 
The antipode is given by the formula $\atpd_{xy}(g) = g^{-1} \in G_{y,x}$.

If $G$ is locally rigid we know from our Larson-Sweedler theorem that there is a global and local Frobenius structure on it. Let us describe these structures explicitly.

The (global) Frobenius $k$-linear category structure is given by $k$-linear category structure described above and the cocategory structure is given by:
\begin{center}
$\comlt_{xyz}(g) = \sum_{h \in G_{y,z}}gh^{-1} \ot h \in G_{x,y} \ot G_{y,z}$

$\couni_{xy}(g) = 
\begin{cases}
1 ~~;~~ g = e\\
0 ~~;~~ g \neq e
\end{cases}$
\end{center}

The (local) Frobenius structure on every $A_{x,y}$ is given by $\lcomlt_{xy}$ as previously described and the local multiplication by:
\begin{center}
$\lmlt_{xy}: G_{x,y} \ot G_{x,y} \to G_{x,y} : g \ot h \mapsto \begin{cases} g ~~;~~ g = h\\ 0 ~~;~~ g \neq h\end{cases}$

$\luni_{xy}: k \to G_{x,y}: 1 \mapsto \sum_{g \in G_{x,y}}g$
\end{center}
and extended linearly.

\printbibliography

\end{document}